\documentclass[10pt]{article}
 
\usepackage{coursE}
\usepackage{url}

\DeclareMathOperator{\diverg}{div}

\title{Hardy and BMO spaces on graphs, application to Riesz transforms}

\author{{\Large Joseph {\sc Feneuil} \footnote{The author is supported by the ANR project ``Harmonic Analysis at its Boundaries'',   ANR-12-BS01-0013-03.}} 
\\
Institut Fourier, UMR 5582\\
100 rue des Math\'ematiques, BP 74, F-38402 
Saint-Martin d'H\`eres, France\\
joseph.feneuil@ujf-grenoble.fr}
\date{\today}

\newcommand{\Ro}{\mathcal R}

\begin{document}

\maketitle

\begin{abstract}
Let $\Gamma$ be a graph with the doubling property  for the volume of balls  and $P$ a reversible random walk on $\Gamma$. 
We introduce  $H^1$ Hardy spaces of functions and $1$-forms  adapted to $P$ and prove various characterizations of  these spaces. 
We  also characterize the dual space of $H^1$ as a $BMO$-type space  adapted to $P$.  As an application, we establish  $H^1$-$H^1$ and $H^1$-$L^1$ boundedness of the Riesz transform.
\end{abstract}


\tableofcontents

\pagebreak

We use the following notations. $A(x) \lesssim B(x)$ means that there exists $C$ independant of $x$ such that $A(x) \leq C \,  B(x)$ for all $x$,  while  $A(x) \simeq B(x)$ means that $A(x) \lesssim B(x)$ and $B(x) \lesssim A(x)$.
The parameters  from which the constant is independant  will be either obvious from context or recalled. 

\noindent Furthermore, if $E,F$ are Banach spaces, $E \subset F$ means that $E$ is continuously included in $F$. In the same way, $E=F$ means that the norms are equivalent.

\section{Introduction and statement of the results}

The study of real variable Hardy spaces in $\R^n$ began in the early 1960's with the paper of Stein and Weiss \cite{SW}.
At the time, the spaces were defined by means of Riesz transforms and harmonic functions. Fefferman and Stein
provided in \cite{FS2} various characterizations (for instance in terms of suitable maximal functions) and developed real
variable methods for the study of Hardy spaces.

In several issues in harmonic analysis, $H^1(\R^n)$ turns out to be the proper substitute of $L^1(\R^n)$. For example, the
Riesz transforms, namely the operators $R_j = \dr_j (-\Delta)^{-\frac12}$, are $L^p(\R^n)$ bounded for all $p \in(1,+\infty)$,  $H^1(\R^n)$-bounded,
but not $L^1(\R^n)$-bounded  (see \cite{Meyer}).

Hardy spaces were defined in the more general context of spaces of homogeneous type by Coifman and Weiss
in \cite{CW}, by means of an atomic decomposition. An atom is defined as a function supported in a ball, with zero
integral and suitable size condition. However, even in the Euclidean context, the definition of the Hardy space $H^1$
given by Coifman and Weiss is not always suited  to  the $H^1$-$L^1$ boundedness of some Calder\`on-Zygmund type
operators. Indeed, the cancellation condition satisfied by atoms does not always match with differential operators
(consider the case of  $-\diverg(A\nabla)$  on $\R^n$, for instance).

To overcome this difficulty, Hardy spaces adapted to operators were developed in various frameworks during the
last decade. In 2005, in \cite{DY2}  and  \cite{DY1}, Duong and Yan defined Hardy and $BMO$  spaces for an operator  $L$ when the
kernel of the semigroup generated by $L$ satisfies a pointwise Gaussian upper bound. It was discovered later that,
together with the doubling condition for the volumes of balls, $L^2$ Davies-Gaffney type estimates for the semigroup
generated by $L$ are enough to develop a quite rich theory of Hardy spaces on Riemannian manifolds (see \cite{AMR}) and
for second order divergence form elliptic operator in $\R^n$ with measurable complex coefficients (see \cite{HoMay}). These
ideas were pushed further in the general context of doubling measure spaces when $L$ is self-adjoint (see \cite{HLMMY}).

The present work is devoted to an analogous theory of Hardy spaces in a discrete context, namely in graphs $\Gamma$
equipped with a suitable discrete Laplace operator,  given by  $I-P$ where $P$ is a Markov operator (see \cite{grigograph} and the
references therein). We define and give various characterizations of the Hardy space $H^1(\Gamma)$ adapted to $P$, under
very weak assumptions on $\Gamma$. The first characterization is formulated in terms of quadratic functionals (of Lusin
type), relying on results and methods developed in \cite{BRuss} and \cite{Fen1}. The second one is the molecular (or atomic) decomposition
of $H^1(\Gamma)$. A description of the dual space of $H^1(\Gamma)$ as a $BMO$-type space is obtained.

We also deal with the Riesz transform on $\Gamma$, namely the operator $d(I-P )^{-\frac12}$ , where $d$ stands for the differential
on $\Gamma$ ({\em i.e. } $df(x,y) := f(y)-f(x)$ for all functions $f$ on $\Gamma$ and all edges $(x,y)$). When $p \in (1,+\infty)$, the $L^p$-boundedness 
of the Riesz transform was dealt with in \cite{BRuss,Russ}. Here, we prove an endpoint boundedness result for
$p = 1$: roughly speaking, the Riesz transform is $H^1$-bounded. In the same spirit as \cite{AMR}, this assertion requires
the definition a Hardy space of ``exact 1-forms`` on the edges of $\Gamma$. We define and give characterizations of this
space by quadratic functionals and molecular decompositions. Finally, the $H^1$-boundedness of the Riesz
transform is established.

Some Hardy spaces associated with $I-P$ were introduced and characterized in \cite{BD}, together with a description
of their duals and the $H^1$-$L^1$ boundedness of Riesz transform was proved. Even if the authors in \cite{BD} also deal
with the case of $H^p$ for $p < 1$, their assumptions on $P$ are stronger than ours (they assume a pointwise Gaussian
upper bound on the iterates of the kernel of $P$, which is not required for most of our results) and they do not
consider Hardy spaces of forms. Moreover, the Hardy spaces introduced in the present work are bigger than the
ones in \cite{BD}.

\subsection{The discrete setting}

Let $\Gamma$ be an infinite set and $\mu_{xy} = \mu_{yx} \geq 0$ a symmetric weight on $\Gamma \times \Gamma$. 
The couple $(\Gamma, \mu)$ induces a (weighted unoriented) graph structure if we define the set of edges by
$$E = \{ (x,y) \in \Gamma \times \Gamma, \, \mu_{xy} >0 \}.$$
We call then $x$ and $y$ neighbours (or $x\sim y$) if $(x,y) \in E$.\par
\noindent We will assume that the graph is connected and locally uniformly finite. 
A graph is connected if for all $x,y \in \Gamma$, there exists a path $x = x_0,x_1, \dots,x_N = y$ such that for all $1\leq i\leq N$, $x_{i-1}\sim x_i$ (the length of such path is then $N$).
A graph is said to be locally uniformly finite if there exists $M_0 \in \N$ such that for all $x\in \Gamma$, $\# \{y\in \Gamma, \, y\sim x\} \leq M_0$ (i.e. the number of neighbours of a   vertex    is uniformly bounded).\par
\noindent The graph is endowed with its natural metric $d$, which is the shortest length of a path joining two points. 
  For all $x\in \Gamma$ and all $r>0$, the ball of center $x$ and radius $r$ is defined as $B(x,r) = \{y\in \Gamma, \, d(x,y) <r\}$. 
In the opposite way, the radius of a ball $B$ is the only integer $r$ such that $B=B(x_B,r)$ (with $x_B$ the center of $B$). 
Therefore, for all balls   $B=B(x,r)$ and all $\lambda  \geq 1$, we set $\lambda B:=B(x,\lambda r)$ and define    $C_j(B) = 2^{j+1}B \backslash 2^jB$ for all $j\geq 2$ and $C_1(B) = 4B$.\par
\noindent  If $E,F\subset \Gamma$, $d(E,F)$ stands for the distance between $E$ and $F$, namely
$$
d(E,F)=\inf_{x\in E,\ y\in F} d(x,y).
$$
\noindent  We define the weight $m(x)$ of a vertex $x \in \Gamma$ by $m(x) = \sum_{x\sim y} \mu_{xy}$.
More generally, the volume of a subset $E \subset \Gamma$ is defined as $m(E) := \sum_{x\in E} m(x)$. We use the notation $V(x,r)$ for the volume of the ball $B(x,r)$, and in the same way, $V(B)$ represents the volume of a ball $B$. \par
\noindent We define now the $L^p(\Gamma)$ spaces. For all $1\leq p < +\infty$, we say that a function $f$ on $\Gamma$ belongs to $L^p(\Gamma,m)$ (or $L^p(\Gamma)$) if
$$\|f\|_p := \left( \sum_{x\in \Gamma} |f(x)|^p m(x) \right)^{\frac{1}{p}} < +\infty,$$
while $L^\infty(\Gamma)$ is the  space  of functions satisfying 
$$\|f\|_{\infty} : = \sup_{x\in\Gamma} |f(x)| <+\infty.$$
Let us define for all $x,y\in \Gamma$ the discrete-time reversible Markov kernel $p$ associated with the measure $m$ by $p(x,y) = \frac{\mu_{xy}}{m(x)m(y)}$.
The discrete kernel $p_l(x,y)$ is then defined recursively for all $l\geq 0$ by
\begin{equation}
\left\{ 
\begin{array}{l}
p_0(x,y) = \frac{\delta(x,y)}{m(y)} \\
p_{l+1}(x,y) = \sum_{z\in \Gamma} p(x,z)p_l(z,y)m(z).
\end{array}
\right.
\end{equation}

\begin{rmq}
Note that this definition of $p_l$ differs from the one of $p_l$ in \cite{Russ}, \cite{BRuss} or \cite{Delmotte1},   because of the $m(y)$ factor.    However,   $p_l$ coincides with $K_l$    in \cite{Dungey}. 
Remark that in the case of the Cayley graphs of finitely generated discrete groups, where $m(x)=1$ for all $x$, the definitions   coincide.   
\end{rmq}

Notice that for all $l\geq 1$, we have
\begin{equation} \label{sum=1}
 \|p_l(x,.)\|_{L^1(\Gamma)} = \sum_{y\in \Gamma} p_l(x,y)m(y) = \sum_{d(x,y) \leq l} p_l(x,y)m(y) = 1 \qquad \forall x\in \Gamma,
\end{equation}
and that the kernel is symmetric:
\begin{equation} \label{symmetry}
 p_l(x,y) = p_l(y,x) \qquad \forall x,y\in \Gamma.
\end{equation}
For all functions $f$ on $\Gamma$, we define $P$ as the operator with kernel $p$, i.e.
\begin{equation}\label{defP} Pf(x) = \sum_{y\in \Gamma} p(x,y)f(y)m(y) \qquad \forall x\in \Gamma.\end{equation}
It is easily checked that $P^l$ is the operator with kernel $p_l$.

Since $p(x,y) \geq 0$ and \eqref{sum=1} holds, one has, for all $p\in [1,+\infty]$ ,
\begin{equation} \label{Pcont}
 \|P\|_{p\to p } \leq 1.
\end{equation}

\begin{rmq} \label{poweri-p}
Let $1\leq p<+\infty$. Since, for all $l\geq 0$, $\left\Vert P^l\right\Vert_{p\rightarrow p}\leq 1$, the operators $(I-P)^{\beta}$ and $(I+P)^{\beta}$ are $L^p$-bounded for all $\beta\geq 0$ (see \cite{CSC}).
\end{rmq}

\noindent We define a nonnegative Laplacian on $\Gamma$ by $\Delta = I-P$. One has then

\begin{rmq}
One can check that $\|\Delta\|_{1\to1} \leq 2$. Moreover, the previous remark states that $\Delta^\beta$ is $L^1(\Gamma)$-bounded. 
Note that the $L^1$-boundedness of the operators $\Delta^\beta$ is not true in the continuous setting (such as Riemannian manifolds), 
and makes some proofs of the present paper easier than in the case of Riemannian manifolds. 
In particular, we did not need then to prove similar results of the ones in \cite{AMM}.
\end{rmq}

\begin{equation}
 \begin{split}
 <(I-P)f,f>_{L^2(\Gamma)} & = \sum_{x,y\in \Gamma} p(x,y)(f(x)-f(y))f(x)m(x)m(y) \\
& = \frac{1}{2} \sum_{x,y\in \Gamma} p(x,y)|f(x)-f(y)|^2m(x)m(y),
 \end{split}
\end{equation}
where we use \eqref{sum=1} for the first equality and \eqref{symmetry} for the second one. The last calculus proves that the following operator
$$\nabla f(x) = \left( \frac{1}{2} \sum_{y\in \Gamma} p(x,y) |f(y)-f(x)|^2 m(y)\right)^{\frac{1}{2}},$$
called ``length of the gradient'' (and the definition of which is taken from \cite{CoulGrigor}), satisfies
\begin{equation}
\|\nabla f\|^2_{L^2(\Gamma)} =  <(I-P)f,f>_{L^2(\Gamma)} = \|\Delta^{\frac12}f\|_{L^2(\Gamma)}.
\end{equation}

\subsection{Assumptions on the graph}

\begin{defi}
 We say that $(\Gamma,\mu)$ satisfies the doubling property if there exists $C >0$ such that
\begin{equation} \label{DV} \tag{DV}
V(x,2r) \leq C V(x,r) \qquad \forall x\in \Gamma, \, \forall r>0.
\end{equation}
\end{defi}

\begin{prop} \label{propDV}
 Let $(\Gamma,\mu)$ satisfying the doubling property. Then there exists $d>0$ such that
\begin{equation} \label{PDV} V(x,\lambda r) \lesssim \lambda^d V(x,r) \qquad \forall x\in \Gamma, \, r>0 \text{ and } \lambda \geq 1.\end{equation}
We denote by $d_0$ the infimum of the $d$ satisfying \eqref{PDV}.
\end{prop}

\begin{defi}
  We say that $(\Gamma,\mu)$ (or $P$) satisfies \eqref{LB} if  there exists  $\epsilon = \epsilon_{LB}>0$  such that
\begin{equation} \label{LB} \tag{LB} p(x,x)m(x) \geq \epsilon \qquad \forall x\in \Gamma.\end{equation}
\end{defi}

\begin{rmq}
 In particular, the condition \eqref{LB} implies that  $-1$ does not belong to the $L^2$-spectrum of $P$, which implies in turn the analyticity of $P$ in $L^p(\Gamma)$, $1<p<+\infty$ (\cite{CSC}). 
\end{rmq}

 From now on,  all the graphs considered ( unless  explicitely stated) satisfy the doubling property and \eqref{LB}. 
In this context, Coulhon, Grigor'yan and Zucca proved in \cite{CGZ} (Theorem 4.1) that the following Davies-Gaffney estimate holds:

\begin{theo} \label{DGTheorem}
 Assume that $(\Gamma,\mu)$ satisfies \eqref{DV}. Then there  exist  $C,c>0$ such that for all subsets $E,F\subset \Gamma$ and all fonctions $f$ supported in $F$, one has
\begin{equation} \label{GUE} \tag{GUE}
\|P^{l-1}f\|_{L^2(E)} \leq  C\exp\left(-c \frac{d(E,F)^2}{l} \right) \|f\|_{L^2(F)} \qquad  \forall l\in \N^*.
\end{equation}
\end{theo}

The estimate \eqref{GUE}, also called Gaffney estimate, will be sufficient to prove most of the results of this paper. 
However, some results proven here can be improved if we assume the following  stronger  pointwise gaussian estimate:

\begin{defi}
We say that $(\Gamma,\mu)$ satisfies \eqref{UE} if there  exist  $C,c>0$ such that
\begin{equation} \label{UE} \tag{UE} p_{l-1}(x,y) \leq C \frac{1}{V(x,\sqrt{l})} \exp\left(-c \frac{d(x,y)^2}{l} \right) \qquad \forall x,y\in \Gamma, \, \forall l\in \N^*.\end{equation}
\end{defi}

\begin{rmq}
 Under \eqref{DV}, property \eqref{UE} is equivalent to
\begin{equation} \label{DUE} \tag{DUE}
 p_{l-1}(x,x) \leq \frac{C}{V(x,\sqrt l)} \qquad \forall x\in \Gamma, \, \forall l\in \N^*.
\end{equation}
The conjonction of \eqref{DV} and \eqref{UE} (or \eqref{DUE}) is also equivalent to some relative Faber-Krahn inequality (see \cite{CoulGrigor}).
\end{rmq}

\subsection{Definition of Hardy spaces on weighted graphs}

We introduce three different definitions for Hardy spaces. The  first two ones rely on  molecular decomposition.

\begin{defi} \label{molec1} Let $M\in \N^*$. When $\epsilon \in (0,+\infty)$, a function $a\in L^2(\Gamma)$ is called a $(BZ_1,M,\epsilon)$-molecule 
if there exist $s\in \N^*$, a $M$-tuple $(s_1,\dots,s_M)\in \bb s,2s\bn^M$, a ball $B$ of radius $\sqrt{s}$ and a function $b \in  L^2(\Gamma)$ such that
\begin{enumerate}[(i)]
 \item $a = (I-P^{s_1})\dots(I-P^{s_M})b$,
 \item $\ds \|b\|_{L^2(C_j(B))} \leq 2^{-j\epsilon} V(2^jB)^{-\frac{1}{2}}$, $\forall j\geq 1$.
\end{enumerate}

A function $a\in L^2(\Gamma)$ is called a $(BZ_1,M,\infty)$-molecule (or a $(BZ_1,M)$-atom) if there exist $s\in \N^*$, a $M$-tuple $ (s_1,\dots,s_M)\in \bb 1,M\bn^M$ a ball $B$ of radius $\sqrt{s}$ and a function $b \in L^2(\Gamma)$ supported in $B$ such that
\begin{enumerate}[(i)]
 \item $a = (I-P^{s_1})\dots(I-P^{s_M}) b$,
 \item $\ds \|b\|_{L^2}  = \|b\|_{L^2(B)} \leq  V(B)^{-\frac{1}{2}}$.
\end{enumerate}
We say that a $(BZ_1,M,\epsilon)$-molecule $a$ is associated with an integer $s$, a $M$-tuple $(s_1,\dots,s_M)$ and a ball $B$ when we want to refer to $s$, $(s_1,\dots,s_M)$ and $B$ given by the definition.
\end{defi}

 The second kind of molecules we consider are defined via the operators $I-(I+s\Delta)^{-1}$: 

\begin{defi} \label{molec2}
Let $M\in \N^*$. When $\epsilon \in (0,+\infty)$, a function $a\in L^2(\Gamma)$ is called a $(BZ_2,M,\epsilon)$-molecule if there exist $s\in \N^*$, a ball $B$ of radius $\sqrt{s}$ and a function $b \in  L^2(\Gamma)$ such that
\begin{enumerate}[(i)]
 \item $ a = [I-(I+s\Delta)^{-1}]^Mb$, 
 \item $\ds \|b\|_{L^2(C_j(B))} \leq 2^{-j\epsilon} V(2^jB)^{-\frac{1}{2}}$, $\forall j\geq 1$.
\end{enumerate}
A function $a\in L^2(\Gamma)$ is called a $(BZ_2,M,\infty)$-molecule (or a $(BZ_2,M)$-atom) if there exist $s\in \N^*$, a ball $B$ of radius $\sqrt{s}$ and a function $b \in L^2(\Gamma)$ supported in $B$ such that
\begin{enumerate}[(i)]
 \item $ a = [I-(I+s\Delta)^{-1}]^M b$, 
 \item $\ds \|b\|_{L^2}  = \|b\|_{L^2(B)} \leq  V(B)^{-\frac{1}{2}}$.
\end{enumerate}
We say that a $(BZ_2,M,\epsilon)$-molecule $a$ is associated with an integer $s$ and a ball $B$ when we want to refer to $s$ and $B$ given by the definition.
\end{defi}

\begin{rmq} \label{sizemolec}
\begin{itemize}
\item[$1.$]
When $b$ is the function occurring in Definition \ref{molec1} or in Definition \ref{molec2}, note that $\|b\|_{L^2} \lesssim V(B)^{-\frac{1}{2}}$.
\item[$2.$]
As will be seen in Proposition \ref{BoundedMolecules3} below, when $a$ is a molecule occurring in Definition \ref{molec1} or in Definition \ref{molec2}, one has $\|a\|_{L^1} \lesssim 1$.
\end{itemize}
\end{rmq}

\begin{defi}
Let $M\in \N^*$ and $\kappa\in \{1,2\}$.

 Let $\epsilon\in (0,+\infty]$. We say that $f$ belongs to $H^1_{BZ\kappa,M,\epsilon}(\Gamma)$ if $f$ admits a molecular $(BZ_\kappa,M,\epsilon)$-representation, 
that is if there exist a sequence $(\lambda_i)_{i\in \N} \in \ell^1$ and a sequence $(a_i)_{i\in \N}$ of $(BZ_\kappa,M,\epsilon)$-molecules such that
\begin{equation} \label{sumf}
f=\sum_{i=0}^\infty \lambda_i a_i
\end{equation}
where the convergence of the  series  to $f$  holds  pointwise. The space is outfitted with the norm
$$\|f\|_{H^1_{BZ\kappa,M,\epsilon}} = \inf\left\{ \sum_{j=0}^\infty |\lambda_j|, \ \sum_{j=0}^\infty \lambda_j a_j, \text{ is a molecular $(BZ_\kappa,M,\epsilon)$-representation of $f$} \right\}.$$
\end{defi}

\begin{prop} \label{H1complete} Let $M\in \N^*$ and $\kappa\in \{1,2\}$. 
Then the space $H^1_{BZ\kappa,M,\epsilon}(\Gamma)$ is complete. Moreover, $H^1_{BZ\kappa,M,\epsilon}(\Gamma)\subset L^1(\Gamma)$.
\end{prop}
\begin{dem}
 That $H^1_{BZ\kappa,M,\epsilon}(\Gamma)\subset L^1(\Gamma)$ follows at once from assertion $2$ in Remark \ref{sizemolec}, which shows that, if $f\in H^1_{BZ\kappa,M,\epsilon}(\Gamma)$, the series \eqref{sumf} converges in $L^1(\Gamma)$, and therefore converges to $f$ in $L^1(\Gamma)$. Moreover,  the space $H^1_{BZ\kappa,M,\epsilon}(\Gamma)$ is complete if it has the property
$$\sum_{j=0}^\infty \|f_j\|_{H^1_{BZ\kappa,M,\epsilon}} < +\infty \Longrightarrow \sum_{j=0}^\infty f_j  \mbox{ converges in }  H^1_{BZ\kappa,M,\epsilon}(\Gamma).$$
\noindent This fact is a straightforward consequence of the fact that $\|a\|_{L^1} \lesssim 1$ whenever $a$ is a molecule (see  Remark \ref{sizemolec} and  Proposition \ref{BoundedMolecules3}). See also the argument for the completeness of $H^1_L$ in \cite{HoMay}, p. 48.
\end{dem}

\begin{rmq} \label{remBZ}
The $BZ_{\kappa}$ molecules are molecules in the sense of Bernicot and Zhao in \cite{BZ} (and then $BZ_{\kappa}$ are Hardy spaces in the sense of Bernicot and Zhao). Note that the definition of molecules is slightly different from the one given in \cite{AMR}, \cite{HoMay} or \cite{HLMMY}.
The article \cite{BZ} provides some properties of the spaces $H^1_{BZ\kappa,M,\epsilon}$. In particular, under the assumption \eqref{UE}, these Hardy spaces are suited for $L^p$ interpolation  (see Remark \ref{remresults} below) . 
\end{rmq}
  
The third Hardy space is defined via quadratic functionals. 

\begin{defi}
 Define, for $\beta>0$, the quadratic functionals $L_\beta$ on $L^2(\Gamma)$ by
$$L_\beta f(x) =  \left( \sum_{(y,l) \in \gamma(x)} \frac{(l+1)^{2\beta-1}}{V(x,\sqrt{l+1})} |\Delta^{\beta} P^{l} f(y)|^2m(y) \right)^{\frac{1}{2}}$$
where $\gamma(x) = \left\{ (y,l) \in \Gamma \times \N, \, d(x,y)^2 \leq l \right\}$.
\end{defi}

\begin{rmq}
One can also use instead of $L_\beta$ the Lusin functional $\tilde L_\beta$ defined by
$$\tilde L_\beta f(x) =  \left( \sum_{(y,k) \in \tilde \gamma(x)} \frac{1}{(k+1)V(x,{k+1})} |(k^2\Delta)^{\beta} P^{k^2} f(y)m(y)|^2 \right)^{\frac{1}{2}}$$
where $\tilde \gamma(x) = \left\{ (y,k) \in \Gamma \times \N, \, d(x,y) \leq k \right\}$.

 The functionals  $L_\beta$ and $\tilde L_\beta$ are two different ways to  discretize  the ``countinuous'' Lusin functional defined  by  
\[\begin{split} L^c_{\beta} f(x) & = \left( \int_0^\infty \int_{d(y,x)^2< s} \frac{1}{sV(x,\sqrt{s})} |(s\Delta)^{\beta} e^{-s\Delta} f(y)|^2 d\mu(y)\, ds \right)^{\frac{1}{2}} \\
   & = \left( \int_0^\infty \int_{d(y,x)< t} \frac{1}{tV(x,t)} |(t^2\Delta)^{\beta} e^{-t^2\Delta} f(y)|^2 d\mu(y)\, dt \right)^{\frac{1}{2}}.
  \end{split} \]
\end{rmq}

\begin{defi}
 The space $E^1_{quad,\beta}(\Gamma)$ is defined by
$$E^1_{quad,\beta}(\Gamma) : = \left\{ f\in L^2(\Gamma), \, \|L_\beta f\|_{L^1} < +\infty\right\}.$$
It is outfitted with the norm
$$\|f\|_{H^1_{quad,\beta}} : = \|L_\beta f\|_{L^1}.$$
\end{defi}

\begin{rmq}
Notice that $\|f\|_{H^1_{quad,\beta}}$ is a norm because the null space of $\Delta$ is reduced to $\{0\}$ (because the set $\Gamma$ is infinite by assumption). So, if $k>\beta$ is an integer and $f\in L^2(\Gamma)$ is such that $\Delta^{\beta}f=0$, then $\Delta^kf=\Delta^{k-\beta}\Delta^{\beta}f=0$, so that $f=0$.
\end{rmq}

\begin{rmq}
 Replacing $L_\beta$ by $\tilde L_\beta$ in the definition of $E^1_{quad,\beta}$ yields an equivalent space $\tilde E^1_{quad,\beta}$,  in the sense that the sets are equal and the norms are equivalent.  
 The proof of this nontrivial fact can be done by adapting the proof of Theorem \ref{Main3} below (details are left to the reader).
\end{rmq}

\subsection{Definition of BMO spaces on weighted graphs}

Fix $x_0 \in \Gamma$ and let $B_0 = B(x_0,1) = \{x_0\}$. 
For $\epsilon>0$ and $M\in \N$,  for all functions $\phi\in L^2(\Gamma)$ which can be written as $\phi = \Delta^M \varphi$ for some function $\varphi \in L^2$, define 
$$\|\phi\|_{\mathcal M_{0}^{M,\epsilon}} : = \sup_{j\geq 1} \left[ 2^{j\epsilon} V(2^jB_0)^{\frac{1}{2}} \|\varphi\|_{L^2(C_j(B_0))}\right] \in [0,+\infty]. $$
We set then
$$\mathcal M_{0}^{M,\epsilon} : = \left\{ \phi = \Delta^M \varphi \in L^2(\Gamma), \ \|\phi\|_{\mathcal M_{0}^{M,\epsilon}} < +\infty \right\}.$$

\begin{defi}
 For any $M\in \N$, we set,
$$\mathcal E_M  = \bigcup_{\epsilon>0} (\mathcal M_{0}^{M,\epsilon})^*$$
and
$$\mathcal F_M  = \bigcap_{\epsilon>0} (\mathcal M_{0}^{M,\epsilon})^*.$$
\end{defi}

\begin{prop}
Let $M \in \N$, $s\in \N^*$ and $(s_1,\dots,s_M)\in \bb s,2s\bn^M$. 
If $f\in \mathcal E_M$, then the functions $(I-P^{s_1})\dots(I-P^{s_M})f$ and $(I-(I+s\Delta)^{-1})^Mf$ can be defined in the sense of distributions and are included in $L^2_{loc}(\Gamma)$.
\end{prop}

\begin{dem}
The proof of this fact is done in Lemma \ref{L2loc}.
\end{dem}

\begin{defi}
 Let $M\in \N$. Let $f\in \mathcal E_M$. 
 
 We say then that $f$ belongs to $BMO_{BZ1,M}(\Gamma)$ if 
\begin{equation} \label{BMO1def}\|f\|_{BMO_{BZ1,M}} : = \sup_{\begin{subarray}{c} s\in \N^*, \\(s_1,\dots,s_M)\in \bb s,2s\bn^M, \\ B \text{ of radius } \sqrt{s}\end{subarray}} \left( \frac{1}{V(B)} \sum_{x\in B} |(I-P^{s_1})\dots (I-P^{s_M}) f(x)|^2 m(x) \right)^{\frac{1}{2}} < +\infty.\end{equation}
 We say then that $f$ belongs to $BMO_{BZ2,M}(\Gamma)$ if 
\begin{equation} \label{BMO2def}\|f\|_{BMO_{BZ2,M}} : = \sup_{\begin{subarray}{c} s\in \N^*, \\ B \text{ of radius } \sqrt{s}\end{subarray}} \left( \frac{1}{V(B)} \sum_{x\in B} |[I-(I+s\Delta)^{-1}]^M f(x)|^2 m(x) \right)^{\frac{1}{2}} < +\infty .\end{equation}
\end{defi}

\subsection{Definition of Hardy spaces of 1-forms}

We define, for all $x\in \Gamma$, the set $T_x = \{(x,y)\in \Gamma^2, \, y\sim x\}$ and 
$$T_\Gamma = \bigcup_{x\in \Gamma} T_x = \{(x,y)\in \Gamma^2, \, y\sim x\}.$$

\begin{defi}
If $x\in \Gamma$, we define, for all $F_x$ defined on $T_x$ the norm 
$$\|F_x\|_{T_x} = \left( \frac{1}{2} \sum_{y\sim x} p(x,y) m(y) |F_x(x,y)|^2 \right)^\frac12.$$
Moreover, a function $F: \, T_\Gamma \to \R$ belongs to $L^p(T_\Gamma)$ if
\begin{enumerate}[(i)]
 \item $F$ is antisymmetric, that is $F(x,y) = -F(y,x)$ for all $x\sim y$,
 \item $\|F\|_{L^p(T_\Gamma)} < +\infty$, with
 $$\|F\|_{L^p(T_\Gamma)} = \left\|x\mapsto \|F(x,.)\|_{T_x} \right\|_{L^p(\Gamma)}.$$
\end{enumerate}
\end{defi}

\begin{defi} Let $f: \Gamma \to \R$ and $F: T_\Gamma \to \R$ be some functions. 
Define the operators $d$ and $d^*$ by
$$d f(x,y) : = f(x) - f(y) \qquad \forall (x,y)\in T_\Gamma$$
and
$$d^* F(x) : = \sum_{y\sim x} p(x,y) F(x,y) m(y) \qquad \forall x\in \Gamma.$$ 
\end{defi}

\begin{rmq}
It is plain to see that $d^*d = \Delta$ and $\|df(x,.)\|_{T_x} = \nabla f(x)$.
\end{rmq}

The definition of Hardy spaces of 1-forms is then similar to the case of functions.
First, we introduce Hardy spaces via molecules.

\begin{defi}
 Let $M\in \N$ and $\epsilon \in (0,+\infty)$. 
A function $a \in L^2(T_\Gamma)$ is called a $(BZ_2,M+\frac{1}{2},\epsilon)$-molecule if there exist $s\in \N^*$, a ball $B$ of radius $\sqrt{s}$ and a function $b \in L^2(\Gamma)$ such that
\begin{enumerate}[(i)]
 \item $a = s^{M+\frac{1}{2}} d\Delta^M(I-s\Delta)^{-M-\frac{1}{2}} b$;
 \item $\|b\|_{L^2(C_j(B))} \leq 2^{-j\epsilon} V(2^jB)^{\frac{1}{2}}$ for all $j\geq 1$.
\end{enumerate}
\end{defi}

\begin{rmq}
As in the case of functions, Corollary \ref{BoundedMolecules4}  below  implies  a uniform bound on the $L^1$ norm of  molecules, 
that is, for all $M\in \N$ and all $\epsilon \in (0,+\infty)$, there exists $C>$ such that each $(BZ_2,M,\epsilon)$-molecule $a$ satisfies
$$ \|a\|_{L^1(T_\Gamma)} \leq C.$$
\end{rmq}

\begin{defi}
 Let $M\in \N$ and $\epsilon \in (0,+\infty)$. We say that $F$ belongs to $H^1_{BZ2,M+\frac{1}{2},\epsilon}(T_\Gamma)$ if $F$ admits a molecular $(BZ_2,M+\frac{1}{2},\epsilon)$-representation, 
that is if there exist a sequence $(\lambda_i)_{i\in \N} \in \ell^1$ and a sequence $(a_i)_{i\in \N}$ of $(BZ_2,M+\frac{1}{2},\epsilon)$-molecules such that
$$F = \sum_{i=0}^\infty \lambda_i a_i$$
where the sum converges  pointwise on $T_{\Gamma}$ . The space is outfitted with the norm
$$\|f\|_{H^1_{BZ2,M+\frac{1}{2},\epsilon}} = \inf \left\{ \sum_{i=0}^\infty |\lambda_i|, \ \sum_{i=0}^\infty \lambda_ia_i \text{ is a molecular $(BZ_2,M+\frac{1}{2},\epsilon)$-representation of $f$}\right\}.$$
\end{defi}

\begin{rmq}
 The space $H^1_{BZ2,M+\frac{1}{2},\epsilon}(T_\Gamma)$ is complete. The argument is analogous to the one of Proposition \ref{H1complete}.
\end{rmq}

In order to define the Hardy spaces  of forms  associated with operators, we introduce the $L^2$ adapted Hardy spaces $H^2(T_\Gamma)$ defined as the  closure  in $L^2(T_\Gamma)$  of
$$ E^2(T_\Gamma) : = \{F \in L^2(T_\Gamma) , \, \exists f\in L^2(\Gamma): \, F = df\}.$$
Notice that $d\Delta^{-1}d^* = Id_{E^2(T_\Gamma)}$. The functional $d\Delta^{-1}d^*$ can be extended to a bounded operator on $H^2(T_\Gamma)$ and 
\begin{equation} \label{dDeltadstar}
d\Delta^{-1}d^* = Id_{H^2(T_\Gamma)}.
\end{equation}

\begin{prop} \label{danddstar}
 For all $p\in [1,+\infty]$,  the operator $d^*$ is bounded from $L^p(T_\Gamma)$ to $L^p(\Gamma)$.\par
\noindent The operator $d\Delta^{-\frac12}$ is an isometry from $L^2(\Gamma)$ to $L^2(T_\Gamma)$ (or $H^2(T_\Gamma)$), 
and the operator $\Delta^{-\frac12}d^*$ is an isometry from $H^2(T_\Gamma)$ to $L^2(\Gamma)$.
\end{prop}

\begin{dem}
 First, the $L^p$-boundedness of $d^*$ is provided by
\[\begin{split}
   \|d^*F\|_{L^p(\Gamma)}^p & = \sum_{x\in \Gamma} \left|\sum_{y\in \Gamma} p(x,y) m(y) F(x,y) \right|^p m(x) \\
& \lesssim \sum_{x\in \Gamma} \|F(x,.)\|_{T_x}^p m(x) = \|F\|_{L^p(T_\Gamma)}^p.
  \end{split}\]
The $L^2$-boundedness of $d\Delta^{-\frac12}$ is obtained by the calculus
\[\begin{split}
   \|d\Delta^{-\frac12}f\|_{L^2(T_\Gamma)}^2 & = \frac{1}{2} \sum_{x\sim y} p(x,y) |\Delta^{-\frac12}f(x)-\Delta^{-\frac12}f(y)|^2 m(x)m(y) \\
& = \|\nabla \Delta^{-\frac12}f\|_{L^2(\Gamma)}^2 = \|\Delta^{\frac{1}{2}}\Delta^{-\frac12}f \|_{L^2(\Gamma)}^2 \\
& = \|f\|_{L^2(\Gamma)}^2. 
  \end{split}\]
The $L^2$-boundedness of $\Delta^{-\frac12}d^*$ is then a consequence of \eqref{dDeltadstar}. Indeed, if $F \in H^2(\Gamma)$,
\[\begin{split}
   \|\Delta^{-\frac12}d^* F\|_{L^2(\Gamma)} & = \|d\Delta^{-\frac12}\Delta^{-\frac12}d^* F\|_{L^2(T_\Gamma)} \\
& = \|F\|_{L^2(T_\Gamma)}.
  \end{split}\]
\end{dem}

\begin{defi}
 The space $E^1_{quad,\beta}(T_\Gamma)$ is defined by
$$E^1_{quad,\beta}(T_\Gamma) := \left\{ F \in H^2(T_\Gamma), \, \|L_\beta [\Delta^{-\frac{1}{2}} d^*F]\|_{L^1} < +\infty\right\}$$
 equipped with  the norm
$$\|F\|_{H^1_{quad,\beta}} : = \|L_\beta [\Delta^{-\frac{1}{2}} d^*F]\|_{L^1}.$$
\end{defi}

Note that, if $\left\Vert F \right\Vert_{H^1_{quad,\beta}}=0$, one has $\Delta^{-1/2}d^{\ast}F=0$, so that $d\Delta^{-1/2}\Delta^{-1/2}d^{\ast}F=0$, which implies that $F=0$ since $F\in H^2(T_{\Gamma})$.
Moreover, check that for all $F\in H^2(T_\Gamma)$, $\|F\|_{H^1_{quad,\beta}} = \|\Delta^{-\frac12} d^* F\|_{H^1_{quad,\beta}}$

\subsection{Main results}

 \label{Mainresults}

In the following results, $\Gamma$ is assumed to satisfy \eqref{DV} and \eqref{LB}.
\begin{theo} \label{Main2}
 Let $M\in \N^*$. Then $BMO_{BZ1,M}(\Gamma) = BMO_{BZ2,M}(\Gamma)$. 
\end{theo}

\begin{theo} \label{Main1}
 Let $M\in \N^*$ and $\kappa\in \{1,2\}$. Let $\epsilon \in (0,+\infty]$.

Then the dual space of $H^1_{BZ\kappa,\epsilon}(\Gamma)$ is $BMO_{BZ1,M}(\Gamma) = BMO_{BZ2,M}(\Gamma)$. In particular, the spaces $H^1_{BZ\kappa,M,\epsilon}(\Gamma)$ depend  neither on $\epsilon$ nor on $\kappa$. 

Moreover, $BMO_{BZ\kappa}(\Gamma)$, initially defined as a subspace of $\mathcal E_M$, is actually included in $\mathcal F_M$.
\end{theo}

\begin{theo} \label{Main3}
 Let $\beta >0$ and $\kappa\in \left\{1,2\right\}$. The completion $H^1_{quad,\beta}(\Gamma)$ of $E^1_{quad,\beta}(\Gamma)$ in $L^1(\Gamma)$ exists. 
 Moreover, if $M\in (\frac{d_0}{4},+\infty)\cap \N^*$ and $\epsilon \in (0,+\infty]$, then the spaces $H^1_{BZ1,M,\epsilon}(\Gamma)$, $H^1_{BZ2,M,\epsilon}(\Gamma)$ and $H^1_{quad,\beta}(\Gamma)$ coincide. More precisely, we have 
 $$E^1_{quad,\beta}(\Gamma) = H^1_{BZ\kappa,M,\epsilon}(\Gamma) \cap L^2(\Gamma).$$
 \end{theo}
\noindent Once the equality $H^1_{BZ1,M,\epsilon}(\Gamma)=H^1_{BZ2,M,\epsilon}(\Gamma)= H^1_{quad,\beta}(\Gamma)$ is  established,  this space will be denoted by $H^1(\Gamma)$.

\begin{cor} \label{Main3b}
 Let $M_1,M_2>\frac{d_0}{4}$. Then we have the equality
$$BMO_{BZ1,M_1}(\Gamma) = BMO_{BZ2,M_2}(\Gamma).$$
\end{cor}

\begin{theo} \label{Main5}
 Let $\beta>0$. The completion $H^1_{quad,\beta}(T_\Gamma)$ of $E^1_{quad,\beta}(T_\Gamma)$ in $L^1(T_\Gamma)$ exists.
 
 Moreover, if $M\in (\frac{d_0}{4}-\frac{1}{2},+\infty)\cap \N$ and $\epsilon \in (0,+\infty)$, then the spaces $H^1_{BZ2,M+\frac{1}{2},\epsilon}(T_\Gamma)$ and $H^1_{quad,\beta}(T_\Gamma)$ coincide. More precisely, we have 
 $$E^1_{quad,\beta}(T_\Gamma) = H^1_{BZ2,M+\frac12,\epsilon}(T_\Gamma) \cap L^2(T_\Gamma). $$
\end{theo}
\noindent Again, the space $H^1_{BZ2,M+\frac{1}{2},\epsilon}(T_\Gamma)=H^1_{quad,\beta}(T_\Gamma)$ will be denoted by $H^1(T_\Gamma)$.

\begin{theo} \label{Main4}
  For this theorem only, assume furthermore that $(\Gamma, \mu)$ satisfies \eqref{UE}. 
Then $M$ can be choosen  arbitrarily  in $\N^*$ in Theorem \ref{Main3} and Corollary \ref{Main3b}, $M$ can be choosen  arbitrarily  in $\N$ in Theorem \ref{Main5}.
\end{theo}

\begin{theo} \label{Main6}
The Riesz transform $d\Delta^{-\frac{1}{2}}$ is bounded from $H^1(\Gamma)$ to $H^1(T_\Gamma)$. 
As a consequence the Riesz transform $\nabla\Delta^{-\frac{1}{2}}$ is bounded from $H^1(\Gamma)$ to $L^1(\Gamma)$.
\end{theo}

\begin{dem}
By definition, 
$$\|d\Delta^{-\frac{1}{2}}f\|_{H^1(T_\Gamma)} \simeq \|d\Delta^{-\frac{1}{2}}f\|_{H^1_{quad,1}(T_\Gamma)} = \|\Delta^{-\frac{1}{2}} d^*d\Delta^{-\frac{1}{2}}f\|_{H^1_{quad,1}(\Gamma)}= \|f\|_{H^1_{quad,1}(\Gamma)} \simeq \|f\|_{H^1(\Gamma)}.$$
Therefore, $d\Delta^{-\frac{1}{2}}$ is $H^1$-bounded. 
Moreover, $\|\nabla\Delta^{-\frac{1}{2}} f\|_{L^1(\Gamma)} = \|d\Delta^{-\frac{1}{2}} f\|_{L^1(T_\Gamma)} \lesssim \|d\Delta^{-\frac12}f\|_{H^1(T_\Gamma)}$.
Indeed, the uniform $L^1$-bound of $(BZ2,M+\frac12,\epsilon)$-molecules (see Corollary \ref{BoundedMolecules4}) yields
$$H^1(T_\Gamma) = H^1_{BZ2,M+\frac{1}{2},\epsilon}(T_\Gamma) \hookrightarrow L^1(T_\Gamma)$$
for any $M>\frac{d_0}{4}-\frac{1}{2}$. 
\end{dem}

\begin{rmq} \label{remresults}\hspace{1cm} 

\begin{enumerate}[(a)]
 \item It is easily checked that under \eqref{UE}, the Hardy space $H^1(\Gamma) = H^1_{BZ2,1,\infty}(\Gamma)$ satisfies the assumption of Theorem 5.3 in \cite{BZ}. 
As a consequence, the interpolation between $H^1(\Gamma)$ and $L^2(\Gamma)$ provides the spaces $L^p(\Gamma)$, $1<p<2$. 
Together with Theorem \ref{Main6}, we can recover the main result of \cite{Russ}, that is: under \eqref{UE}, the Riesz transform $\nabla\Delta^{-\frac12}$ is $L^p$-bounded for all $p\in (1,2]$.

\item An interesting byproduct of Theorem \ref{Main3} is the equality, for any $\epsilon \in (0,+\infty]$ and any $M>\frac{d_0}{4}$, between the spaces 
$H^1_{BZ\kappa,M,\epsilon}(\Gamma) \cap L^2(\Gamma)$ and $E^1_{BZ\kappa,M,\epsilon}(\Gamma)$ defined by
{\small $$E^1_{BZ\kappa,M,\epsilon}(\Gamma) := \left\{ f\in L^2(\Gamma), \ \sum_{j=0}^\infty \lambda_j a_j \text{ is a molecular $(BZ_\kappa,M,\epsilon)$-representation of $f$ and the  series  converges in $L^2(\Gamma)$ } \right\}$$}
and outfitted with the norm
{\small $$\|f\|_{E^1_{BZ\kappa,M,\epsilon}} = \inf\left\{ \sum_{i\in \N}|\lambda_i|, \ \sum_{j=0}^\infty \lambda_j a_j \text{ is a molecular $(BZ_\kappa,M,\epsilon)$-representation of $f$ and the  series  converges in $L^2(\Gamma)$ } \right\}.$$}
We have similar byproducts of Theorems \ref{Main5} and \ref{Main4}.  Precise  statements and proofs are done in Corollary \ref{Main30}. 

As a consequence, the completion of $E^1_{BZ\kappa,M,\epsilon}(\Gamma)$ in $L^1(\Gamma)$ exists and is equal to $H^1_{BZ\kappa,M,\epsilon}$. 
On Riemannian manifolds  or in more general contexts,  the proof of this fact is  much  more complicated and is the main result of \cite{AMM}.  Let us emphasize that the proofs of our main results does not go through the $E^1_{BZ\kappa,M,\epsilon}$ spaces.

\item We may replace (i) in the definition of $(BZ_2,M,\epsilon)$-molecules by

       (i') $a = (I-(I+s_1\Delta)^{-1})\dots ((I-(I+s_M\Delta)^{-1}) b$, where $(s_1,\dots,s_M)\in [s,2s]^M$

       or

       (i'')  $a = (I-(I+r^2\Delta)^{-1})^M b$, where $r$ is the radius of the ball $B$ (or the smallest integer greater than $\sqrt{s}$)

       and still get the same space $H^1_{BZ2,M,\epsilon}(\Gamma)$.

 \item However, when $M\geq 3$, it is unclear whether replacing  item $(i)$ of  the definition of $(BZ_1,M,\epsilon)$-molecules by

       (i') $a= (I-P^s)^M$ 

       yields the same space $H^1_{BZ1,M,\epsilon}(\Gamma)$.       
\end{enumerate}
\end{rmq}

Section 2 is devoted to the proof of auxiliary results that will be useful for the next sections.
The proof of Theorem \ref{Main2} is treated in paragraph 3.2 and the proof of Theorem \ref{Main1} is done in paragraph 3.3. 
In the last section, we establish Theorems \ref{Main3}, \ref{Main5} and \ref{Main4}.

\subsection{Comparison with other papers}

\begin{itemize}
 \item {\bf Comparison with \cite{AMR}}: In \cite{AMR}, the authors proved  analogous results  (that is the $H^1$ boundedness of the Riesz transform under very weak assumptions  and the various characterizations of $H^1$ ) on Riemannian manifolds. Some differences between the two papers can be noted. First, $BMO$ spaces are not considered there. They also choose to define some Hardy spaces via  tent  spaces (while we  prefer to use  Lusin functionals).
Contrary to us, they introduced the spaces $H^p$,  for all $p\in [1,+\infty]$ , and proved  that these spaces form an interpolation scale for the complex method. 

 \item {\bf Comparison with \cite{HLMMY}}: This article develops Hardy and BMO spaces adapted to a  symmetric  operator $L$ in a general context of doubling measure spaces when the semigroup generated by $L$ satisfies $L^2$ Gaffney estimates.
However,  on graphs,  it is unclear whether these $L^2$ Gaffney estimates for the semigroup generated by the Laplacian hold or not.  Yet, Coulhon, Grigor'yan and Zucca proved in \cite{CGZ} that we have $L^2$ Gaffney type estimates for the discrete iterates of Markov operators and we only rely on these estimates in the present paper.

 \item {\bf Comparison with \cite{BD}}: First of all, as in \cite{HLMMY}, there are no results about Hardy spaces on 1-forms  and the authors do not prove  the $H^1$ boundedness of the Riesz transforms.
Then, as said in the introduction, they assume in all their paper a pointwise gaussian bound of the Markov kernel while it is not required for most of our results. 
Moreover, the results of the present paper stated under \eqref{UE} are stronger that those stated in \cite{BD}. 
Indeed, in the results stated in \cite{BD}, the constant $M$ need to be greater than $\frac{d_0}{2}$ 
while, in the present paper, we used the pointwise gaussian bound in order to get rid of the dependance of $M$ on the ``dimension'' $d_0$.

Besides, the definitions of their Hardy spaces and ours {\it a priori } differ. Let us begin with the Hardy spaces defined via molecules. 
For convenience, we introduce a new definition of molecules.

\begin{defi} \label{molec3}
Let $M\in \N^*$ and $\epsilon \in (0,+\infty)$. A function $a\in L^2(\Gamma)$ is called a $(HM,M,\epsilon)$-molecule if there exist a ball $B$ of radius $r \in \N^*$ and a function $b \in  L^2(\Gamma)$ such that
\begin{enumerate}[(i)]
 \item $a = [r^2\Delta]^Mb$,
 \item $\ds \|[r^2\Delta]^k b\|_{L^2(C_j(B))} \leq 2^{-j\epsilon} V(2^jB)^{-\frac{1}{2}}$, $\forall j\in \N^*$, $\forall k\in \bb 0,M\bn$.
\end{enumerate}
The space $H^1_{HM,M,\epsilon}(\Gamma)$ is then defined in the same way  as  $H^1_{BZ\kappa,M,\epsilon}(\Gamma)$.
\end{defi}

Using methods developed in \cite{HLMMY} and in the present paper, it can be proved that, if $M>\frac{d_0}{4}$ (or if $M\in \N^*$ if we assume the extra condition \eqref{UE}), 
there is equality between the spaces $H^1_{HM,M,\epsilon}(\Gamma)$ and $H^1_{quad,1}(\Gamma) = H^1(\Gamma)$.  The proofs are similar to those of the present paper.  
The molecules introduced by Bui and Duong - we call them $(BD,M,\epsilon)$-molecules - are the $(HM,M,\epsilon)$-molecules where we replaced $r^2$ by $r$ in (i) and (ii).
It is easily checked that a $(BD,M,\epsilon)$-molecule is a $(HM,M,\epsilon)$-molecule and hence, under assumption \eqref{UE}, our Hardy spaces are bigger than theirs.

Since they  proved (as we do here) that Hardy spaces defined with molecules and with quadratic functionals coincide, the Hardy spaces via quadratic functionals in \cite{BD} are also different from ours. 
Indeed, 
our Hardy spaces are of parabolic type (heat kernel) while those of \cite{BD} are modelled on the Poisson semigroup.
Furthermore, they  only consider  one Lusin functional,  while  we consider a family of Lusin functionals (indexed by $\beta>0$), 
and the independance of Hardy spaces $H^1_{quad,\beta}(\Gamma)$ with respect to $\beta$ is a key point of the proof of the boundedness of Riesz transforms.\par

\noindent  {\bf Acknowledgements: } the author is grateful to E. Russ for comments and suggestions that improved the paper.
He would also like to thank P. Auscher and A. Morris for interesting discussions. 

\end{itemize}

\section{Preliminary results}
\subsection{$L^2$-convergence}

\begin{prop} \label{L2convergence}
 Let $\beta >0$. Let $P$ satisfying \eqref{LB}. One has the following convergence: for all $f\in L^2(\Gamma)$,
$$\sum_{k=0}^N a_k (I-P)^\beta P^k f \xrightarrow{N\to + \infty} f \qquad \text{ in } L^2(\Gamma)$$
where $\sum a_k z^k$ is the Taylor series of the function $(1-z)^{-\beta}$.
\end{prop}

\begin{rmq}
This result extends Lemma 1.13 in \cite{BRuss}. It provides a discrete version of the identity
$$f = c_\beta \int_0^\infty (t\Delta)^\beta e^{-t\Delta} f dt.$$
\end{rmq}

\begin{cor}  \label{L2convergenceCor}
 Let $(\Gamma,\mu)$ a weighted graph. One has the following convergence: for all $f\in L^2(\Gamma)$,
$$\sum_{k=0}^N a_k (I-P^2)^\beta P^{2k} f \xrightarrow{N\to + \infty} f \qquad \text{ in } L^2(\Gamma)$$
\end{cor}

\begin{dem} (Proposition \ref{L2convergence})

First, notice that Corollary \ref{L2convergenceCor} is an immediate consequence of Proposition \ref{L2convergence} since $P^2$  is a Markov operator satisfying \eqref{LB} (see \cite{CGZ}).
 
Let $f \in L^2(\Gamma)$. Let us check the behavior of
\begin{equation} \label{L2converg}\left\| \left[\sum_{k=0}^N a_k (I-P)^\beta P^k -I\right] f \right\|_{L^2}\end{equation}
when $N \to +\infty$. Since $\|P\|_{2\to2} = 1$ and $P$ satisfies \eqref{LB}, there exists $a>-1$ such that
$$P = \int_{a}^1 \lambda dE(\lambda).$$
Thus
\begin{equation} \label{FunctCalculus}\begin{split}
   \left\| \left[\sum_{k=0}^N a_k (I-P)^\beta P^k -I\right] f \right\|_{L^2}^2  
& = \int_a^1 \left[\sum_{k=0}^N a_k (1-\lambda)^\beta \lambda^k -1\right]^2 dE_{ff}(\lambda).  \end{split} \end{equation}
However, 
$$\sum_{k=0}^N a_k (1-\lambda)^\beta \lambda^k \xrightarrow{N\to \infty} \left\{\begin{array}{ll} 1 &  \text{ for all }\lambda \in [a,1) \\ 0 & \text{ if } \lambda = 1\end{array} \right. $$
and since the sum is nonnegative and increasing in $N$, then
$$ \left|\sum_{k=0}^N a_k (1-\lambda)^\beta \lambda^k -1\right| \leq 1 \qquad \forall \lambda \in [a,1].$$
We use this result in \eqref{FunctCalculus} to get the uniform bound
\begin{equation}\label{UniformBoundedness} \begin{split}
 \left\| \left[\sum_{k=0}^N a_k (I-P)^\beta P^k -I\right] f \right\|_{L^2}^2 & \leq \int_a^{1} dE_{ff}(\lambda) = \|f\|_{L^2}^2 .               
\end{split} \end{equation}

Let us  focus on  \eqref{L2converg} when we furthermore  assume that $f \in \Ro(\Delta)$, that is $f = \Delta g$ for some $g\in L^2(\Gamma)$.
The identity \eqref{FunctCalculus} reads as
\[\begin{split} 
 \left\| \left[\sum_{k=0}^N a_k (I-P)^\beta P^k -I\right] f \right\|_{L^2}^2 
& = \int_a^1 \left[\sum_{k=0}^N a_k (1-\lambda)^{\beta+1} \lambda^k -(1-\lambda)\right]^2 dE_{gg}(\lambda).
  \end{split} \]
Yet, $\sum_{k=0}^N a_k (1-\lambda)^{\beta+1} \lambda^k -(1-\lambda)$ converges uniformly to 0 for all $\lambda \in [a,1]$.

Consequently, for all $\epsilon >0$, there exists $N_0$ such that, for all $N>N_0$,
$$\left\| \left[\sum_{k=0}^N a_k (I-P)^\beta P^k -I\right] f \right\|_{L^2}^2 \leq \epsilon \int_a^1 dE_{gg}(\lambda) = \epsilon \|g\|_{L^2}^2.$$
This implies
\begin{equation} \label{L2conv}\sum_{k=0}^N a_k (I-P)^\beta P^k f  \xrightarrow{N\to\infty} f \quad \text{ in $L^2$ and for all } f\in \Ro(\Delta). \end{equation}

Since $L^2 = \overline{\mathcal R(\Delta)}$, the combination of \eqref{UniformBoundedness} and \eqref{L2conv} provides the desired conclusion. 
Indeed, \eqref{L2conv} provides the $L^2$-convergence on the dense space $\mathcal R(\Delta)$ and the uniform boundedness \eqref{L2conv} allows us to extend the convergence to $L^2(\Gamma)$.
\end{dem}

\subsection{Davies-Gaffney estimates}

\label{DGEstimates}

\begin{defi}
 We say that a family of operators $(A_s)_{s\in \N}$ satisfies Davies-Gaffney estimates 
if there exist three constants $C,c,\eta>0$ such that for all subsets $E,F\subset \Gamma$ and all functions $f$ supported in $F$, there holds
\begin{equation} \label{DGEdef}\left\|A_s f \right\|_{L^2(E)} \leq C \exp\left(-c\left[\frac{d(E,F)^2}{s}\right]^\eta\right) \|f\|_{L^2}. \end{equation}
\end{defi}

Hofmann and Martell proved in \cite[Lemma 2.3]{HoMar} the following result about Davies-Gaffney estimates:

\begin{prop} \label{compositionGDE}
 If $A_s$ and $B_t$ satisfy Davies-Gaffney estimates, then there exist $C,c,\eta>0$ such that for all subsets $E,F\subset \Gamma$ and all functions $f$ supported in $F$, there holds
\begin{equation} \label{DGEasbt}\left\|A_sB_t f \right\|_{L^2(E)} \leq C \exp\left(-c\left[\frac{d(E,F)^2}{s+t}\right]^\eta\right) \|f\|_{L^2} \end{equation}
In particular, $(A_sB_s)_{s\in \N}$ satisfies Davies-Gaffney estimates.

More precisely, if $\eta_A$ and $\eta_B$ are the constants involved in \eqref{DGEdef} respectively for $A_s$ and $B_t$ , then the constant $\eta$ that occurs in \eqref{DGEasbt} can be choosen equal to $ \min\{\eta_A,\eta_B\}$. 
\end{prop}

\begin{prop} \label{GaffneyEstimatesResolvant}
Let $M\in \N$. The following families of operators satisfy the Davies-Gaffney estimates 
\begin{enumerate}[(i)]
 \item $\ds \prod_{i=1}^M \left(\frac{1}{t^i_s}\sum_{k=0}^{t^i_s} P^k \right)$, where for all $i\in \bb 1,M\bn$, $t^i_s \in \bb 1,2s\bn$,
 \item $\ds \prod_{i=1}^M (I-P^{t^i_s})$, where for all $i\in \bb 1,M\bn$, $t^i_s \in \bb s,2s\bn$, 
 \item $(s\Delta)^MP^s$,
 \item $(I+s\Delta)^{-M}$,
 \item $(I-(I+s\Delta)^{-1})^M = (s\Delta)^M(I+s\Delta)^{-M}$.
\end{enumerate}
In (i), (ii) and (iii), the parameter $\eta$ is equal to 1
 and in (iv) and (v), $\eta$ is equal to $\frac12$.
\end{prop}

\begin{dem} 
(i) and (ii) are  direct consequences  of \eqref{GUE} and Proposition \ref{compositionGDE}. 
 Assertion  (iii) is the consequence of \eqref{GUE} and \eqref{LB} and a proof can be found in \cite{Dungey}.

We turn now to the proof of (iv) and (v). According to Proposition \ref{compositionGDE}, it remains to show the Davies-Gaffney estimates for $(I+s\Delta)^{-1}$,  and since $s\Delta(I+s\Delta)^{-1} = I - (I+s\Delta)^{-1}$, 
it is enough to deal with  $(I+s\Delta)^{-1}$.  The  $L^2$-functional calculus provides the identity
\begin{equation} \begin{split}
   (I+s\Delta)^{-1} f&  = \frac{1}{1+s} \left(I-\frac{s}{1+s}P \right)^{-1}f \\
& =  \sum_{k=0}^{+\infty} \frac{1}{1+s}\left(\frac{s}{1+s}\right)^k P^kf,
  \end{split}\end{equation}
where the convergence  holds  in $L^2(\Gamma)$.

Let $f$ be a function supported in $F$. Then, one has with the Gaffney-Davies estimates \eqref{GUE}:
\[\begin{split}
   \|(I+s\Delta)^{-1} f\|_{L^2(E)} & \lesssim \sum_{k=0}^{+\infty}  \frac{1}{1+s}  \left(\frac{s}{1+s}\right)^k \|P^kf\|_{L^2(E)} \\
& \lesssim \sum_{k=0}^{+\infty}  \frac{1}{1+s} \left(\frac{s}{1+s}\right)^k \exp\left(-c \frac{d(E,F)^2}{1+k} \right) \|f\|_{L^2(F)} \\
& \lesssim \|f\|_{L^2(F)} \left[ \sum_{k=0}^{s} \frac{1}{1+s} \exp\left(-c \frac{d(E,F)^2}{1+k} - c' \frac{k}{1+s} \right) \right. \\
& \qquad + \left. \sum_{k=s}^{+\infty} \frac{1+s}{(1+k)^2} \exp\left(-c \frac{d(E,F)^2}{1+k} - c' \frac{k}{1+s} \right) \right].
  \end{split}\]
Yet, the function $\psi: k\in \R^+ \mapsto c \frac{d(E,F)^2}{1+k} + c' \frac{k}{1+s}$ is bounded from below and
$$\psi(k) \gtrsim \frac{d(E,F)}{\sqrt{1+s}}.$$
Hence, the use of Lemma \ref{ExpDecay} proved in the appendix yields
\[\begin{split}
   \|(I+s\Delta)^{-1} f\|_{L^2(E)} & \lesssim \|f\|_{L^2(F)} \exp\left(-c \frac{d(E,F)}{\sqrt{1+s}} \right) \left[ \sum_{k=0}^{s} \frac{1}{1+s} + \sum_{k=s}^{+\infty} \frac{1+s}{(1+k)^2}\right] \\
& \lesssim  \|f\|_{L^2(F)} \exp\left(-c \frac{d(E,F)}{\sqrt{1+s}} \right).
  \end{split}\]
\end{dem}

\begin{prop} \label{BoundedMolecules3}
 Let $\kappa \in \{1,2\}$. Let $a$ be a $(BZ_\kappa,M,\epsilon)$-molecule. Then
$$\|a\|_{L^1} \lesssim 1 \qquad \text{ and } \qquad \|a\|_{L^2(C_j(B))} \lesssim \frac{2^{-j\epsilon}}{V(2^jB)^{\frac{1}{2}}} \quad \forall j\in \N^*.$$
\end{prop}

\begin{dem} We will only prove the case where $\kappa = 1$. The case $\kappa = 2$ is proven similarly and will therefore be skipped.

 Since 
$$\|a\|_{L^1} \leq \sum_{j\geq 1} V(2^{j+1}B)^{\frac12}\|a\|_{L^2(C_j(B))},$$
we only need to check the second fact.
Let $s\in \N$, $(s_1,\dots,s_M)\in \bb s,2s\bn^M$ and a ball $B$ associated with the molecule $a$.

Define $\ds \tilde C_j(B) = \bigcup_{k=j-1}^{j+1} C_k(B)$ and observe that $d(C_j(B),\Gamma\backslash \tilde C_j(B)) \gtrsim 2^j\sqrt{s}$. Then Proposition \ref{GaffneyEstimatesResolvant} provides
\[\begin{split}
   \|a\|_{L^2(C_j(B))} & \leq \|(I-P^{s_1})\dots(I-P^{s_M}) [b\1_{\tilde C_j(B)}]\|_{L^2(C_j(B))} + \|(I-P^{s_1})\dots(I-P^{s_M}) [b\1_{\Gamma \backslash \tilde C_j(B)}]\|_{L^2(C_j(B))} \\
& \lesssim \|b\|_{L^2(\tilde C_j(B))} + e^{-c4^{j}} \|b\|_{L^2} \\
& \lesssim \frac{2^{-j\epsilon}}{V(2^jB)^\frac12} + \frac{e^{-c4^{j}}}{V(B)^\frac12} \\
& \lesssim \frac{2^{-j\epsilon}}{V(2^jB)^\frac12}.
  \end{split}\]
\end{dem}

\subsection{Gaffney estimates for the gradient} 

\begin{prop} \label{DGGradient1}
 Let $(\Gamma,\mu)$ satisfying \eqref{LB} (note that \eqref{DV} is not assumed here). Let $c>0$ such that
\begin{equation} \label{cFormula} \frac{8ce^{8c}}{\epsilon_{LB}} \leq 1. \end{equation}
There exists $C>0$ such that for all subsets $F\subset \Gamma$ and all $f$ supported in $F$, one has
$$\left\|P^k f e^{c\frac{d^2(.,F)}{k+1}}\right\|_{L^2} \leq C \|f\|_{L^2}.$$
\end{prop}

The proof of Proposition \ref{DGGradient1}    is based on the following result of Coulhon, Grigor'yan and Zucca:

\begin{lem} \label{DGGradient2}
  Let $(\Gamma,\mu)$ satisfying \eqref{LB}. Let $(k,x) \mapsto g_k(x)$ be a positive function on $\N \times \Gamma$. 
Then, for all finitely supported functions $f\in L^2(\Gamma)$ and for all $k\in \N$,
$$\left\|\sqrt{g_{k+1}} P^{k+1}f \right\|_{L^2}^2 - \left\|\sqrt{g_{k}} P^{k}f \right\|_{L^2}^2 
\leq \sum_{x\in \Gamma} |P^kf(x)|^2 \left(g_{k+1}(x) - g_k(x) + \dfrac{|\nabla g_{k+1}(x)|^2}{4\epsilon_{LB} g_{k+1}(x)} \right)m(x)$$
\end{lem}

\begin{dem}
This fact is actually  established in the proof of \cite[Theorem 2.2, pp. 566-567]{CGZ}. 
\end{dem}

\begin{dem} (Proposition \ref{DGGradient1}).

First, let us prove the result for $f$ supported in a finite set $F \subset \Gamma$. Let $f$ (finitely supported and) supported in $F$.
We wish to use Lemma \ref{DGGradient2} with 
$$g_{k}(x) = e^{2c\frac{d^2(x,F)}{k+1}}.$$
Check that, with Taylor-Lagrange inequality
\[\begin{split} 
  g_{k+1}(x) - g_k(x) & \leq \max_{t\in [k,k+1]} \left\{ - \frac{2cd^2(x,F)}{(t+1)^2} e^{2c\frac{d^2(x,F)}{t+1}}\right\} \\
& \quad = - 2c \left(\frac{d(x,F)}{k+2}\right)^2 g_{k+1}(x).
  \end{split}\]
In the same way, one has
$$\nabla g_{k+1}(x) \leq \frac{4c[d(x,F)+1]}{k+2} e^{2c\frac{[d(x,F)+1]^2}{k+2}}.$$
Since $f$ is supported in $F$, then $P^k f$ is supported in $\{x\in \Gamma, \, d(x,F) \leq k\}$. 
As a consequence, we can assume in the previous calculus that $d(x,F) \leq k$ and thus
$$\frac{[d(x,F)+1]^2}{k+2} \leq \frac{d^2(x,F)}{k+2} + 2.$$
Then
$$\frac{|\nabla g_{k+1}(x)|^2}{4\epsilon_{LB} g_{k+1}(x)} \leq \left(\frac{[d(x,F)+1]}{k+2}\right)^2 \frac{4c^2e^{8c}}{\epsilon_{LB}} g_{k+1}(x).$$
First case: $d(x,F) \geq 1$, then
$$\frac{|\nabla g_{k+1}(x)|^2}{4\epsilon_{LB} g_{k+1}(x)} \leq \left(\frac{d(x,F)}{k+2}\right)^2 \frac{16c^2e^{8c}}{\epsilon_{LB}} g_{k+1}(x)$$
and by \eqref{cFormula},
$$g_{k+1}(x) - g_k(x) + \dfrac{|\nabla g_{k+1}(x)|^2}{4\epsilon_{LB} g_{k+1}(x)}  \leq 0.$$
Second case, $d(x,F) =0$, then
$$g_{k+1}(x) - g_k(x) + \dfrac{|\nabla g_{k+1}(x)|^2}{4\epsilon_{LB} g_{k+1}(x)} \leq \frac{1}{(k+2)^2} \frac{16c^2e^{8c}}{\epsilon_{LB}^2} \leq  \frac{2c}{(k+2)^2}$$
In all cases, one has then $P^k f(x) = 0$ or
$$g_{k+1}(x) - g_k(x) + \dfrac{|\nabla g_{k+1}(x)|^2}{4\epsilon_{LB} g_{k+1}(x)} \leq  \frac{2c}{(k+2)^2}.$$
Lemma \ref{DGGradient2} yields
$$\left\| P^{k+1}f e^{c\frac{d^2(.,F)}{k+2}}\right\|_{L^2}^2 - \left\|P^{k}f e^{c\frac{d^2(.,F)}{k+1}}\right\|_{L^2}^2 
\leq \frac{2c}{(k+2)^2} \|P^k f\|_{L^2}^2,$$
and hence, by induction,
$$\left\| P^{k}f e^{c\frac{d^2(.,F)}{k+1}}\right\|_{L^2}^2 \leq \|f\|_{L^2}^2 +  \sum_{l=0}^{k-1} \frac{2c}{(l+2)^2} \|P^l f\|_{L^2}^2 \lesssim \|f\|_{L^2}^2.$$

\smallskip

Consider now a general $f \in L^2(\Gamma)$. Without loss of generality, we can assume that $f$ is nonnegative. 
Let $(\Gamma_i)_{i\in \N}$ an increasing sequence of finite subsets of $\Gamma$ such that $\bigcup_{i=0}^\infty \Gamma_i = \Gamma$.
Let $f_i = f\1_{\Gamma_i}$. One has then for any $x\in \Gamma$ and $k\in \N$,
$$f_i \uparrow f \qquad \text{ and } P^k f_i \uparrow P^kf.$$
By the monotone convergence theorem, we obtain,
$$\left\| P^{k}f_i e^{c\frac{d^2(.,F)}{k+1}}\right\|_{L^2}^2  \uparrow \left\| P^{k}f e^{c\frac{d^2(.,F)}{k+1}}\right\|_{L^2}^2$$
so that
\[\begin{split}
   \left\| P^{k}f e^{c\frac{d^2(.,F)}{k+1}}\right\|_{L^2}^2 & = \lim_{i\to \infty} \left\| P^{k}f_i e^{c\frac{d^2(.,F)}{k+1}}\right\|_{L^2}^2 \\
& \lesssim \sup_{i\in \N} \|f_i\|_{L^2}^2 \\
& \quad = \|f\|_{L^2}^2.
  \end{split}\]
\end{dem}

\begin{prop} \label{DGGradient3}
 Let $(\Gamma,\mu)$ satisfying \eqref{LB} (note that \eqref{DV} is not assumed here). Let $c>0$ as in Proposition \ref{DGGradient1}.
There exists $C>0$ such that for all subsets $F\subset \Gamma$ and all functions $f$ supported in $F$, one has
$$\left\|\nabla P^k f e^{\frac{c}{2}\frac{d^2(.,F)}{k+1}}\right\|_{L^2} \leq C \frac{\|f\|_{L^2}}{\sqrt{k+1}}.$$
\end{prop}

\begin{dem}
 The proof of this proposition is very similar to the one of Lemma 7 in \cite{Russ}. We define
$$I = I_k(f) : = \left\|\nabla P^k f e^{\frac{c}{2}\frac{d^2(.,F)}{k+1}} \right\|_{L^2}^2.$$
One has then
\[\begin{split}
   I & = \sum_{x,y\in \Gamma} p(x,y) |P^k(x)-P^k(y)|^2 e^{c\frac{d^2(x,F)}{k+1}} m(x)m(y) \\
& = \sum_{x,y\in \Gamma} p(x,y) [P^kf(x)-P^kf(y)]P^kf(x) e^{c\frac{d^2(x,F)}{k+1}} m(x)m(y) \\
& \qquad - \sum_{x,y\in \Gamma} p(x,y) [P^kf(x)-P^kf(y)]P^kf(y) e^{c\frac{d^2(x,F)}{k+1}} m(x)m(y) \\
& = 2\sum_{x,y\in \Gamma} p(x,y) [P^kf(x)-P^kf(y)]P^kf(x) e^{c\frac{d^2(x,F)}{k+1}} m(x)m(y) \\
& \qquad + \sum_{x,y\in \Gamma} p(x,y) [P^kf(x)-P^kf(y)]P^kf(x) \left[e^{c\frac{d^2(y,F)}{k+1}}- e^{c\frac{d^2(x,F)}{k+1}}\right]m(x)m(y) \\
& : = 2I_1 + I_2.
  \end{split}\]
We first estimate $I_1$. One has
\[\begin{split}
   I_1 & = \sum_{x\in \Gamma} P^kf(x) e^{c\frac{d^2(x,F)}{k+1}} m(x) \sum_{y\in \Gamma} p(x,y) [P^kf(x)-P^kf(y)] m(y) \\
& = \sum_{x\in \Gamma} (I-P)P^kf(x) P^kf(x) e^{c\frac{d^2(x,F)}{k+1}} m(x).
  \end{split}\]
Consequently, with the analyticity of $P$ and Proposition \ref{DGGradient1}, we get
\begin{equation} \label{I1Estimate}
\begin{split}
I_1 & \leq  \|(I-P)P^k f\|_{L^2} \left\|P^k f e^{c\frac{d^2(.,F)}{k+1}}\right\|_{L^2} \\
& \lesssim \frac{1}{k+1} \|f\|^2_{L^2}.
\end{split}
\end{equation}
We now turn to the estimate of $I_2$. One has, since $d(x,y) \leq 1$ (otherwise $p(x,y) = 0$),
\[\begin{split}
   \left| e^{c\frac{d^2(y,F)}{k+1}}- e^{c\frac{d^2(x,F)}{k+1}} \right| & \leq 2\frac{[d(x,F)+1]}{k+1} e^{c\frac{d^2(x,F)}{k+1}} \\
& \lesssim \frac{1}{\sqrt{k+1}} e^{\frac{3c}{2}\frac{d^2(x,F)}{k+1}}.
  \end{split} \]
Since $f$ is supported in $F$, $P^kf$ is supported in $\{x\in \Gamma, \, d(x,F)  \leq k\}$. 
Consequently, we can assume that $d(x,F) \leq k+1$ so that
$$\frac{[d(x,F)+1]^2}{k+1} \leq \frac{d^2(x,F)}{k+1} + 2.$$
Therefore, the term $I_2$ can be estimated by
\begin{equation}\label{I2Estimate} \begin{split}
   |I_2| & \lesssim \frac{1}{\sqrt{k+1}} \sum_{x,y \in \Gamma} |P^kf(x)-P^kf(y)| |P^kf(x)| e^{\frac{3c}{2}\frac{d^2(x,F)}{k+1}} m(x)m(y) \\
& \lesssim \frac{1}{\sqrt{k+1}} \left(\sum_{x,y \in \Gamma} |P^kf(x)-P^kf(y)|^2  e^{c\frac{d^2(x,F)}{k+1}} m(x)m(y)\right)^\frac12 
 \left(\sum_{x,y \in \Gamma} |P^kf(x)|^2 e^{2c\frac{d^2(x,F)}{k+1}} m(x)m(y)\right)^\frac12 \\
& \qquad = \frac{1}{\sqrt{k+1}} \sqrt{I} \left\|P^k f e^{c\frac{d^2(.,F)}{k+1}}\right\|_{L^2} \\
& \lesssim \sqrt{\frac{I}{k+1}} \|f\|_{L^2},
  \end{split}\end{equation}
where we used again Proposition \ref{DGGradient1} for the last line.

The estimates \eqref{I1Estimate} and \eqref{I2Estimate} yield
$$I \lesssim \frac{1}{k+1} \|f\|_{L^2}^2 + \sqrt{\frac{I}{k+1}} \|f\|_{L^2},$$
that is
$$I \lesssim \frac{1}{k+1} \|f\|_{L^2}^2,$$
which is the desired conclusion.
\end{dem}

\begin{cor} \label{GaffneyEstimatesNabla}
  Let $(\Gamma,\mu)$ satisfying \eqref{LB} (note that \eqref{DV} is not assumed here). Let $M\in \N$. The following families of operators satisfy the Davies-Gaffney estimates
\begin{enumerate}[(i)]
 \item $\ds s^{M+\frac{1}{2}} \nabla \Delta^M P^s$,
 \item $s^{M+\frac{1}{2}}\nabla\Delta^M(I+s\Delta)^{-M-\frac{1}{2}}$.
\end{enumerate}
\end{cor}

\begin{dem}
 According to Propositions \ref{compositionGDE} and \ref{GaffneyEstimatesResolvant}, 
it is enough to check that $\sqrt{s}\nabla P^s$ and $\sqrt{s}\nabla(I+s\Delta)^{-\frac{1}{2}}$ satisfy Davies-Gaffney estimates.

Indeed, Proposition \ref{DGGradient3} yields, if $E,F\subset \Gamma$, $f$ supported in $F$ and $c>0$ satisfy \eqref{cFormula}
\[\begin{split}
   \|\sqrt{s}\nabla P^sf\|_{L^2(E)} e^{\frac{c}{2}\frac{d(E,F)}{s+1}} & \leq \sqrt{s}\left\|\nabla P^s f e^{\frac{c}{2}\frac{d^2(.,F)}{k+1}}\right\|_{L^2} \\
& \lesssim \frac{\sqrt{s}}{\sqrt{s+1}} \leq 1.
  \end{split}\]
It suffices now to check that $\sqrt{s}\nabla(I+s\Delta)^{-\frac{1}{2}}$ satisfies Davies-Gaffney estimates. First notice that
\[\begin{split} \|\sqrt{s}\nabla(I+s\Delta)^{-\frac{1}{2}} f\|_{L^2} & = \|(s\Delta)^{\frac{1}{2}}(I+s\Delta)^{-\frac{1}{2}}f\|_{L^2}\\
   & = \left\|(I-(I+s\Delta)^{-1})^{\frac{1}{2}}f\right\|_{L^2} \\
& \leq \|f\|_{L^2}.
\end{split}\]
Then the family of operators $\sqrt{s}\nabla(I+s\Delta)^{-\frac{1}{2}}$ is $L^2$-uniformly bounded. 
Hence, we can suppose without loss of generality that $d(E,F)^2 \geq 1+s$.
Write,
\[\begin{split}
   (I+s\Delta)^{-\frac{1}{2}}f & = \frac{1}{\sqrt{1+s}} \left(I-\frac{s}{1+s} P\right)^{-\frac{1}{2}}f \\
& = \frac{1}{\sqrt{1+s}} \sum_{k=0}^\infty a_k \left( \frac{s}{1+s}\right)^k P^kf
  \end{split}\]
where $\sum a_k z^k$ is the Taylor serie of the function $(1-z)^{-\frac{1}{2}}$ and the convergence  holds  in $L^2(\Gamma)$. 
Note that $a_k \simeq \sqrt{k+1}$ (see for example \cite{Fen1}, Lemma B.1) and
\[\begin{split}
   \|\sqrt{s}\nabla(I+s\Delta)^{-\frac{1}{2}}f\|_{L^2(E)} 
& \lesssim \frac{\sqrt{s}}{\sqrt{1+s}} \sum_{k=0}^\infty \frac{1}{\sqrt{1+k}} \left( \frac{s}{1+s}\right)^k \|\nabla P^kf\|_{L^2(E)} \\
& \lesssim \|f\|_{L^2} \sum_{k=0}^\infty \frac{1}{1+k} \left( \frac{s}{1+s}\right)^k e^{-c\frac{d(E,F)^2}{1+k}} \\
& \lesssim \|f\|_{L^2} \frac{1}{d(E,F)^2} \sum_{k=0}^\infty \left( \frac{s}{1+s}\right)^k e^{-c\frac{d(E,F)^2}{1+k}}\\
& \lesssim \|f\|_{L^2} \frac{1}{d(E,F)^2} \left[ \sum_{k=0}^s  e^{-c\left[\frac{d(E,F)^2}{1+k} + \frac{k}{1+s}\right]} \right.\\
& \qquad \left. + \sum_{k=s+1}^\infty \left(\frac{1+s}{1+k}\right)^2 e^{-c\left[\frac{d(E,F)^2}{1+k} + \frac{k}{1+s}\right]} \right]
\end{split}\]
where we used (i) for the second estimate and Lemma \ref{ExpDecay} for the last one. 

Arguing as in the proof of Proposition \ref{GaffneyEstimatesResolvant}, we find
\[\begin{split}
   \|\sqrt{s}\nabla(I+s\Delta)^{-\frac{1}{2}}\|_{L^2(E)} 
& \lesssim \|f\|_{L^2} \frac{1+s}{d(E,F)^2} e^{-c\frac{d(E,F)}{\sqrt{1+s}}}\\
& \lesssim \|f\|_{L^2},
\end{split}\]
since we assumed that $d(E,F)^2 \geq 1+s$.
\end{dem}

\begin{cor} \label{BoundedMolecules4}
Let $M\in \N$. Then if $a=s^{M+\frac{1}{2}}d\Delta^M(I+s\Delta)^{-M-\frac{1}{2}}b$ is a $(BZ_2,M+\frac{1}{2},\epsilon)$-molecule associated with the ball $B$, then
$$\|a\|_{L^1(T_\Gamma)} \lesssim 1 \qquad \text{ and } \qquad  \|a\|_{L^2(T_{C_j(B)})} \lesssim \frac{2^{-j\epsilon}}{V(2^jB)^{\frac{1}{2}}} \quad \forall j\in \N^*.$$
\end{cor}

\begin{dem}
 First, notice that
$$\|a\|_{L^1(T_\Gamma)} \leq \sum_{j\geq 1} V(2^{j+1}B)^\frac12 \|x\mapsto\|a(x,.)\|_{T_x}\|_{L^2(C_j(B))}.$$
Then it remains to check the last claim, that is 
$$\|a\|_{L^2(T_{C_j(B)})} : = \|x\mapsto\|a(x,.)\|_{T_x}\|_{C_j(B)} \lesssim \frac{2^{-j\epsilon}}{V(2^jB)^{\frac{1}{2}}}.$$
Since $a = s^{M+\frac{1}{2}}d\Delta^M(I+s\Delta)^{-M-\frac{1}{2}} b$, then 
$$x\mapsto\|a(x,.)\|_{T_x} = s^{M+\frac{1}{2}} \nabla \Delta^M(I+s\Delta)^{-M-\frac{1}{2}}b(x).$$
We conclude as in Proposition \ref{BoundedMolecules3}, using the Davies-Gaffney estimates provided by Corollary \ref{GaffneyEstimatesNabla}.
\end{dem}

\subsection{Off diagonal decay for Littlewood-Paley functionals} 

\begin{lem} \label{l2l1}
Let $M>0$ and $\alpha \in [0,1]$. Define $\mathcal A = \{ (A_l^{d,u})_{l\in \N^*}, \ d\in \R_+, \,  u\in \N \}$,  where, for all $l\geq 1$, 
$$A_l^{d,u} = l^{\alpha} \dfrac{\exp\left( -\frac{d}{l+u}\right)}{(l+u)^{1+M}}.$$
Then there exists $C = C_{M,\alpha}$ such that
$$\left( \sum_{l\in \N^*} \frac{1}{l} a_l^2 \right) \leq C \sum_{l\in \N^*} \frac{1}{l} a_l \qquad \forall (a_l)_l \in \mathcal A.$$
\end{lem}

\begin{dem}
 The proof is similar to Proposition C.2 in \cite{Fen1}.
\end{dem}

\begin{lem} \label{OffDiagonalDecay}
 Let $M\in \N^*$ and $\beta>0$. Then there exists $C_{M,\beta}$ such that for all sets $E,F\subset \Gamma$, all $f$ supported in $F$, all $s\in \N$ and all $M$-tuples $(s_1,\dots,s_M) \in \bb s,2s\bn^M$, one has
$$\left\| L_\beta (I-P^{s_1})\dots(I-P^{s_M})f \right\|_{L^2(E)} \leq C_{M,\beta} \left(1+\frac{d(E,F)^2}{s}\right)^{-M} \|f\|_{L^2}.$$
\end{lem}

\begin{dem}
The proof follows the ideas of \cite{Fen1} Lemma 3.3 (or \cite{BRuss} Lemma 3.2 if $\beta = 1$). 
First, since $L_\beta$ and $(I-P^{s_1})\dots(I-P^{s_M})$ are $L^2$-bounded (uniformly in $s$) and without loss of generality, we can assume that $s\leq d(E,F)^2$.

Denote by $\eta$ the only integer such that $\eta + 1 \geq \beta + M > \eta \geq 0$. Notice that $M-\eta \leq 1-\beta <1$ and thus $M-\eta \leq 0$.

 We use the following fact, which is an immediate consequence of Proposition \ref{L2convergence}
$$ \Delta^{\beta+M} f =  (I-P)^{\beta+M} f = \sum_{k\geq 0} a_k P^k(I-P)^{\eta+1}f \qquad \forall f\in L^2(\Gamma)$$
where $\sum a_k z^k$ is the Taylor serie of the function $(1-z)^{\beta+M-\eta -1}$ 
Notice that if $\beta+M$ is an integer, then $a_k = \delta_0(k)$.

By the use of the generalized Minkowski inequality, we get
\[ \begin{split}
 &   \left\| L_\beta (I-P^s)^M f\right\|_{L^2(E)} \\
& \qquad \leq \sum_{k\geq 0} a_k  \left( \sum_{l\geq 1} l^{2\beta-1} \sum_{x\in E} \frac{m(x)}{V(x,\sqrt l)} \sum_{y\in B(x,\sqrt{l})} m(y) |\Delta^{1 + \eta-M}(I-P^{s_1})\dots(I-P^{s_M})P^{k+l-1}f(y)|^2 \right)^{\frac{1}{2}} \\
& \qquad \lesssim s^M \sup_{t\in \bb 0,2Ms\bn}  \sum_{k\geq 0} a_k  \left( \sum_{l\geq 1} l^{2\beta-1} \sum_{x\in E} \frac{m(x)}{V(x,\sqrt l)} \sum_{y\in B(x,\sqrt{l})} m(y) |\Delta^{1 + \eta}P^{k+l+t-1} f(y)|^2 \right)^{\frac{1}{2}} \\
& \qquad \lesssim s^M \sup_{t\in \bb 0,2Ms\bn} \sum_{k\geq 0} a_k  \left( \sum_{l\geq 1} l^{2\beta-1} \sum_{y\in D_{l}(E)}  m(y) |\Delta^{1 + \eta}P^{k+l+t-1} f(y)|^2 \sum_{x\in B(y,\sqrt{l})} \frac{m(x)}{V(x,\sqrt l)} \right)^{\frac{1}{2}} \\
& \qquad \lesssim s^M \sup_{t\in \bb 0,2Ms\bn} \sum_{k\geq 0} a_k  \left( \sum_{l\geq 1} l^{2\beta-1} \|\Delta^{1 + \eta}P^{k+l+t-1}f\|^2_{L^2(D_{l}(E))} \right)^{\frac{1}{2}} \\
& \qquad  \qquad : =  s^M \sup_{t\in \bb 0,Ms\bn} \Lambda(t) 
   \end{split}\]
where $D_{l}(E) = \{y \in \Gamma, \, \dist(y,E) < \sqrt{l}\}$, and where we notice that $\ds \sum_{x\in B(y,\sqrt{l})} \frac{m(x)}{V(x,\sqrt{l})} \lesssim 1$ with the doubling property.

\ms

\begin{itemize}
 \item[{\bf 1-}]{\bf Estimate when} $\mathbf{ l < \frac{d(E,F)^2}{4}}$ 

\smallskip

The important point here is to notice that $\dist(F,D_{l}(E)) \geq \frac{1}{2} d(E,F) \gtrsim d(E,F)$. 
Then, using Davies-Gaffney estimates (Proposition \ref{GaffneyEstimatesResolvant}, (iii) ) , we may obtain

\begin{equation}\begin{split} \label{ODDcase1}
\|\Delta^{1 + \eta}P^{k+l+t-1}f\|_{L^2(D_{l}(E))} & \lesssim  \dfrac{\exp\left( -c \frac{d(E,F)^2}{l+k+t}\right)}{(l+k+t)^{(1+\eta)}} \|f\|_{L^2} \\
& \quad \leq  l^{M-\eta} \dfrac{\exp\left( -c \frac{d(E,F)^2}{l+k+t}\right)}{(l+k+t)^{1+M}} \|f\|_{L^2}
\end{split}\end{equation}
since $M-\eta \leq 0$.

\item[{\bf 2-}]{\bf Estimate when} $\mathbf{ l \geq \frac{d(E,F)^2}{4}}$ 

\smallskip

We use the analyticity of $P$ to obtain,

\begin{equation}\begin{split} \label{ODDcase2}
\|\Delta^{1 + \eta}P^{k+l+t-1}f\|_{L^2(D_{l}(E))} & \leq \|(I-P)^{1+\eta} P^{k+l+t-1} f\|_{L^2(\Gamma)} \\
 & \lesssim \frac{1}{(k+l+t)^{1+\eta}} \|f\|_{L^2} \\
   & \lesssim l^{M-\eta} \dfrac{1}{(k+l+t)^{1+M}} \|f\|_{L^2} \\
& \lesssim l^{M-\eta}  \dfrac{\exp\left( -c \frac{d(E,F)^2}{l+k+t}\right)}{(l+k+t)^{1+M}} \|f\|_{L^2} 
  \end{split}\end{equation}
where the third line is due to  $M-\eta \leq 0$ and the last one holds because $l+k \gtrsim d(E,F)^2$.

\item[{\bf 3-}]{\bf Conclusion}

The first two steps imply the following estimate on $\Lambda(t)$:
\[\begin{split}
  \Lambda(t) & \lesssim \|f\|_{L^2(F)} \sum_{k\geq 0} a_k  \left( \sum_{l\geq 1} \frac{1}{l} \left| l^{(\beta+M-\eta)} \dfrac{\exp\left( -c \frac{d(E,F)^2}{l+k+t}\right)}{(k+l+t)^{1+M}}\right|^2 \right)^{\frac{1}{2}} \\
  & \lesssim \|f\|_{L^2(F)} \sum_{k\geq 0} a_k \sum_{l\geq 1} \frac{1}{l} l^{(\beta+M-\eta)} \dfrac{\exp\left( -c \frac{d(E,F)^2}{l+k+t}\right)}{(k+l+t)^{1+M}} 
  \end{split}\]
where we used Lemma \ref{l2l1} for the last line (indeed, $\beta+M-\eta \in (0,1]$). Check that Thus since
$$\sum_{k=0}^{m-1} a_k (m-k)^{\beta+M-\eta-1} \lesssim 1.$$
Indeed, when  $\beta+M-\eta=1$, the result is obvious. Otherwise, it is a consequence of the fact that $a_k \simeq k^{\eta-M-\beta}$ (see Lemma B.1 in \cite{Fen1}).
Hence, one has
\[\begin{split}
   \Lambda(t) & \lesssim \|f\|_{L^2(F)} \sum_{m\geq 1} \dfrac{\exp\left( -c \frac{d(E,F)^2}{m+t}\right)}{(m+t)^{1+M}} \\
& \qquad = d(E,F)^{-2(1+M)}\|f\|_{L^2(F)} \sum_{m\geq 1} \dfrac{d(E,F)^{2(1+M)}}{(m+t)^{1+M}} \exp\left( -c \frac{d(E,F)^2}{m+t}\right)\\
& \lesssim d(E,F)^{-2(1+M)} \|f\|_{L^2(F)} \left[ \sum_{m= 1}^{d(E,F)^2} 1 + \sum_{m> d(E,F)^2}  \dfrac{d(E,F)^{2(1+M)}}{(m+t)^{1+M}}\right] \\
& \lesssim d(E,F)^{-2M} \|f\|_{L^2(F)}.
  \end{split}\]
As a consequence,
$$\left\| L_\beta (I-P^s)^M f\right\|_{L^2(E)} \lesssim \left(\frac{d(E,F)^2}{s}\right)^{-M} \|f\|_{L^2(F)}$$
which is the desired conclusion.
\end{itemize}
\end{dem}

\begin{lem} \label{OffDiagonalDecayTer}
 Let $M\in \R_+^*$ and $\beta >0$ such that either  $M\in \N$ or $\beta \geq 1$. Define the Littlewood-Paley functional $G_\beta$ on  $L^2(\Gamma)$  by
$$G_\beta f(x) = \left(\sum_{l\leq 1} l^{2\beta-1}|\Delta^\beta P^{l-1} f(x)|^2\right)^\frac12 \qquad  \forall x\in \Gamma . $$
Then there exists $C_{M}>0$ such that for all sets $E,F\subset \Gamma$, all functions $f$ supported in $F$ and all $s\in \N$, one has
$$\left\| G_\beta (s\Delta)^M f \right\|_{L^2(E)} \leq C_{M} \left(\frac{d(E,F)^2}{s}\right)^{-M} \|f\|_{L^2}.$$
\end{lem}

\begin{dem}
The proof is similar to the one of Lemma \ref{OffDiagonalDecay}. Notice that $\ds \sup_{t\in \bb 0,Ms\bn}\Lambda(t)$ is replaced by $\Lambda := \Lambda(0)$. 
Then the end of the calculus is the same provided that $\eta -M \geq 0$, which is the case under our assumption on $M$ and $\beta$. See also \cite[Lemma 3.3]{Fen1}.
\end{dem}

\begin{lem} \label{OffDiagonalDecay4}
Let $M\in \N$. Then there exists $C_{M}>0$ such that for all sets $E,F\subset \Gamma$, all $f$ supported in $F$ and all $s\in \N$, one has
$$\left\| L_{\frac12} (I-(I+s\Delta)^{-1})^{M+\frac{1}{2}} f \right\|_{L^2(E)} \leq C_{M} \left(1+\frac{d(E,F)^2}{s}\right)^{-M-\frac{1}{2}} \|f\|_{L^2}.$$
\end{lem}

\begin{dem}
Since $L_\frac12$ and $(I-(I+s\Delta)^{-1})^{M+\frac{1}{2}}$ are $L^2$-bounded (uniformly in $s$) and without loss of generality, we can assume that $s\leq d(E,F)^2$.

 We use the following computation,
\begin{equation} \label{ResolvantEquality} \begin{split}
   (I+s\Delta)^{-M-\frac{1}{2}}f & = ((1+s)I-sP)^{-M-\frac{1}{2}} = (1+s)^{-M-\frac{1}{2}} \left(I-\frac{s}{1+s}P\right)^{-M-\frac{1}{2}}f \\
& = (1+s)^{-M-\frac{1}{2}} \sum_{k\geq 0} a_k \left( \frac{s}{1+s}\right)^k P^kf
  \end{split}\end{equation}
where $\sum a_k z^k$ is the Taylor series of the function $(1-z)^{-M-\frac{1}{2}}$ and the convergence  holds  in $L^2(\Gamma)$.

By the use of the generalized Minkowski inequality, we get
\[ \begin{split}
 &   \left\| L_\frac12 (I-(I+s\Delta)^{-1})^{M+\frac{1}{2}} f\right\|_{L^2(E)} \\
& \qquad \leq \frac{s^{M+\frac{1}{2}}}{(1+s)^{M+\frac{1}{2}}}\sum_{k\geq 0} a_k \left( \frac{s}{1+s}\right)^k \left( \sum_{l\geq 1}  \sum_{x\in E} \frac{m(x)}{V(x,\sqrt l)} \sum_{y\in B(x,\sqrt{l})} m(y) |\Delta^{1 + M}P^{k+l-1}f(y)|^2 \right)^{\frac{1}{2}} \\
& \qquad \leq \frac{s^{M+\frac{1}{2}}}{(1+s)^{M+\frac{1}{2}}}\sum_{k\geq 0} a_k \left( \frac{s}{1+s}\right)^k  \left( \sum_{l\geq 1}  \sum_{y\in D_{l}(E)}  m(y) |\Delta^{1 + M}P^{k+l-1} f(y)|^2 \sum_{x\in B(y,\sqrt{l})} \frac{m(x)}{V(x,\sqrt l)} \right)^{\frac{1}{2}} \\
& \qquad \lesssim \frac{s^{M+\frac{1}{2}}}{(1+s)^{M+\frac{1}{2}}}\sum_{k\geq 0} a_k \left( \frac{s}{1+s}\right)^k  \left( \sum_{l\geq 1}  \|\Delta^{1 + M}P^{k+l-1}f\|^2_{L^2(D_{l}(E))} \right)^{\frac{1}{2}}.
   \end{split}\]

When $l<\frac{d^2(E,F)}{4}$, notice that $d(F,D_l(E)) \gtrsim d(E,F)$ so that
\begin{equation} \label{ODDCase5} \|\Delta^{1 + M}P^{k+l-1}f\|_{L^2(D_{l}(E))} \lesssim \frac{\exp(-c\frac{d^2(E,F)}{l+k})}{(l+k)^{M+1}} \|f\|_{L^2}.\end{equation}
Moreover, when $l\geq\frac{d^2(E,F)}{4}$, one has
\begin{equation} \label{ODDCase6} \begin{split}
   \|\Delta^{1 + M}P^{k+l-1}f\|_{L^2(D_{l}(E))} & \leq \|\Delta^{1 + M}P^{k+l-1}f\|_{L^2} \\
& \lesssim \frac{1}{(l+k)^{M+1}} \|f\|_{L^2} \\
& \lesssim \frac{\exp(-c\frac{d^2(E,F)}{l+k})}{(l+k)^{M+1}} \|f\|_{L^2}.
  \end{split}\end{equation}
As a consequence
\[\begin{split}
   \left\| L_\frac12 (I-(I+s\Delta)^{-1})^{M+\frac{1}{2}} f\right\|_{L^2(E)} 
& \lesssim \frac{s^{M+\frac{1}{2}}}{(1+s)^{M+\frac{1}{2}}}\|f\|_{L^2} \sum_{k\geq 0} a_k \left( \frac{s}{1+s}\right)^k  \left( \sum_{l\geq 1} \frac{\exp(-c\frac{d^2(E,F)}{l+k})}{(l+k)^{2(M+1)}}  \right)^{\frac{1}{2}} \\
& \lesssim \frac{s^{M+\frac{1}{2}}}{(1+s)^{M+\frac{1}{2}}} \|f\|_{L^2} \sum_{k\geq 0} a_k \left( \frac{s}{1+s}\right)^k  \left( \sum_{n\geq 1} \frac{\exp(-c\frac{d^2(E,F)}{n})}{n^{2(M+1)}}  \right)^{\frac{1}{2}} \\
& \lesssim \frac{s^{M+\frac{1}{2}}}{d(E,F)^{2M+1}}\frac{1}{(1+s)^{M+\frac{1}{2}}} \|f\|_{L^2} \sum_{k\geq 0} a_k \left( \frac{s}{1+s}\right)^k  \\
& \quad = \left(\frac{s}{d(E,F)^2}\right)^{M+\frac{1}{2}} (1+s(1-1))^{-M-\frac{1}{2}}\|f\|_{L^2} \\
& \quad = \left(\frac{s}{d(E,F)^2}\right)^{M+\frac{1}{2}}\|f\|_{L^2}.
  \end{split}\]
\end{dem}

Let us now  recall  a result that can be found in \cite{Fen1}, Theorem 1.4.

\begin{prop} \label{GTinequalities}
Assume that $(\Gamma,\mu)$ satisfy \eqref{UE}. Let $K>0$ and $j\in \N$. There exist $C,c>0$ such that for all sets $E,F\in \Gamma$ and all $x_0\in \Gamma$ all $l\in \N^*$ satisfying
\begin{equation}\label{assumption1GT}\sup_{y\in F} d(x_0,y) \leq K d(E,F)\end{equation}
or
\begin{equation}\label{assumption2GT}\sup_{y\in F} d(x_0,y) \leq K \sqrt{l}\end{equation}
 and all functions $f$ supported in $F$, there holds
$$\|\Delta^j P^{l-1} f\|_{L^2(E)} \leq \frac{C}{l^j} \frac{1}{V(x_0,\sqrt{l})^{\frac{1}{2}}} e^{-c\frac{d(E,F)^2}{l}} \|f\|_{L^1(F)}$$
and
$$\|\nabla\Delta^j P^{l-1} f\|_{L^2(E)} \leq \frac{C}{l^{j+\frac{1}{2}}} \frac{1}{V(x_0,\sqrt{l})^{\frac{1}{2}}} e^{-c\frac{d(E,F)^2}{l}} \|f\|_{L^1(F)}.$$
\end{prop}

\begin{lem} \label{OffDiagonalDecayBis}
Assume that $(\Gamma,\mu)$ satisfy \eqref{UE}. For all $M\in \N^*$ and all $\beta >0$, there exists $C_{M}>0$ such that for  all disjoint  sets $E,F\in \Gamma$ and all $x_0$ satisfying \eqref{assumption1GT}, all $f$ supported in $F$ and all $s\in \N^*$, one has
$$\left\| L_\beta (I-P^s)^M f \right\|_{L^2(E)} \leq \frac{C_{M}}{V(x_0,d(E,F))^{\frac{1}{2}}} \left(\frac{d(E,F)^2}{s}\right)^{-M} \|f\|_{L^1}.$$
\end{lem}

\begin{dem}
The proof of this Lemma is similar to the one of Lemma \ref{OffDiagonalDecay} and we only indicate the main changes.

When $l < \frac{d(E,F)^2}{4}$, replace first \eqref{ODDcase1} by
\begin{equation}\begin{split} \label{ODDcase1bis}
\|\Delta^{1 + \eta}P^{k+l+t-1}f\|_{L^2(D_{l}(E))} & \lesssim  \frac{1}{V(x_0,\sqrt{k+l+t})^\frac12}\dfrac{\exp\left( -c \frac{d(E,F)^2}{l+k+t}\right)}{(l+k+t)^{(1+\eta)}} \|f\|_{L^1} \\
& \lesssim  \frac{l^{M-\eta} }{V(x_0,d(E,F))^\frac12} \dfrac{\exp\left( -c \frac{d(E,F)^2}{l+k+t}\right)}{(l+k+t)^{1+M}} \|f\|_{L^1}
\end{split}\end{equation}
where the second line holds because $M-\eta \leq 0$ and the first one  holds by Proposition \ref{GTinequalities}.  
Indeed, there exists $K>0$ such that 
$$\sup_{y\in F} d(x_0,y) \leq K d(E,F).$$
Thus, $x_0$, $D_{l}(E)$ and $F$ satisfy \eqref{assumption1GT} with constant $4K$ .

When $l \geq \frac{d(E,F)^2}{4}$, replace also \eqref{ODDcase2} by
\begin{equation}\begin{split} \label{ODDcase2bis}
\|\Delta^{1 + \eta}P^{k+l+t-1}f\|_{L^2(D_{l}(E))} & \lesssim \frac{1}{V(x_0,\sqrt{k+l+t})^\frac12}\frac{1}{(k+l+t)^{1+\eta}} \|f\|_{L^1} \\
   & \lesssim \frac{l^{M-\eta}}{V(x_0,d(E,F))^\frac12} \dfrac{\exp\left( -c \frac{d(E,F)^2}{l+k+t}\right)}{(k+l+t)^{1+M}} \|f\|_{L^1}
  \end{split}\end{equation}
where the first line follows from Proposition \ref{GTinequalities}, since $x_0$, $F$ and $k+l+t$ satisfy \eqref{assumption2GT}, 
and the second line to the facts that $k+l \gtrsim d(E,F)^2$ and $M-\eta \leq 0$.
\end{dem}

\begin{lem} \label{OffDiagonalDecay5}
Assume that $(\Gamma,\mu)$ satisfy \eqref{UE}. For all $M>0$ , there exists $C_{M}$ such that for all sets $E,F\in \Gamma$ and all $x_0$ satisfying \eqref{assumption1GT}, all $f$ supported in $F$ and all $s\in \N^*$, one has
$$\left\| L_\frac12  (I-(I+s\Delta)^{-1})^{M+\frac{1}{2}} f \right\|_{L^2(E)} \leq \frac{C_{M}}{V(x_0,d(E,F))^{\frac{1}{2}}} \left(\frac{d(E,F)^2}{s}\right)^{-M-\frac12} \|f\|_{L^1}.$$
\end{lem}

\begin{dem}
The proof of this Lemma is similar to the one  of Lemma \ref{OffDiagonalDecay4} and we only indicate the main changes. 

When $l < \frac{d(E,F)^2}{4}$, replace \eqref{ODDCase5} by
\begin{equation}\begin{split} 
\|\Delta^{1+M}P^{k+l-1}f\|_{L^2(D_{l}(E))} & \lesssim  \frac{1}{V(x_0,\sqrt{k+l})^\frac12}\dfrac{\exp\left( -c \frac{d(E,F)^2}{l+k}\right)}{(l+k)^{(1+M)}} \|f\|_{L^1} \\
& \lesssim  \frac{1}{V(x_0,d(E,F))^\frac12} \dfrac{\exp\left( -c \frac{d(E,F)^2}{l+k}\right)}{(l+k)^{1+M}} \|f\|_{L^1}
\end{split}\end{equation}
where the first line holds due to Lemma \ref{GTinequalities} since $x_0$, $D_{l}(E)$ and $F$ satisfy \eqref{assumption1GT}.

When $l \geq \frac{d(E,F)^2}{4}$, replace also \eqref{ODDCase6} by
\begin{equation}\begin{split}
\|\Delta^{1+M}P^{k+l-1}f\|_{L^2(D_{l}(E))}  & \lesssim \frac{1}{V(x_0,\sqrt{k+l})^\frac12}\frac{1}{(k+l)^{1+M}} \|f\|_{L^1} \\
   & \lesssim \frac{1}{V(x_0,d(E,F))^\frac12} \dfrac{\exp\left( -c \frac{d(E,F)^2}{l+k}\right)}{(k+l)^{1+M}} \|f\|_{L^1}
  \end{split}\end{equation}
where the second line follows from Lemma \ref{GTinequalities}, since $x_0$, $F$ and $k+l+t$ satisfy \eqref{assumption2GT}, 
and the third line to the fact that $k+l \gtrsim d(E,F)^2$.
\end{dem}

\section{BMO spaces} 

\label{Duality}

\subsection{Dense sets in Hardy spaces}

\begin{lem} \label{M0subsetH1}
Let $M\in \N$ and $\kappa \in \{1,2\}$.

For all $\epsilon \in (0,+\infty)$, we have the following inclusion
$$\mathcal M_0^{M,\epsilon}(\Gamma) \hookrightarrow H^1_{BZ\kappa,M,\infty}(\Gamma)$$
and for all $\phi \in \mathcal M_0^{M,\epsilon}(\Gamma)$,
$$\|\phi\|_{H^1_{BZ\kappa,M,\infty}} \leq C_{M,\epsilon} \|\phi\|_{\mathcal M_0^{M,\epsilon}}.$$
\end{lem}

\begin{dem}
Let $\phi$ in $\mathcal M_0^{M,\epsilon}(\Gamma)$. Then there exists $\varphi \in L^2(\Gamma)$ such that $\phi = \Delta^M \varphi$ and for all $j\geq 1$,
$$\|\varphi\|_{L^2(C_j(B_0))} 2^{j\epsilon} \lesssim \|\phi\|_{\mathcal M_0^{M,\epsilon}}.$$
Observe that 
\begin{equation} \label{ConvPtoP} \varphi(x) = \sum_{y\in \Gamma} a_y \frac{\1_{\{y\}}(x)}{m(y)} \qquad \forall x\in \Gamma \end{equation}
where $a_y = \varphi(y)m(y)$. In order to prove that $\phi \in H^1_{BZ\kappa,M,\infty}$, it suffices to prove 
\begin{enumerate}[(i)]
 \item for every $y\in \Gamma$, $\Delta^M\frac{\1_{\{y\}}}{m(y)}$ is, up to a harmless multiplicative constant, a $(BZ_\kappa,M)$-atom,
 \item $\ds \sum_{y\in \Gamma} |a_y| \lesssim \|\phi\|_{\mathcal M_0^{M,\epsilon}}$,
 \item $\ds \phi = \sum_{y\in \Gamma} a_y  \Delta^M\frac{\1_{\{y\}}}{m(y)}$ where the convergence  holds  in $L^1(\Gamma)$.
\end{enumerate}
It is easy to check that $\Delta^M\frac{\1_{\{y\}}}{m(y)}$ is a $(BZ_1,M)$-atom associated with $s=1$, $(1,\dots,1)$ and the ball $B(y,1)$.
When $\kappa = 2$, notice that
\[\begin{split}
   \Delta^M\frac{\1_{\{y\}}}{m(y)} & = \left(I - (I+(M^2+1)\Delta)^{-1}\right)^M \left(\frac{I+(M^2+1)\Delta}{M^2+1}\right)^M \frac{\1_{\{y\}}}{m(y)}.
  \end{split}\]
Moreover, $\left(\frac{I+(M^2+1)\Delta}{M^2+1}\right)^M \frac{\1_{\{y\}}}{m(y)}$ is supported in $B(y,M+1)$ and
\[\begin{split}
   \left\| \left(\frac{I+(M^2+1)\Delta}{M^2+1}\right)^M \frac{\1_{\{y\}}}{m(y)} \right\|_{L^2} & \leq \left(\frac{2M^2+3}{M^2+1}\right)^M  \left\| \frac{\1_{\{y\}}}{m(y)} \right\|_{L^2} \\
& \lesssim \frac{1}{m(y)^{\frac{1}{2}}} \\
& \lesssim \frac{1}{V(y,M+1)^{\frac{1}{2}}}.
  \end{split}\]
For point (ii), remark that
\[\begin{split}
 \sum_{y\in \Gamma} |a_y| & = \sum_{j\geq 1} \sum_{y\in C_j(B_0)} |a_y| \\
& \leq \sum_{j\geq 1}  \left( \sum_{y\in C_j(B_0)} \frac{|a_y|^2}{m(y)}\right)^{\frac12}  \left( m(C_j(B_0))\right)^{\frac12} \\
& \lesssim \sum_{j\geq 1}  V(2^jB_0)^{\frac{1}{2}}  \left( \sum_{y\in C_j(B_0)} |\varphi(y)|^2 m(y)\right)^{\frac12} \\
& \qquad = \sum_{j\geq 1}  V(2^jB_0)^{\frac{1}{2}} \|\varphi\|_{L^2(C_j(B_0))} \\
& \quad \leq \sum_{j\geq 1} 2^{-j\epsilon} \|\phi\|_{\mathcal M_0^{M,\epsilon}} \\
& \lesssim \|\phi\|_{\mathcal M_0^{M,\epsilon}}.
  \end{split}\]
For point (iii), notice that (ii) implies the $L^1$-convergence in \eqref{ConvPtoP}. 
The result is then a consequence of the $L^1$-boundedness of $\Delta$. 
\end{dem}

\begin{lem} \label{L2loc} Let $M\in \N^*$ and let $B \subset \Gamma$ be a ball. For all $s\in \N^*$, define $A_s$ as either $(I-P^{s_1})\dots(I-P^{s_M})$ with $(s_1,\dots,s_M)\in \bb 1,2s\bn^M$, or $(I-(I+s\Delta)^{-1})^M$.
If $\varphi \in L^2(B)$ then, for all $s\in \N^*$, $\epsilon>0$ and $M\in \N^*$, $A_s\varphi \in \mathcal M_0^{M,\epsilon}(\Gamma)$.

As a consequence, if $f \in \mathcal E_M$ for some $M\in \N$, then for all $s\in \N$ we can define $A_s f$ as a linear form on finitely supported functions and
$$A_s f \in L^2_{loc}(\Gamma)$$
\end{lem}

\begin{rmq}
 In the case of graphs, a distribution $g$ is in $L^2_{loc}(\Gamma)$ means that we can write $g(x)$ for all $x\in \Gamma$, that is $g$ is a function.
On the contrary, notice that each function on $\Gamma$ belongs to $L^2_{loc}(\Gamma)$ and we use then the notation $L^2_{loc}(\Gamma)$ only by analogy to the case of continuous spaces. 
\end{rmq}

\begin{dem} Fix $\epsilon >0$ and let $\varphi \in L^2(B)$ for some ball $B$ and $k\in \N$ such that $B \subset 2^{k+2}B_0$. The uniform $L^2$-boundedness of $A_s(s\Delta)^{-M}$ yields
$$\sup_{j\in \bb 1,k+1\bn} 2^{j\epsilon} V(2^jB_0)^{\frac{1}{2}} \left\|A_s\Delta^{-M} \varphi \right\|_{L^2(C_j(B_0))} \lesssim s^M 2^{k\epsilon} V(2^kB_0)^{\frac{1}{2}}\|\varphi\|_{L^2(B)}.$$
Moreover,  Proposition \ref{GaffneyEstimatesResolvant} implies, for $j \geq k+2$
\[\begin{split}
2^{j\epsilon} V(2^jB_0)^\frac12 \left\|A_s\Delta^{-M} \varphi \right\|_{L^2(C_j(B_0))} & \lesssim s^M 2^{j\epsilon} V(2^jB_0)^\frac12 e^{-c\frac{2^j}{\sqrt{s}}} \|\varphi\|_{L^2(B)} \\
&  \lesssim s^{M+\frac{\epsilon}{2} + \frac{d_0}{4}+1} V(B_0)^{\frac{1}{2}}\|\varphi\|_{L^2(B)} 
  \end{split}\]
where $d_0$ is given by Proposition \ref{propDV}.
One concludes that $A_s \varphi \in \mathcal M_0^{M,\epsilon}(\Gamma)$ and
\begin{equation} \label{SoUsefull} \left\|A_s \varphi \right\|_{\mathcal M_0^{M,\epsilon}} \lesssim  s^{M+\frac{\epsilon}{2} + \frac{d_0}{4}+1} 2^{k\epsilon} V(2^{k}B_0)^{\frac{1}{2}} \|\varphi\|_{L^2(B)}. \end{equation}

\bs

Let us prove the second claim of the lemma. Let $\epsilon$ such that $f \in (\mathcal M_0^{M,\epsilon}(\Gamma))^*$.
For all balls $B$ and all functions $\varphi$ supported in $B$, one has
\[\begin{split}
   \left|\left<A_s f,\varphi \right>\right| & := \left|\left<f,A_s \varphi \right>\right| \\
& \lesssim \|f\|_{(\mathcal M_0^{M,\epsilon})^*}  \left\|A_s \varphi \right\|_{\mathcal M_0^{M,\epsilon}} \\
& \lesssim \|f\|_{(\mathcal M_0^{M,\epsilon})^*} \|\varphi\|_{L^2(B)},
  \end{split}\]
which proves the lemma since the estimate works for any ball $B$ and any $\varphi \in L^2(B)$.
\end{dem}

\begin{defi}
 Let $\kappa \in \{1,2\}$ and $M\in \N^*$. Define $\mathbb H^1_{BZ\kappa,M,\epsilon}(\Gamma)$ as the subset of $H^1_{BZ\kappa,M,\epsilon}(\Gamma)$ made of the functions $g$ that can be written as $g = \sum_{i=0}^N \lambda_i a_i$ where $\lambda_i \in \R$ and $a_i$ is a  $(BZ_\kappa,M,\epsilon)$-molecule and 
$$\sum_{i=0}^N |\lambda_i| \lesssim 2 \|g\|_{H^1_{A,\epsilon}}.$$
\end{defi}

\begin{lem} \label{DensityFiniteSum}
 For $\kappa \in \{1,2\}$ and $M\in \N^*$, the set $\mathbb H^1_{BZ\kappa,M,\epsilon}(\Gamma)$ is dense in $H^1_{BZ\kappa,M,\epsilon}(\Gamma)$.
\end{lem}

\begin{rmq}
This lemma is  identical to  Lemma 4.5 in \cite{BZ}. However,  we present here a different proof. \end{rmq}

\begin{dem}
Let $\kappa \in \{1,2\}$ and $M\in \N^*$.

Let $f\in H^1_{BZ\kappa,M,\epsilon}(\Gamma)$. There exist a numerical sequence $(\lambda_i)_{i\in \N} \in \ell^1(\N)$ and a sequence $(a_i)_{i\in \N}$ of $(BZ_\kappa,M,\epsilon)$-molecules such that $f = \sum \lambda_i a_i$ and
$$\sum_{i\in \N} |\lambda_i| \leq \frac{3}{2} \|f\|_{H^1_{BZ\kappa,M,\epsilon}}.$$
Let $\eta\in \left(0,\frac{1}{4} \right)$. There exists $N \in \N$ such that $\sum_{i> N} |\lambda_i| \leq \eta \|f\|_{H^1_{BZ\kappa,M,\epsilon}}$. 
We set $g = \ds \sum_{i=0}^N \lambda_i a_i$. Then
\[\begin{split}
\|f-g\|_{H^1_{BZ\kappa,M,\epsilon}} & = \left\| \sum_{i>N} \lambda_i a_i \right\|_{H^1_{BZ\kappa,M,\epsilon}}  \\
& \leq \sum_{i>N} |\lambda_i| \leq \eta \|f\|_{H^1_{BZ\kappa,M,\epsilon}}
  \end{split}\]
and, therefore $\|f\|_{H^1_{BZ\kappa,M,\epsilon}} \leq \|g\|_{H^1_{BZ\kappa,M,\epsilon}} +  \eta \|f\|_{H^1_{BZ\kappa,M,\epsilon}}$, which implies
\[\begin{split}
\sum_{i=0}^N |\lambda_i|  & \leq \frac{3}{2}\|f\|_{H^1_{BZ\kappa,M,\epsilon}} \\
& \leq \frac{3}{2(1-\eta)} \|g\|_{H^1_{BZ\kappa,M,\epsilon}} \\
& \leq  2\|g\|_{H^1_{BZ\kappa,M,\epsilon}}.
  \end{split}\]
\end{dem}

\begin{lem} \label{HH1inM0} Let $\kappa \in \{1,2\}$  and $M \in \N$. Let $0<\epsilon < \bar\epsilon \leq +\infty$. Then $\mathbb H^1_{BZ\kappa,M,\bar\epsilon} (\Gamma) \subset \mathcal M_0^{M,\epsilon}$.
\end{lem}

\begin{dem} Since $\mathcal M_0^{M,\epsilon}$ is a vector space, it is enough to prove that for each $(BZ_\kappa,M,\bar\epsilon)$-molecule $a$, one has $a \in \mathcal M_0^{M,\epsilon}$.

Notice that the case $\bar\epsilon = \infty$ is proven in Lemma \ref{L2loc}. 
Let $\bar \epsilon < +\infty$ and $a= A_s b$ be a $(BZ_\kappa,M,\bar\epsilon)$-molecule associated with $s\in \N^*$ and the ball $B$ of radius $\sqrt{s}$.
For all $j\geq 1$, Corollary \ref{corVitali2} provides a covering  of $C_j(B)$ with balls of radius $\sqrt{s}$ and with bounded overlapping. 
We  label these balls as $(B_i)_{i\in I_j}$. . 
Consequently,
\[\begin{split}
  \|a\|_{\mathcal M_0^{M,\epsilon}} & \leq \sum_{j\geq 1} \|A_s (b\1_{C_j(B)})\|_{\mathcal M_0^{M,\epsilon}} \\
& \leq \sum_{j\geq 1} \sum_{i\in I_j}  \|A_s (b\1_{B_i})\|_{\mathcal M_0^{M,\epsilon}}. 
  \end{split}\]
Moreover, $d(B_i,B_0) \lesssim 2^{j+k}$ where $k$ is such that $B \subset 2^{k+2}B_0$. Thus Lemma \ref{L2loc} implies
\[\begin{split}
  \|a\|_{\mathcal M_0^{M,\epsilon}} & \leq  C_s \sum_{j\geq 1} \sum_{i\in I_j} 2^{(j+k)\epsilon} V(2^{j+k}B_0)^{\frac{1}{2}} \|b\|_{L^2(B_i)} \\
& \leq C_s \sum_{j\geq 1} 2^{(j+k)\epsilon} V(2^{j+k}B_0)^{\frac{1}{2}} \|b\|_{L^2(\tilde C_j(B))} \\
& \leq C_s \sum_{j\geq 1} \frac{2^{(j+k)\epsilon}}{2^{j\bar\epsilon}} \left(\frac{V(2^{j+k}B_0)}{V(2^{j-1}B)}\right)^{\frac{1}{2}} \\
& \leq C_s 2^{k(\epsilon+\frac{d}{2})} \sum_{j\geq 1} 2^{j(\epsilon-\bar\epsilon)} \\
& < +\infty 
  \end{split}\]
where $\tilde C_j$ denote $C_j(B) \cup C_{j-1}(B) \cup C_{j+1}(B)$, and where we use the definition of a $(BZ_\kappa,M,\bar\epsilon)$-molecule for the third line and the fact that $2^{j+k}B_0 \subset 2^{j+k+2}B$.
\end{dem}

\subsection{Inclusions between $BMO$ spaces} 

\begin{lem} \label{lemBMO2inBMO1}
There exists $C>0$ such that for all $s\in \N^*$, all $M$-tuples $(s_1,\dots,s_M)\in \bb s,2s\bn^M$, all balls $B$ of radius $\sqrt{s}$ and all functions $f\in BMO_{BZ2,M}$, one has
$$ \|(I-P^{s_1})\dots(I-P^{s_M}) f\|_{L^2(B)} \leq C V(B)^{\frac{1}{2}} \|f\|_{ BMO_{BZ2,M} }.$$
\end{lem}

\begin{dem} For $s\in \N^*$, the operator $Q_s$ stands for
$$Q_s := \frac{1}{s} \sum_{k=0}^{s-1} P^k = (I-P^s)(s\Delta)^{-1}.$$
For all $s \in \N^*$, all $s_0\in \bb s,2s\bn$ and all $f\in \mathcal E_0$, one has
\[\begin{split}
   (I-P^{s_0}) f  &  =  (I-P^{s_0}) (I + s\Delta)(I+s\Delta)^{-1}f  \\
& =  \left(\frac{1}{s} \sum_{k=0}^{s_0-1} P^k f \right) (I+s\Delta) s\Delta(I+s\Delta)^{-1} f \\
& = \left[ \frac{s_0}{s} Q_{s_0} + (I-P^{s_0}) \right] (I-(I+s\Delta)^{-1})  f
  \end{split}\]
Recall that all terms make sense and are in $L^2_{loc}(\Gamma)$, according to Lemma \ref{L2loc}.
As a consequence, for $(s_1,\dots,s_M)\in \bb s,2s\bn^M$, one has
\begin{equation}\label{molBZ1aremolBZ2}\begin{split}
   (I-P^{s_1})\dots(I-P^{s_M}) f & = \prod_{i=1}^M  \left[ \frac{s_i}{s} Q_{s_i} + (I-P^{s_i}) \right] (I-(I+s\Delta)^{-1})^M f  \end{split} \end{equation}
  Since $\frac{s_i}{s} \leq 2$, Proposition \ref{GaffneyEstimatesResolvant} yields that $\ds  \prod_{i=1}^M  \left[ \frac{s_i}{s} Q_{s_i} + (I-P^{s_i}) \right]$ satisfies Gaffney-Davies estimates. Hence,
 \[ \begin{split}
    \left\|(I-P^{s_1})\dots(I-P^{s_M}) f \right\|_{L^2(B)} & 
   \leq \sum_{j\geq 1} \left\|\prod_{i=1}^M  \left[ \frac{s_i}{s} Q_{s_i} + (I-P^{s_i}) \right] \left[ \1_{C_j(B)} (I-(I+s\Delta)^{-1})^M f \right] \right\|_{L^2(B)} \\
& \lesssim \sum_{j\geq 1} e^{-c4^j} \| (I-(I+s\Delta)^{-1})^M f\|_{L^2(C_j(B))} \\
& \lesssim \sum_{j\geq 1} e^{-c4^j} \|(I-(I+s\Delta)^{-1})^M f\|_{L^2(2^{j+1}B)} \\
& \lesssim \sum_{j\geq 1} e^{-c4^j} V(2^{j+1}B)^{\frac{1}{2}} \|f\|_{BMO_{BZ2,M}} \\
& \lesssim V(B)^{\frac{1}{2}} \|f\|_{ BMO_{BZ2,M}}
\end{split}\]
where the last line holds thanks to Proposition \ref{propDV}.
\end{dem}

\begin{cor} \label{Main22}
Let $M \in \N^*$. Then $BMO_{BZ2,M}(\Gamma) \subset BMO_{BZ1,M}(\Gamma)$. More precisely, for all $f\in BMO_{BZ2,M}(\Gamma)$,
$$\|f\|_{BMO_{BZ1,M}} \lesssim \|f\|_{BMO_{BZ2,M}}.$$
\end{cor}

\begin{dem}
Immediate consequence of Lemma \ref{lemBMO2inBMO1}. 
\end{dem}

We want now to prove the  converse  inclusion, that is $BMO_{BZ1,M}(\Gamma) \subset BMO_{BZ2,M}(\Gamma)$. 
We begin  with  the next proposition, inspired from Proposition 2.6 in \cite{DY1}.

\begin{prop} \label{DYprop}
 Let $M\in \N^*$. There exists $C>0$ only depending on $\Gamma$ and $M$ such that for all $f\in BMO_{BZ1,M}(\Gamma)$, for all balls  $B = B(x_0,\sqrt s)$  and all integers $(a,b_1,\dots,b_M)\in \N \times \bb 0,2s\bn^M$,
$$\|P^a(I-P^{b_1})\dots(I-P^{b_M}) f\|_{L^2(B)} \leq C a_s^{\frac{d_0+1}{2}} V(B)^{\frac{1}{2}} \|f\|_{BMO_{BZ1,M}}$$
where $a_s = \max\left\{1,\frac{a}{s}\right\}$.
\end{prop}

\begin{rmq}
We can replace $a_s^{d_0+1}$ by $a_s^{d_0+\epsilon}$ with $\epsilon>0$ in the conclusion of the Proposition \ref{DYprop} (in this case, $C$ depends on $\epsilon$).
\end{rmq}

\begin{dem} (Proposition \ref{DYprop})
\begin{enumerate}[(1)]
\item Let us prove the proposition when $\ds s \leq \min_{i\in \bb1,M\bn}b_i$. The case where $a=0$ is a consequence of the definition of $BMO_{BZ1,M}$ and will therefore be skipped. 
Let $(B_i)_{i\in I_j}$ be the covering of $C_j(B)$ provided by Corollary \ref{corVitali2}. Then,
\begin{equation} \label{DYdem}
 \begin{split}
   \|P^a & (I-P^{b_1})\dots(I-P^{b_M})f\|_{L^2(B)} \\
   & \lesssim \|(I-P^{b_1})\dots(I-P^{b_M})f\|_{L^2(C_1(B))} +  \sum_{j\geq 2} \exp\left(-c\frac{4^jb}{a}\right) \|(I-P^{b_1})\dots(I-P^{b_M}) f\|_{L^2(C_j(B))} \\
& \quad \leq \|(I-P^{b_1})\dots(I-P^{b_M}) f\|_{L^2(4B)} + \sum_{j\geq 2} \exp\left(-c\frac{4^jb}{a}\right) \|(I-P^{b_1})\dots(I-P^{b_M})f\|_{L^2(2^{j+1}B)} \\
& \lesssim  V(4B)^{\frac{1}{2}} \|f\|_{BMO_{BZ1,M}} + \sum_{j\geq 2} \sum_{i\in I_{j}}\exp\left(-c\frac{4^jb}{a}\right) \|(I-P^{b_1})\dots(I-P^{b_M}) f\|_{L^2(B_i)} \\
& \lesssim V(B)^{\frac{1}{2}} \|f\|_{BMO_{BZ1,M}} \left[ 1 + \sum_{j\geq 2} 2^{j{d_0+1}} \exp\left(-c\frac{4^jb}{a}\right) \right] \\
& \lesssim  V(B)^{\frac{1}{2}} \|f\|_{BMO_{BZ1,M}}  a_s^{\frac{d_0+1}{2}}  
\end{split}\end{equation}
where we use the Davies-Gaffney estimates for the first line and the doubling property for the last but one line.

\item General case. For each $b_i <s$, write
$$(I-P^{b_i}) = (I-P^{2s}) - P^{b_i}(I-P^{2s-b_i}).$$
Hence, $P^a(I-P^{b_1})\dots(I-P^{b_M})$ can be written as a sum of terms
$$P^{\tilde a}(I-P^{\tilde b_1})\dots(I-P^{\tilde b_M})$$
where $\tilde b_i \in \bb s,2s\bn $ and $\tilde a \in \bb a, a+Ms\bn$. 
The general case can be then deduced from the previous case.
\end{enumerate}
\end{dem}

\begin{prop} \label{lemBMO1inBMO2}
Let $M\in \N^*$. There exists $C>0$ such that for all balls $B$ of radius $\sqrt{s}$, all integers $b\in \bb 0,2s\bn$ and all $f\in BMO_{BZ1,M}$, one has
$$ \|(I-(I+b\Delta)^{-1})^M f\|_{L^2(B)} \leq C V(B)^{\frac{1}{2}} \|f\|_{BMO_{BZ1,M}}.$$
\end{prop}

\begin{dem} 
Let $\varphi \in L^2(\Gamma)$ supported in $B$. Recall that Lemma \ref{L2loc} states that $\varphi$, $(I-P^{k_1})\dots(I-P^{k_M})\varphi$ and $(I-(I+b\Delta)^{-1})^M \varphi$ are in $\mathcal M_0^{M,\epsilon}$ for all $\epsilon >0$.
Moreover, for all $b \in \N$, one has
\[\begin{split}
   (I+b\Delta)^{-1}\varphi  &  = (1+b)^{-1} \left(I-\frac{b}{1+b}P \right)^{-1}\varphi \\
& =  \sum_{k=1}^{+\infty} \left(\frac{1}{1+b}\right) \left(\frac{b}{1+b}\right)^k P^k \varphi
  \end{split}\]
where the convergence  holds  in $L^2(\Gamma)$. Consequently,
\[\begin{split}
   (I-(I+b\Delta)^{-1})\varphi  & = \sum_{k=1}^{+\infty} \left(\frac{1}{1+b}\right) \left(\frac{b}{1+b}\right)^k (I- P^k) \varphi
  \end{split}\]
and thus,
\[\begin{split}
   (I-(I+b\Delta)^{-1})^M\varphi  & = \left(\frac{1}{1+b}\right)^M \sum_{k=1}^{+\infty}  \left(\frac{b}{1+b}\right)^k \sum_{k_1+\dots+k_M = k} (I-P^{k_1})\dots(I-P^{k_M}) \varphi
  \end{split}\]
where the convergence  still holds  in $L^2(\Gamma)$.

 In order to prove that the convergence  holds  in $\mathcal M_0^{M,\epsilon}$ for all $\epsilon>0$, it suffices to show that
$$S := \left(\frac{1}{1+b}\right)^M \sum_{k=1}^{+\infty}  \left(\frac{b}{1+b}\right)^k \sum_{k_1+\dots+k_M = k} \|(I-P^{k_1})\dots(I-P^{k_M}) \varphi\|_{\mathcal M_0^{M,\epsilon}} < +\infty.$$
Indeed, according to \eqref{SoUsefull}, one has
\[\begin{split}
   S & \lesssim \left(\frac{1}{1+b}\right)^M \sum_{k=1}^{+\infty}  \left(\frac{b}{1+b}\right)^k \sum_{k_1+\dots+k_M = k} k^{\frac{\epsilon}{2}+ \frac{d}{4}+1} \|\varphi\|_{L^2(B)} \\
   & \lesssim  \sum_{k=1}^{+\infty} \frac{(k+1)^{M + \frac{\epsilon}{2}+ \frac{d}{4}}}{(1+b)^M}  \left(\frac{b}{1+b}\right)^k \|\varphi\|_{L^2(B)} \\
& \lesssim b^{\frac{\epsilon}{2} + \frac{d}{4} + 2}\sum_{k=1}^{+\infty} \frac{1}{(1+k)^2}  \|\varphi\|_{L^2(B)} \\
& < +\infty
  \end{split}\]
where the  third  line comes from Lemma \ref{ExpDecay}.

For $f\in \mathcal E_M$, there exists $\epsilon>0$ such that $(\mathcal M_0^{M,\epsilon})^*$. Moreover, Lemma \ref{L2loc} states that $(I-(I+s\Delta)^{-1})^Mf$ and $(I-P^{k_1})\dots(I-P^{k_M})f$ (for all $(k_1,\dots,k_M)\in \N^M$) are in $L^2_{loc}(\Gamma)$. As a consequence,
\[\begin{split}
   \|(I-(I+b\Delta)^{-1})^M f\|_{L^2(B)} & = \sup_{\begin{subarray}{c} \|\varphi\|_2=1 \\ \Supp \varphi \subset B \end{subarray}} |\left< f,(I-(I+b\Delta)^{-1})^M \varphi \right>| \\
& \leq \left(\frac{1}{1+b}\right)^M \sum_{k=1}^{+\infty}  \left(\frac{b}{1+b}\right)^k \sum_{k_1+\dots+k_M = k} \sup_{\begin{subarray}{c} \|\varphi\|_2=1 \\ \Supp \varphi \subset B  \end{subarray}} |\left< f, (I-P^{k_1})\dots(I-P^{k_M}) \varphi \right>| \\
& \qquad = \left(\frac{1}{1+b}\right)^M \sum_{k=1}^{+\infty}  \left(\frac{b}{1+b}\right)^k \sum_{k_1+\dots+k_M = k} \|(I-P^{k_1})\dots(I-P^{k_M}) f\|_{L^2(B)}
  \end{split}\]
where the pairing is between $\mathcal M_0^{M,\epsilon}$ and its dual. Therefore
\[\begin{split}
   \|(I-(I+b\Delta)^{-1})^Mf\|_{L^2(B)} & \lesssim \left(\frac{1}{1+b}\right)^M \sum_{k=1}^{+\infty}  \left(\frac{b}{1+b}\right)^k \sum_{k_1+\dots+k_M = k} \|(I-P^{k_1})\dots(I-P^{k_M}) f\|_{L^2(B)} \\
& \quad \leq  \left(\frac{1}{1+b}\right)^M \sum_{k=1}^{b}  \left(\frac{b}{1+b}\right)^k \sum_{k_1+\dots+k_M = k} \|(I-P^{k_1})\dots(I-P^{k_M}) f\|_{L^2(B)} \\
& \qquad +   \left(\frac{1}{1+b}\right)^M \sum_{k=b+1}^{\infty}  \left(\frac{b}{1+b}\right)^k \sum_{k_1+\dots+k_M = k} \|(I-P^{k_1})\dots(I-P^{k_M}) f\|_{L^2(B)} \\
& \qquad : =  I_1 + I_2. 
  \end{split}\]
We estimate the first term with Proposition \ref{DYprop} and Lemma \ref{ExpDecay}:
\[\begin{split}
   I_1 & \lesssim  \sum_{k=1}^{b} \frac{(1+k)^{M-1}}{(1+b)^M} \left(\frac{b}{1+b}\right)^k \|f\|_{BMO_{BZ1,M}} V(B)^{\frac{1}{2}} \\
& \lesssim (1+b)^{-1} \sum_{k=0}^{b-1}  \|f\|_{BMO_{BZ1,M}} V(B)^{\frac{1}{2}} \\
& \lesssim \|f\|_{BMO_{BZ1,M}} V(B)^{\frac{1}{2}}.
  \end{split}\]
We turn now to the estimate of the second term. One has,  using Proposition  \ref{DYprop} and  Lemma \ref{ExpDecay} again, 
\[ \begin{split}
    I_2 & \lesssim \left(\frac{1}{1+b}\right)^M \sum_{k=b+1}^{\infty}  \left(\frac{b}{1+b}\right)^k \sum_{k_1+\dots+k_M = k} \|(I-P^{k_1})\dots(I-P^{k_M}) f\|_{L^2(\sqrt{\frac{k}{b}}B)}  \\
    & \lesssim  \sum_{k=b+1}^{\infty} \frac{1}{1+k} \left(\frac{1+k}{1+b}\right)^M \left(\frac{b}{1+b}\right)^k  \|f\|_{BMO_{BZ1,M}} V\left(\sqrt{\frac{k}{b}}B\right)^{\frac{1}{2}} \\
& \lesssim\sum_{k=b+1}^{\infty} \frac{1}{1+k} \left(\frac{1+k}{1+b}\right)^{M+\frac{d_0}{2}+1} \left(\frac{b}{1+b}\right)^k  \|f\|_{BMO_{BZ1,M}} V(B)^{\frac{1}{2}} \\
& \lesssim\sum_{k=b+1}^{\infty} \frac{1}{1+k} \left(\frac{1+k}{1+b}\right)^{-1}  \|f\|_{BMO_{BZ1,M}} V(B)^{\frac{1}{2}} \\
& \lesssim \|f\|_{BMO_{BZ1,M}} V(B)^{\frac{1}{2}}, \end{split}\]
where we used Proposition \ref{propDV} for the third line. 
\end{dem}

\begin{cor} \label{Main21}
Let $M \in \N$. Then $BMO_{BZ1,M}(\Gamma) \subset BMO_{BZ2,M}(\Gamma)$. More precisely, for all $f\in BMO_{BZ1,M}(\Gamma)$,
$$\|f\|_{BMO_{BZ2,M}} \lesssim \|f\|_{BMO_{BZ1,M}}.$$
\end{cor}

\begin{dem}
Immediate consequence of Proposition \ref{lemBMO1inBMO2}.
\end{dem}

\subsection{Duals of Hardy spaces}

\begin{prop} \label{Main11}
Let $\kappa \in \{1,2\}$ and $M \in \N^*$. 

Let $\ell$ be a bounded linear functional on $H^1_{BZ\kappa,M,\infty}(\Gamma)$.
Then $\ell$ actually belongs to $BMO_{BZ\kappa,M}(\Gamma) \cap \mathcal F_M$ and for all $g\in \mathbb H^1_{BZ\kappa,M,\infty}(\Gamma)$, there holds
\begin{equation} \label{dualityM0} \ell(g) = \left< \ell, g\right> \end{equation}
where the pairing is between $\mathcal M_0^{M,\epsilon}(\Gamma)$ and its dual. Moreover,
$$\|\ell\|_{BMO_{BZ\kappa,M}} \lesssim \|\ell\|_{(H^1_{BZ\kappa,M,\infty})^*}$$
\end{prop}

\begin{dem} Let $\kappa \in \{1,2\}$ and $M \in \N$.

Let $\ell$ in $\left[H^1_{BZ\kappa,M,\infty}(\Gamma)\right]^*$. 
According to Lemma \ref{M0subsetH1}, $\ell \in \bigcap_{\epsilon >0}\left[\mathcal M_0^{M,\epsilon}\right]^* = \mathcal F_M$. 
The following two claims
\begin{enumerate}[(i)]
 \item $\mathbb H^1_{BZ\kappa,M,\infty}(\Gamma) \subset \mathcal M_0^{M,\epsilon}$, \item $\mathbb H^1_{BZ\kappa,M,\infty}(\Gamma)$ is dense in $H^1_{BZ\kappa,M,\infty}(\Gamma)$,
\end{enumerate}
are respectively a consequence of Lemma \ref{HH1inM0} and of Lemma \ref{DensityFiniteSum}.  They imply that \eqref{dualityM0} makes sense and uniquely describes $\ell$. 

It remains to check the last claim, that is
$$\|\ell\|_{BMO_A} \lesssim \|\ell\|_{(H^1_{BZ\kappa,M,\infty})^*}$$
Fix $s\in \N^*$, a $M$-tuple $(s_1,\dots,s_M) \in \bb s,2s \bn^M$, and a ball $B$ of radius $\sqrt{s}$. We wrote $A_s$ for $(I-P^{s_1})\dots(I-P^{s_M})$ if $\kappa = 1$ and for $(I-(I+s\Delta)^{-1})^M$ if $\kappa = 2$.

Let $\varphi \in L^2(B)$ with norm 1. Then 
$$\frac{1}{V(B)^{\frac{1}{2}}} A_s\varphi$$
is a $(BZ_\kappa,M)$-atom. Thus,
$$\left\|\frac{1}{V(B)^\frac12} A_s \varphi \right\|_{H^1_{BZ\kappa,M,\infty}} \leq 1,$$
i.e.,
\[\begin{split}
   \frac{1}{V(B)^{\frac{1}{2}}} |\left< A_s\ell,\varphi \right>| & = \frac{1}{V(B)^{\frac{1}{2}}} |\left< \ell, A_s\varphi \right>| \\
& \lesssim \|\ell\|_{(H^1_{BZ\kappa,M,\infty})^*}.
  \end{split}\]
Lemma \ref{L2loc} provides that $A_s \ell \in L^2_{loc}(\Gamma)$. Taking the supremum over all $\varphi$ supported in $B$, we obtain
$$\left(\frac{1}{V(B)} \sum_{x\in B} |A_s \ell(x)|^2 m(x)\right)^{\frac{1}{2}} \lesssim \|\ell\|_{(H^1_{BZ\kappa,M,\infty})^*}.$$
Finally, taking the supremum over all $s\in \N^*$, all $M$-tuples $(s_1,\dots,s_M) \in \bb s,2s \bn^M$ and all balls $B$ of radius $\sqrt{s}$ leads us to the result.
\end{dem}

\begin{prop} \label{Main12}
Let $\kappa \in \{1,2\}$ and $M \in \N^*$. 

Let $\epsilon >0$ and $f \in BMO_{BZ\kappa,M}(\Gamma) \cap \mathcal F_M$. The linear functional given by
$$\ell(g) : = \left<f,g \right>$$
initially defined on $\mathbb H^1_{BZ\kappa,M,2\epsilon}(\Gamma)$, 
and where the pairing is between $\mathcal M_0^{M,\epsilon}$ and its dual,
 has a unique bounded extension to $H^1_{BZ\kappa,M,2\epsilon}(\Gamma)$ with
$$\|\ell \|_{(H^1_{BZ\kappa,M,2\epsilon})^*} \lesssim \|f\|_{ BMO_{BZ\kappa,M}(\Gamma)}.$$
\end{prop}

\begin{dem} Let $\kappa \in \{1,2\}$ and $M \in \N^*$. In the proof, $A_s$ will denote $(I-P^{s_1})\dots(I-P^{s_M})$ (for some $(s_1,\dots,s_M)\in \bb s,2s\bn^M$) or $(I-(I+s\Delta)^{-1})^M$, depending  whether  $\kappa$ is equal to 1 or 2.

Let us prove that for every $(BZ_\kappa,M,2\epsilon)$-molecule $a$, one has
\begin{equation}
|\left< f,a\right>| \lesssim \|f\|_{BMO_{BZ\kappa,M}}.
\end{equation}
Since $f\in \mathcal F_M$, then $f\in \left(\mathcal M_0^{M,\epsilon}\right)^*$. 
In particular, Lemma \ref{L2loc} provides that $A_s f \in L^2_{loc}(\Gamma)$.
Thus, if $a = A_s b$ is a $(BZ_\kappa,M,2\epsilon)$-molecule associated with a ball $B$ of radius $\sqrt{s}$, we may write
\[\begin{split}
   |\left<f,a \right>| & = \left| \sum_{x\in \Gamma}  A_sf(x) b(x) m(x) \right| \\
& \leq  \sum_{j\geq 1} \|A_sf\|_{L^2(C_j(B))} \|b\|_{L^2(C_j(B))} \\
& \leq  \sum_{j\geq 1} 2^{-2j\epsilon} V(2^jB)^{-\frac{1}{2}} \|A_sf\|_{L^2(2^{j+1}B)} \\ 
& \lesssim \sum_{j\geq 1} 2^{-2j\epsilon} V(2^jB)^{-\frac{1}{2}} V(2^{j+1}B)^{\frac12} \|f\|_{BMO_{BZ\kappa,M}} \\
& \lesssim \|f\|_{BMO_{BZ\kappa,M}},
  \end{split}\]
  where we used for the last but one line Proposition \ref{DYprop} (if $\kappa =1$) or Proposition \ref{lemBMO2inBMO1} and Corollary \ref{Main22} (if $\kappa =2$).

Our next step is to show that for every $g\in \mathbb H^1_{BZ\kappa,M,2\epsilon}$, we have
$$|\left< f,g\right>| \lesssim \|g\|_{H^1_{BZ\kappa,M,2\epsilon}} \|f\|_{BMO_{BZ\kappa,M}}.$$
Indeed, let $N \in \N$,  $(\lambda_i)_i\in \bb 0,N\bn \in \R^N$ and $(a_i = A_{s_i}b_i)_{i \in \bb0,N\bn}$ a sequence of $(BZ_\kappa,M,2\epsilon)$-molecules 
that satisfies $g = \sum \lambda_i a_i$ and $\sum |\lambda_i| \lesssim 2 \|g\|_{H^1_{BZ\kappa,M,2\epsilon}}$, then
\[\begin{split}
   \left|\ell\left( g\right)\right| & \leq \sum_{i=0}^N \left|\lambda_i\right| \left| \ell(a_i)\right| \\
& \lesssim  \|f\|_{BMO_{BZ\kappa,M}} \sum_{i=0}^N \left|\lambda_i\right| \\
& \lesssim \|f\|_{BMO_{BZ\kappa,M}} \|g\|_{H^1_{BZ\kappa,M,2\epsilon}}.
  \end{split}\]
Since $\mathbb H^1_{BZ\kappa,M,2\epsilon}$ is dense in $H^1_{BZ\kappa,M,2\epsilon}$, $\ell$ has an unique bounded extension that satisfies
$$\|\ell\|_{(H^1_{BZ\kappa,M,2\epsilon})^*} \lesssim  \|f\|_{BMO_{BZ\kappa,M}}.$$
\end{dem}

\begin{prop} \label{Main13}
Let $\kappa \in \{1,2\}$ and $M \in \N^*$.

Let $f \in BMO_{BZ\kappa,M}(\Gamma)$ and let $\epsilon>0$ such that $f\in (\mathcal M_0^{M,\epsilon}(\Gamma))^*$.
The linear functional given by
$$\ell(g) : = \left<f,g \right>$$
initially defined on $\mathbb H^1_{BZ\kappa,M,\infty}(\Gamma)$ which is a dense subset of $\mathcal M_0^{M,\epsilon}$, 
and where the pairing is that between $\mathcal M_0^{M,\epsilon}$ and its dual,
 has a unique extension to $H^1_{BZ\kappa,M,\infty}(\Gamma)$ with
$$\|\ell \|_{(H^1_{BZ\kappa,M,\infty})^*} \lesssim \|f\|_{BMO_{BZ\kappa,M}}.$$
\end{prop}

\begin{dem}
Same proof than Proposition \ref{Main12} with obvious modifications. 
The only difference is: in Proposition \ref{Main12}, $\epsilon>0$ is given by the Hardy space $H^1_{BZ\kappa,M,2\epsilon}$ and in Proposition \ref{Main13}, $\epsilon>0$ is given by the  functional  $f \in \mathcal E_M$. 
\end{dem}

We turn now to the proof of Theorem \ref{Main1}.

\begin{dem} Let $\kappa \in \{1,2\}$ and $M\in \N^*$.

Proposition \ref{Main11} and Corollary \ref{Main13} provide the continuous embeddings $$(H^1_{BZ\kappa,M,\infty})^* \hookrightarrow  BMO_{BZ\kappa,M} \cap \mathcal F_M \hookrightarrow BMO_{BZ\kappa,M} \hookrightarrow (H^1_{BZ\kappa,M,\infty})^*.$$
As a consequence, $BMO_{BZ\kappa,M}$ is the dual space of $H^1_{BZ\kappa,M,\infty}$ and is actually included in $\mathcal F_M$.

Besides, Propositions \ref{Main11} and \ref{Main13} yield, for any $\epsilon >0$
$$(H^1_{BZ\kappa,M,\infty})^* \hookrightarrow  BMO_{BZ\kappa,M} \cap \mathcal F_M \hookrightarrow (H^1_{BZ\kappa,M,\epsilon})^*.$$
Since the inclusion $(H^1_{BZ\kappa,M,\epsilon})^*\hookrightarrow (H^1_{BZ\kappa,M,\infty})^*$ is obvious, we obtain that $BMO_{BZ\kappa,M} \cap \mathcal F_M = BMO_{BZ\kappa,M}$ is also the dual space of $H^1_{BZ\kappa,M,\epsilon}$.

The last claim of the Theorem, that is for a fixed $M\in \N^*$, the spaces $H^1_{BZ\kappa,M,\epsilon}(\Gamma)$ for $\kappa \in \{1,2\}$ and $\epsilon \in (0,+\infty]$ are all equivalent, is only a consequence of the proposition \ref{Rudin} below.
Indeed, for $m\in \N^*$ and $\kappa \in \{1,2\}$, the inclusion $H^1_{BZ\kappa,M,\epsilon} \subset H^{1}_{BZ\kappa,M,\eta}$ when $0<\eta<\epsilon\leq +\infty$ is obvious 
and then Proposition \ref{Rudin} yields the  equality between the spaces $H^1_{BZ\kappa,M,\epsilon}$ for $\epsilon \in (0,+\infty]$, together with the equivalence of norms.  It remains to check that, for example, $H^1_{BZ1,M,\infty} \subset H^1_{BZ2,M,1}$. 
For this, notice first that similarly to \eqref{molBZ1aremolBZ2}, for a $(BZ_1,M)$-atom $a $ associated with $s\in \N^*$, $(s_1,\dots,s_m)\in \bb s,2s\bn^M$, a ball $B$ of radius $\sqrt{s}$ and a function $b\in L^2(B)$, one has
\[\begin{split}
 a & = (I-P^{s_1})\dots(I-P^{s_M}) b \\
& = (I-(I+s\Delta)^{-1})^M \prod_{i=1}^M  \left[ \frac{s_i}{s} Q_{s_i} + (I-P^{s_i}) \right]  b.
  \end{split}\]
We have to check that $\prod_{i=1}^M  \left[ \frac{s_i}{s} Q_{s_i} + (I-P^{s_i}) \right]  b$ satisfies, up to a multiplicative constant, the estimates given by (ii) of the definition of a $(BZ_2,M,1)$-molecule. 
This calculus, which is a straightforward consequence of the Gaffney estimates provided by Proposition \ref{GaffneyEstimatesResolvant}, is left to the reader.
 \end{dem}

\begin{prop} \label{Rudin}
If $(E,\|.\|_E)$ and $(F,\|.\|_F)$ are two Banach spaces with the same dual $(G, \|.\|_G)$ and moreover if we have the continuous inclusion $E \subset F$, then $E=F$ with equivalent norms.
\end{prop}

\begin{dem}
Let $T$ be the linear operator defined by 
$$T: \, e\in E \mapsto e \in F.$$
$T$ is bounded and its adjoint $T^*$ is 
$$T^*: g\in G \mapsto g\in G,$$
that is the identity on $G$. 
Theorem 4.15 in \cite{RudinFA} implies that $E = F$, and then, by the open mapping theorem, we deduce that the norm of $E$ is dominated by the norm of $F$.
\end{dem}

\section{Inclusions between Hardy spaces}

\label{HardyEquiv}

\subsection{$H^1_{BZ1,M,\epsilon} \cap L^2 \subset E^1_{quad,\beta}$: the case of functions}

\begin{prop} \label{Main31}
Let $\epsilon>0$ , $M \in (\frac{d_0}{4},+\infty)\cap \N$ and $\beta>0$. Then $H^1_{BZ1,M,\epsilon}(\Gamma) \cap L^2(\Gamma) \subset E^1_{quad,\beta}(\Gamma) $ and
$$\|f\|_{H^1_{quad,\beta}} \lesssim \|f\|_{H^1_{BZ1,M,\epsilon}}$$
\end{prop}

\begin{dem}
Let $f\in H^1_{BZ1,M,\epsilon} \cap L^2(\Gamma)$. 
Then there exist $(\lambda_i)_{i\in \N} \in \ell^1$ and $(a_i)_{i\in \N}$ a sequence of $(BZ1,M,\epsilon)$-molecules such that $f = \sum \lambda_i a_i$ and 
$$\sum_{i\in \N} |\lambda_i| \simeq \|f\|_{H^1_{BZ1,M,\epsilon}}.$$

First, since  $\|P^k\|_{1\to 1}\leq 1$  for all $k\in \N$, the operators $\Delta^\beta$ and then $\Delta^\beta P^{l-1}$ are $L^1$-bounded for $\beta >0$ (see \cite{CSC}). Consequently,
$$\Delta^\beta P^{l-1} \sum_{i\in \N} \lambda_i a_i = \sum_{i\in \N} \lambda_i \Delta^\beta P^{l-1} a_i.$$
Since the space $\Gamma$ is discrete, the $L^1$-convergence implies the pointwise convergence, that is, for all $x\in \Gamma$,
\[\begin{split}
   \left|\Delta^\beta P^{l-1} \sum_{i\in \N} \lambda_i a_i(x)\right| & = \left|\sum_{i\in \N} \lambda_i \Delta^\beta P^{l-1} a_i(x) \right| \\
& \leq \sum_{i\in \N} |\lambda_i| \left|\Delta^\beta P^{l-1} a_i(x) \right|.
  \end{split}\]
From here, the estimate
$$\|L_\beta f\|_{L^1}  = \left\| L_\beta \sum_{i\in \N} \lambda_i a_i \right\|_{L^1} \lesssim \sum_{i\in \N} |\lambda_i| \|L_\beta a_i\|_{L^1}$$
is  just  a consequence of the generalized Minkowski inequality.

It remains to prove that there exists a constant $C$ such that for all $(BZ1,M,\epsilon)$-molecules $a$, one has
\begin{equation} \label{BoundenessOnAtoms}
 \|L_\beta a\|_{L^1} \leq C.
\end{equation}

Let $s\in \N^*$, $(s_1,\dots,s_M)\in \bb s,2s\bn^M$ and a ball $B$ associated with the molecule $a$. By H\"older inequality and the doubling property, we may write
\begin{equation} \label{BoundenessOnAtoms2}
 \|L_\beta a\|_{L^1} \lesssim \sum_{j=1}^\infty V(2^jB)^{\frac{1}{2}} \|L_\beta a\|_{L^2(C_j(B))}.
\end{equation}

We will estimate now each term $\|L_\beta a\|_{L^2(C_j(B))}$. 

The result is then a consequence of Lemma \ref{OffDiagonalDecay} which can be reformulated as follows
\begin{equation} \label{ODDlem2} \left\| L_\beta (I-P^{s_1})\dots(I-P^{s_M}) [f\1_F] \right\|_{L^2(E)} \leq C_M \left(1+\frac{d(E,F)^2}{s}\right)^{-M} \|f\|_{L^2(F)}. \end{equation}

Notice that 
$$d(C_k(B),C_j(B)) \simeq \left\{ \begin{array}{ll} 0 & \text{ si } |j-k| \leq 1  \\
                                                    2^j\sqrt{s} & \text{ si } k \leq j-2 \\
                                                    2^k\sqrt{s} & \text{ si } k \geq j+2
                             \end{array}
\right. .$$
Thus,
\[\begin{split}
   \|L_\beta a\|_{L^2(C_j(B))} & \leq \sum_{k\geq 1} \|L_\beta(I-P^{s_1})\dots(I-P^{s_M}) [b\1_{C_k(B)}]\|_{L^2(C_j(B))} \\
& \lesssim \sum_{k \leq j-2} 4^{-jM} \|b\|_{L^2(C_k(B))} + \sum_{k=j-1}^{j+1} \|b\|_{L^2(C_k(B))} + \sum_{k\geq j+2} 4^{-kM} \|b\|_{L^2(C_k(B))} \\
& \lesssim \sum_{k \leq j-2} 4^{-jM} 2^{-\epsilon k} V(2^kB)^{-\frac{1}{2}} + 2^{-\epsilon j} V(2^jB)^{-\frac{1}{2}} + \sum_{k\geq j+2} 4^{-kM} 2^{-\epsilon j} V(2^kB)^{-\frac{1}{2}} \\
& \lesssim 2^{-\bar\epsilon j} V(2^jB)^{-\frac{1}{2}}
  \end{split}\]
where $\bar\epsilon = \min\{\epsilon, 2M-\frac{d_0}{2}\}$.

As a consequence, one has
\[\begin{split}
   \|L_\beta a\|_{L^1} & \lesssim  \sum_{j\geq 1} 2^{-\bar\epsilon j} \left(\frac{V(2^jB)}{V(2^jB)}\right)^{\frac{1}{2}} \\
& < +\infty .
  \end{split}\]
\end{dem}

\begin{prop} \label{Main41}
Let $(\Gamma,\mu)$ satisfying \eqref{UE}, $M\in \N^*$, $\epsilon>0$ and $\beta>0$. Then $H^1_{BZ1,M}(\Gamma) \cap L^2(\Gamma) \subset E^1_{quad,\beta}(\Gamma) $ and
$$\|f\|_{H^1_{quad,\beta}} \lesssim \|f\|_{H^1_{BZ1,M,\epsilon}}$$
\end{prop}

\begin{dem}
As in the proof of Proposition \ref{Main31}, it remains to check that for all $(BZ1,M,\epsilon)$-molecules 
$a = (I-P^{s_1})\dots (I-P^{s_M}) b$ associated with $s\in \N$, $(s_1,\dots,s_M)$ and $B = B(x_B,r_B)$, one has
$$\sum_{j=1}^\infty V(2^jB)^{\frac{1}{2}} \|L_\beta a\|_{L^2(C_j(B))} \lesssim 1.$$
The case $j=1$ follows from the $L^2$-boundedness of $L_\beta$ and of $(I-P^s)^M$, thus
$$\|L_\beta a\|_{L^2(C_j(B))} \lesssim \|a\|_{L^2} \lesssim \frac{1}{V(B)^{\frac{1}{2}}}.$$
For the case $j\geq 2$, we introduce $\tilde C_j(B)$ defined by
$$\tilde C_j(B) = \bigcup_{1\leq k \leq j-2} C_k(B).$$
Check that $\tilde C_j(B)$, $C_j(B)$, and $x_B$ satisfy \eqref{assumption1GT}, since $d(\tilde C_j(B),C_j(B)) \gtrsim 2^{j}r_B$. 
Thus, Lemma \ref{OffDiagonalDecayBis} yields
\[\begin{split}
   \|L_\beta a\|_{L^2(C_j(B))} & \leq \|L_\beta (I-P^s)^M [b\1_{\tilde C_j(B)}]\|_{L^2(C_j(B))} + \|L_\beta (I-P^s)^M [b\1_{\Gamma \backslash \tilde C_j(B)}]\|_{L^2(C_j(B))} \\
& \lesssim \frac{4^{-jM}}{V(x_B,2^{j}r_B)^{\frac{1}{2}}}  \|b\|_{L^1(\tilde C_j(B))} +  \|b\|_{L^2(\Gamma\backslash \tilde C_j(B))} \\
& \lesssim \frac{2^{-j\bar\epsilon}}{V(2^jB)^{\frac{1}{2}}} 
  \end{split}\]
where $\bar \epsilon = \min\{2M,\epsilon\}$. Summing in $j\geq 1$ ends the proof.
\end{dem}

\subsection{$H^1_{BZ2,M+\frac12,\epsilon} \cap H^2 \subset E^1_{quad,\beta} $: the case of 1-forms}

\begin{prop} \label{Main51}
Let $\epsilon>0$ , $M \in (\frac{d_0}{4}-\frac{1}{2},+\infty)\cap \N$. Then $H^1_{BZ2,M+\frac{1}{2}}(T_\Gamma) \cap H^2(T_\Gamma) \subset E^1_{quad,\frac{1}{2}}(T_\Gamma) $ and
$$\|f\|_{H^1_{quad,\frac{1}{2}}} \lesssim \|f\|_{H^1_{BZ2,M+\frac{1}{2},\epsilon}}.$$
\end{prop}

\begin{dem}
Let $F\in H^1_{BZ2,M+\frac{1}{2},\epsilon}(T_\Gamma) \cap H^2(T_\Gamma)$. 
Then there exist $(\lambda_i)_{i\in \N} \in \ell^1$ and $(a_i)_{i\in \N}$ a sequence of $(BZ2,M+\frac{1}{2},\epsilon)$-molecules such that $F = \sum \lambda_i a_i$ and 
$$\sum_{i\in \N} |\lambda_i| \simeq \|f\|_{H^1_{BZ1,M,\epsilon}}.$$

First, by $L^1$-boundedness of the operators $P$ and $d^*$ (see Proposition \ref{danddstar}) and by the Minkowski inequality, one has
\[\begin{split}
 \|L_\frac12 \Delta^{-\frac{1}{2}} d^*F\|_{L^1}  
& = \sum_{x\in \Gamma} m(x) \left( \sum_{l\geq 1} \sum_{y\in B(x,\sqrt{l})} m(y) |P^{l-1} d^* F(y)|^2 \right)^\frac12 \\
 & = \sum_{x\in \Gamma} m(x) \left( \sum_{l\geq 1} \sum_{y\in B(x,\sqrt{l})} m(y) |P^{l-1} d^* \sum_{i\in \N} \lambda_i a_i(y)|^2 \right)^\frac12 \\
 & \lesssim \sum_{i\in \N} |\lambda_i| \sum_{x\in \Gamma} m(x) \left( \sum_{l\geq 1} \sum_{y\in B(x,\sqrt{l})} m(y) |P^{l-1} d^*  a_i(y)|^2 \right)^\frac12.
\end{split} 
 \]

It remains to prove that there exists a constant $C$ such that for all $(BZ2,M+\frac{1}{2},\epsilon)$-molecules $a$, one has
\begin{equation} \label{z39}
 \sum_{x\in \Gamma} m(x) \left( \sum_{l\geq 1} \sum_{y\in B(x,\sqrt{l})} m(y) |P^{l-1} d^*  a_i(y)|^2 \right)^\frac12 \lesssim 1.
\end{equation}
Let $a = d\Delta^{-\frac{1}{2}} (I-(I+s\Delta)^{-1})^{M+\frac12} b$ be a $(BZ_2,M+\frac12,\epsilon)$-molecule associated with $s\in \N^*$ and the ball $B$. Since $d^*d\Delta^{-\frac12} = \Delta^\frac12$, \eqref{z39} becomes
$$\|L_\frac12 (I-(I+s\Delta)^{-1})^{M+\frac12} b\| \lesssim 1.$$

We end the proof as we did for Proposition \ref{Main31}, using Lemma \ref{OffDiagonalDecay4} instead of Lemma \ref{OffDiagonalDecay}.
\end{dem}

\begin{prop} \label{Main42}
Let $(\Gamma,\mu)$ satisfying \eqref{UE}. Let $\epsilon>0$ , $M \in \N$. Then $H^1_{BZ2,M+\frac{1}{2}}(T_\Gamma) \cap H^2(T_\Gamma) \subset E^1_{quad,\frac{1}{2}}(T_\Gamma)$ and
$$\|f\|_{H^1_{quad,\frac{1}{2}}} \lesssim \|f\|_{H^1_{BZ2,M+\frac{1}{2},\epsilon}}.$$
\end{prop}

\begin{dem}
We begin the proof as the one of Proposition \ref{Main51}. 
We end the proof as Proposition \ref{Main41} instead of Proposition \ref{Main31}, using Lemma \ref{OffDiagonalDecay5} instead of Lemma \ref{OffDiagonalDecay4}.
\end{dem}

\subsection{$E^1_{quad,\beta} \subset H^1_{BZ2,M,\epsilon}  \cap L^2$: the case of functions}

In this paragraph, we will need  a few  results on tents spaces (see \cite{CMS}, \cite{Russ2}, \cite{HLMMY}). However, we need in our proofs some "discrete" tent spaces, defined below:

\begin{defi}
For $x\in \Gamma$, we recall
$$\gamma(x) = \left\{ (y,k) \in \Gamma \times \N, \, d(x,y)^2 \leq k \right\}$$
and for a set $O \subset \Gamma$, we define
$$\hat O = \{(y,k) \in \Gamma \times \N, \, d(x,O^c)^2 > k \}.$$
For a function $F$ defined on $\Gamma \times \N$, consider for all $x\in \Gamma$
$$\mathcal A F(x) = \left( \sum_{(y,k)\in \gamma(x)} \frac{1}{k+1}\frac{m(y)}{V(x,\sqrt{k+1})} |F(y,k)|^2 \right)^{\frac{1}{2}}.$$
For $p\in [1,+\infty)$, the tent space $T^p_2(\Gamma)$ is defined as the space of functions $F$ on $\Gamma \times \N$ for which $ \mathcal AF \in L^p(\Gamma)$, 
and is outfitted with the norm $\|F\|_{T_2^p} = \|\mathcal A F\|_{L^p}$ (the space $T^p_2$ is then complete).
\end{defi}

%

\begin{defi}
 A function $A$ on $\Gamma \times \N$ is said to be a $T^1_2$-atom if there exists a ball $B\subset \Gamma$ such that $A$ is supported in $\hat B$ and
$$ \|A\|_{T^2_2}^2 := \sum_{(x,k)\in\hat B } \frac{m(x)}{k+1}|A(x,k)|^2 \leq \frac{1}{V(B)}.$$
\end{defi}

\begin{prop} \label{AtomicDecompTentSpaces}
For every element $F \in T^1_2(\Gamma)$, there exist a scalar sequence $(\lambda_i)_{i\in \N} \subset \ell^1$ and a sequence of $T^1_2$-atoms $(A_i)_{i\in \N}$ such that
\begin{equation} \label{decompositionT1}
 F = \sum_{i=0}^{+\infty} \lambda_i A_i \qquad \text{ in $T^1_2(\Gamma)$}.
\end{equation}

Moreover,
$$\sum_{i\geq 0} |\lambda_i| \simeq \|F\|_{T^1_2}$$
where the implicit constants only depend on the  constant in \eqref{DV}. 
Finally, if $F \in T^1_2(\Gamma) \cap T^2_2(\Gamma)$, then the decomposition \eqref{decompositionT1} also converges in $T^2_2(\Gamma)$. 
\end{prop}

\begin{dem}
 This proof is analogous to the one of Theorem 1.1 in \cite{Russ2} and of Theorem 4.10 in \cite{HLMMY} with obvious modifications.
\end{dem}

We introduce the functional $\pi_{\eta,\beta}:  \ T^2_2(\Gamma)\to L^2(\Gamma)$ defined for any real $\beta>0$ and any integer $\eta \geq \beta$ by
$$\pi_{\eta,\beta} F(x) = \sum_{l\geq 1} \frac{c_l^\eta}{l^\beta} \left[\Delta^{\eta-\beta}(I+P)^\eta P^{l-1} F(.,l-1)\right](x)$$
where $\ds \sum_{l\geq 1} c_{l}^\eta z^{l-1}$ is the Taylor series of the function $(1-z)^{-\eta}$.

\begin{lem}
 The operator $\pi_{\eta,\beta}$ is bounded from $T^2_2(\Gamma)$ to $L^2(\Gamma)$.
\end{lem}

\begin{dem}
Let $g\in L^2(\Gamma)$. Then, for all $F \in T_2^2(\Gamma)$,
\[\begin{split}
   \left<\pi_{\eta,\beta} F,g \right> & = \sum_{l\geq 1}  \frac{c_l^\eta}{l^\beta} \left< \Delta^{\eta-\beta}(I+P)^\eta P^{l-1} F(.,l) , g \right> \\
& =  \sum_{l\geq 1}  \frac{c_l^\eta}{l^\beta} \left<  F(.,l-1) , \Delta^{\eta-\beta}(I+P)^\eta P^{l-1} g \right> \\
& \leq \sum_{l\geq 1}  \frac{c_l^\eta}{l^\beta} \|F(.,l-1)\|_{L^2} \|\Delta^{\eta-\beta}(I+P)^\eta P^{l-1} g\|_{L^2} \\
& \leq \left(\sum_{l\geq 1} \frac{1}{l} \|F(.,l-1)\|_{L^2}^2 \right)^\frac{1}{2} \left(\sum_{l\geq 1} l^{1-2\beta}(c^\eta_l)^2 \|\Delta^{\eta-\beta}(I+P)^\eta P^{l-1} g\|_{L^2}^2 \right)^\frac{1}{2} \\
& \lesssim \|F\|_{T_2^2} \|(I+P)^\eta g\|_{L^2} \\
& \lesssim \|F\|_{T^2_2} \|g\|_{L^2} 
  \end{split}\]
where the last but one line comes from the $L^2$-boundedness of Littlewood-Paley functionals (since $l^{1-2\beta}(c^\eta_l)^2 \simeq l^{2(\eta-\beta)-1}$, see \cite{Fen1}, Lemma B.1).
\end{dem}

\begin{lem} \label{molecules}
 Suppose that $A$ is a $T^1_2(\Gamma)$-atom associated with a ball $B \subset \Gamma$. 
Then for every $M\in \N^*$, $\beta>0$ and $\epsilon\in (0,+\infty)$, there exist an integer $\eta = \eta_{M,\beta,\epsilon}$ and a uniform constant $C_{M,\beta,\epsilon}>0$ such that 
$C_{M,\beta,\epsilon}^{-1} \pi_{\eta,\beta}(A)$ is a $(BZ_2,M,\epsilon)$-molecule associated with the ball $B$.
\end{lem}

\begin{dem}
Let $\eta = \lceil \frac{d_0}{4} + \frac{\epsilon}{2} + \beta \rceil + M + 1$, that is the only integer such that 
$$\eta \geq \frac{d_0}{4} + \frac{\epsilon}{2} + \beta + M + 1 > \eta -1.$$

Let $A$ be a $T^1_2$-atom associated with a ball $B$ of radius $r$. We write
$$a : = \pi_{\eta,\beta}(A) = (I-(I+r^2\Delta)^{-1})^M b$$
where
$$b:= \sum_{l\geq 1} \frac{c_l^\eta}{l^\beta} \left(\frac{I+r^2\Delta}{r^2}\right)^M\Delta^{\eta-\beta-M}(I+P)^\eta P^{l-1} A(.,l-1)$$
Let us check that $a$ is a $(BZ_2,M,\epsilon)$-molecule associated with $B$, up to multiplication by some harmless constant $C_{M,\epsilon}$. 
First, one has, for all $g\in L^2(4 \eta B)$,
\[\begin{split}
   \left| \left< b,g \right> \right|& \leq   \sum_{m=0}^M \frac{c_m}{r^{2(M-m)}} \sum_{l\geq 1} \frac{c_l^\eta}{l^\beta} \left|\left<\Delta^{\eta-\beta-M+m}(I+P)^\eta P^{l-1} A(.,l-1),g \right> \right| \\
& = \sum_{m=0}^M \frac{c_m}{r^{2(M-m)}} \sum_{l\geq 1} \frac{c_l^\eta}{l^\beta} \left|\left< A(.,l-1),\Delta^{\eta-\beta-M+m}(I+P)^\eta P^{l-1} g \right>\right| \\
& \lesssim \sum_{m=0}^M \frac{1}{r^{2(M-m)}} \sum_{l\geq 1} l^{\eta-\beta-1} \|A(.,l-1)\|_{L^2(B)}\|\Delta^{\eta-\beta-M+m}(I+P)^\eta P^{l-1}g\|_{L^2(B)}\\ 
& \lesssim \sum_{m=0}^M \frac{1}{r^{2(M-m)}} \|A\|_{T^2_2} \left(  \sum_{l= 1}^{r^2} l^{2(\eta-\beta)-1} \|\Delta^{\eta-\beta-M+m}(I+P)^\eta P^{l-1}g\|_{L^2(B)}^2  \right)^{\frac12} \\
& \lesssim \sum_{m=0}^M \|A\|_{T^2_2} \left(  \sum_{l= 1}^{r^2} l^{2(\eta-\beta-M+m)-1} \|\Delta^{\eta-\beta-M+m}(I+P)^\eta P^{l-1}g\|_{L^2(B)}^2  \right)^{\frac12} \\
& \lesssim  \|A\|_{T^2_2} \sum_{m=0}^M \left\| G_{\eta-\beta-M+m}(I+P)^\eta g \right\|_{L^2} \\
& \lesssim  \|A\|_{T^2_2} \left\| (I+P)^\eta g \right\|_{L^2} \\
& \lesssim \frac{1}{V(B)^{\frac{1}{2}}} \left\|g \right\|_{L^2}
\end{split}\]
where we used the $L^2$-boundedness of the quadratic Littlewood-Paley functional for the last but one line (see \cite{Fen1}, \cite{al}).

Let $j> \log_2(\eta) + 1$ and $g\in L^2(C_j(B))$. Since $\Supp (I+P)^\eta g \in C_{j,\eta}(B) = \{x\in \Gamma, d(x,C_j(B))\leq \eta\}$ and $d(C_{j,\eta}(B),B) \gtrsim 2^jr$, 
\[\begin{split}
   \left| \left< b,g \right> \right|&  \lesssim \sum_{m=0}^M \frac{1}{r^{2(M-m)}} \|A\|_{T^2_2} \left(  \sum_{l= 1}^{r^2} l^{2(\eta-\beta)-1} \|\Delta^{\eta-\beta-M+m}(I+P)^\eta P^{l-1}g\|_{L^2(B)}^2  \right)^{\frac12} \\
& \lesssim \sum_{m=0}^M r^{2(\eta-\beta-M-1)} \|A\|_{T^2_2} \left(  \sum_{l= 1}^{r^2} l^{2(1+m)-1} \|\Delta^{\eta-\beta-M+m}(I+P)^\eta P^{l-1}g\|_{L^2(B)}^2  \right)^{\frac12} \\
& \lesssim r^{2(\eta-\beta-M-1)} \|A\|_{T^2_2} \sum_{m=0}^M \left\| G_{1+m} \Delta^{\eta-\beta-M-1} (I+P)^\eta g \right\|_{L^2(B)} \\
& \lesssim \frac{r^{2(\eta-\beta-M-1)}}{(4^j r^2)^{\eta-\beta-M-1}} \|A\|_{T^2_2} \left\| (I+P)^\eta g \right\|_{L^2} \\
& \lesssim 2^{-j(\frac{d_0}{2} + \epsilon)} \|A\|_{T^2_2} \left\|g \right\|_{L^2} \\
& \lesssim \frac{2^{-j\epsilon}}{V(2^jB)^{\frac{1}{2}}} \left\|g \right\|_{L^2}
\end{split}\]
where we used Lemma \ref{OffDiagonalDecayTer} for the last but two line and Proposition \ref{propDV} for the last one.
We conclude that, up to multiplication by some harmless constant, $b$ is a $(BZ_2,M,\epsilon)$-molecule.
\end{dem}

\begin{prop} \label{Main32}
Let $M \in \N^*$, $\epsilon>0$ and $\beta>0$. Then $E^1_{quad,\beta}(\Gamma)  \subset H^1_{BZ2,M,\epsilon}(\Gamma) \cap L^2(\Gamma)$ and
$$\|f\|_{H^1_{quad,\beta}} \lesssim \|f\|_{H^1_{BZ2,M,\epsilon}}.$$
\end{prop}

\begin{dem} 
Let $f\in E^1_{quad,\beta}(\Gamma)$. We set
$$F(.,l) = [(l+1)\Delta]^\beta P^{l}f.$$
By definition of $H^1_{quad,\beta}(\Gamma)$, one has that $F \in T^1_2(\Gamma)$. Moreover, since $f\in L^2(\Gamma)$, $L^2$-boundedness of Littlewood-Paley functionals (see \cite{BRuss}, \cite{Fen1}) yields that $F \in T^2_2(\Gamma)$.
Thus, according to Lemma \ref{AtomicDecompTentSpaces}, there exist a numerical sequence $(\lambda_i)_{i\in \N}$ and a sequence of $T^1_2$-atoms $(A_i)_{i\in \N}$ such that
$$F = \sum_{i=0}^\infty \lambda_i A_i \qquad \text{ in } T^1_2(\Gamma) \text{ and } T^2_2(\Gamma)$$
and
$$\sum_{i\in \N} |\lambda_i| \lesssim \|F\|_{T^1_2} = \|f\|_{H^1_{quad,\beta}}.$$
Choose $\eta$ as in Lemma \ref{molecules}. Using Corollary \ref{L2convergenceCor}, since $f\in L^2(\Gamma)$,
\begin{equation} \label{L2DecomMol}\begin{split}
f & =  \pi_{\eta,\beta} F(.,l)  \\
& = \sum_{i = 0}^{+\infty} \lambda_i \pi_{\eta,\beta} (A_i)
  \end{split} \end{equation}
where the sum converges in $L^2(\Gamma)$. 
According to Lemma \ref{molecules}, $\pi_{\eta,\beta} (A_i)$ are molecules and 
then \eqref{L2DecomMol} would provide a $(M,\epsilon)$-representation of $f$ if the convergence  held  in $L^1(\Gamma)$.
By uniqueness of the limit, it remains to prove that $\sum \lambda_i \pi_{\eta,\beta} (A_i)$ converges in $L^1$. 
Indeed,
\[\begin{split}
   \sum_{i\in \N} |\lambda_i| \left\|\pi_{\eta,\beta} (A_i)\right\|_{L^1} & \lesssim \sum_{i\in \N} |\lambda_i| \\
& < +\infty
 \end{split}\]
where the first line comes from Proposition \ref{BoundedMolecules3} and the second one from the fact that $(\lambda_i)_{i\in\N} \in \ell^1(\N)$.
\end{dem}

\subsection{$E^1_{quad,\beta}  \subset H^1_{BZ2,M+\frac12,\epsilon} \cap H^2$: the case of 1-forms}

\begin{lem} \label{molecules2}
 Suppose that $A$ is a $T^1_2(\Gamma)$-atom associated with a ball $B \subset \Gamma$. 
Let $M \in \N$ and $\epsilon>0$, there exist an integer $\eta = \eta_{M,\epsilon}$ and a uniform constant $C_{M,\epsilon}>0$ such that 
$C_{M,\epsilon}^{-1} d\Delta^{-\frac{1}{2}}\pi_{\eta,\frac{1}{2}}(A)$ is a $(BZ_2,M+\frac{1}{2},\epsilon)$-molecule associated with the ball $B$.
\end{lem}

\begin{dem}
Let $\eta = \lceil \frac{d_0}{4} + \frac{\epsilon}{2}\rceil + M + 2$. We will also write $t$ for $\lceil \frac{d_0}{4} + \frac{\epsilon}{2}\rceil \in \N^*$.

Let $A$ be a $T^1_2$-atom associated with a ball $B$ of radius $r$. We write
$$a : = d\Delta^{-\frac{1}{2}}\pi_{\eta,\frac{1}{2}}(A) = r^{2M+1 }d\Delta^{M}(I+r^2\Delta)^{-M-\frac{1}{2}} b$$
where
\begin{equation}\label{bexpression}\begin{split} b & := \sum_{l\geq 1} \frac{c_l^\eta}{\sqrt{l}} \left(\frac{I+r^2\Delta}{r^2}\right)^{M+\frac{1}{2}}\Delta^{\eta-1-M}(I+P)^\eta P^{l-1} A(.,l-1) \\
   & = \sqrt{\frac{r^2}{1+r^2}}\sum_{k=0}^\infty a_k \left(\frac{r^2}{1+r^2}\right)^k \sum_{l\geq 1} \frac{c_l^\eta}{\sqrt{l}} \left(\frac{I+r^2\Delta}{r^2}\right)^{M+1}\Delta^{1+t}(I+P)^\eta P^{l+k-1} A(.,l-1) \\
  \end{split}\end{equation}
where $\sum a_k z^k$ is the Taylor serie of the function $(1-z)^{-\frac{1}{2}}$ (cf \eqref{ResolvantEquality}).

Let us check that $a$ is a $(BZ_2,M+\frac{1}{2},\epsilon)$-molecule associated with $B$, up to multiplication by some harmless constant $C_{M,\epsilon}$. 

\smallskip

Let  $g\in L^2(4\eta B)$. One has with the first equality in \eqref{bexpression},
\[\begin{split}
   \left|\left< b,g \right>\right| & \leq  r^{-2M-1} \sum_{l\geq 1} \frac{c_l^\eta}{\sqrt{l}} \left|\left< A(.,l-1),\left(I+r^2\Delta\right)^{M+\frac{1}{2}}\Delta^{1+t}(I+P)^\eta P^{l-1} g \right>\right| \\
& \lesssim   r^{-2M-1} \sum_{l\geq 1} \frac{c_l^\eta}{\sqrt{l}} \|A(.,l-1)\|_{L^2(B)}\|\left(I+r^2\Delta\right)^{M+\frac{1}{2}}\Delta^{1+t}(I+P)^\eta P^{l-1}g\|_{L^2}\\ 
& \lesssim  \|A\|_{T^2_2} r^{-2M-1} \left(  \sum_{l= 1}^{r^2} l^{2(\eta-1)} \|\left(I+r^2\Delta\right)^{M+\frac{1}{2}}\Delta^{1+t}(I+P)^\eta P^{l-1}g\|_{L^2}^2  \right)^{\frac12} \\
& \lesssim  \|A\|_{T^2_2} r^{-2M-1} \left(  \sum_{l= 1}^{r^2} l^{2(1+t+M)} \|\left(I+r^2\Delta\right)^{M+\frac{1}{2}}\Delta^{1+t}(I+P)^\eta P^{l-1}g\|_{L^2}^2  \right)^{\frac12} \\
& \lesssim  \|A\|_{T^2_2} \|(I+P)^\eta g\|_{L^2} \\
& \lesssim  \|A\|_{T^2_2} \|g\|_{L^2}
 \end{split}\]
where we use that the functionals $\ds g \mapsto r^{-2M-1}\left( \sum_{l=1}^{r^2} l^{2(1+t+M)} |\left(I+r^2\Delta\right)^{M+\frac{1}{2}}\Delta^{1+t}P^{l-1}g|^2 \right)^{1/2}$ are $L^2$-bounded uniformly in $r$. 
Indeed, since $(-1) \notin Sp(P)$, functional calculus provides, for some $a>-1$,
\[\begin{split} 
\|\left(I+r^2\Delta\right)^{M+\frac{1}{2}}\Delta^{1+t}P^{l-1}g\|_{L^2}^2 & = \int_{a}^1 (1+r^2(1-\lambda))^{2M+1}(1-\lambda)^{2(1+t)} \lambda^{2(l-1)} dE_{gg}(\lambda) \\
& \lesssim \int_{a}^1 \left[1+r^{2(2M+1)}(1-\lambda)^{2M+1}\right](1-\lambda)^{2(1+t)} \lambda^{2(l-1)} dE_{gg}(\lambda).
\end{split}\]
Thus,
\[\begin{split}
   r^{-2(2M+1)} & \sum_{l=1}^{r^2} l^{2(1+t+M)}  \|\left(I+r^2\Delta\right)^{M+\frac{1}{2}}\Delta^{1+t}P^{l-1}g\|_{L^2}^2 \\
& \lesssim \int_{a}^1 (1-\lambda)^{2(1+t)} \sum_{l=1}^{r^2} l^{2(1+t)-1}\lambda^{2(l-1)} dE_{gg}(\lambda) +  \int_{a}^1 (1-\lambda)^{2(1+t+M)+1} \sum_{l=1}^{r^2} l^{2(1+t+M)} \lambda^{2(l-1)} dE_{gg}(\lambda) \\
& \lesssim \int_{a}^1 (1-\lambda)^{2(1+t)} \sum_{l=1}^{\infty} l^{2(1+t)-1}\lambda^{2(l-1)} dE_{gg}(\lambda) +  \int_{a}^1 (1-\lambda)^{2(1+t+M)+1} \sum_{l=1}^{\infty} l^{2(1+t+M)}\lambda^{2(l-1)} dE_{gg}(\lambda) \\
& \lesssim \int_{a}^1 \frac{(1-\lambda)^{2(1+t)}}{(1-\lambda^2)^{2(1+t)}} dE_{gg}(\lambda) +  \int_{a}^1 \frac{(1-\lambda)^{2(1+t+M)+1}}{(1-\lambda^2)^{2(1+t+M)+1}} dE_{gg}(\lambda) \\
& \quad = \int_{a}^1 \left[(1+\lambda)^{-2(1+t)} + (1+\lambda)^{-2(1+t+M)-1} \right] dE_{gg}(\lambda) \\
& \lesssim \int_{a}^1 dE_{gg}(\lambda) = \|g\|_{L^2}^2
  \end{split}\]
  where the third inequality comes from the fact that $l^{\xi-1} \sim c^{\xi}_l$ (see Lemma B.1 in \cite{Fen1}).

\smallskip

Let $j> \log_2(\eta) + 1$ and $g\in L^2(C_j(B))$. One has $d(C_{j,\eta}(B),B) \gtrsim 2^jr$ (cf Lemma \ref{molecules}). The second identity in \eqref{bexpression} provides

\[\begin{split}
   \left|\left< b,g \right>\right| & \leq  \sum_{m=0}^{M+1} \frac{c_m}{r^{2(M+1-m)}} \sqrt{\frac{r^2}{1+r^2}} \sum_{k=0}^\infty a_k \left(\frac{r^2}{1+r^2}\right)^k \sum_{l\geq 1} \frac{c_l^\eta}{\sqrt{l}} \left|\left< A(.,l-1),\Delta^{1+t+m}(I+P)^\eta P^{l+k-1} g \right>\right| \\
& \lesssim \sum_{m=0}^{M+1} \frac{c_m}{r^{2(M+1-m)}} \sum_{k=0}^\infty a_k \left(\frac{r^2}{1+r^2}\right)^k \sum_{l\geq 1} \frac{c_l^\eta}{\sqrt{l}} \|A(.,l-1)\|_{L^2(B)}\|\Delta^{1+t+m}(I+P)^\eta P^{l+k-1}g\|_{L^2(B)}\\ 
& \lesssim  \|A\|_{T^2_2} \sum_{m=0}^{M+1} \frac{1}{r^{2(M+1-m)}} \sum_{k=0}^\infty a_k \left(\frac{r^2}{1+r^2}\right)^k\left(  \sum_{l= 1}^{r^2} l^{2(\eta-1)} \|\Delta^{1+t+m}(I+P)^\eta P^{l+k-1}g\|_{L^2(B)}^2  \right)^{\frac12} \\
& \lesssim  \|A\|_{T^2_2} \|(I+P)^\eta g\|_{L^2} \sum_{m=0}^{M+1} r^{2(\eta-M-2+m)} \sum_{k=0}^\infty a_k \left(\frac{r^2}{1+r^2}\right)^k \left(  \sum_{l= 1}^{r^2} \dfrac{e^{-c\frac{4^jr^2}{l+k}}}{(l+k)^{2(1+t+m)}}  \right)^{\frac12} \\
& \lesssim  \|A\|_{T^2_2} \|g\|_{L^2} \sum_{m=0}^{M+1} r^{2(t+m)} \sum_{k=0}^\infty a_k \left(\frac{r^2}{1+r^2}\right)^k \left(  \sum_{l= 1}^{\infty} \dfrac{e^{-c\frac{4^jr^2}{l+k}}}{(l+k)^{2(1+t+m)}}  \right)^{\frac12}  \\
& \lesssim  \|A\|_{T^2_2} \|g\|_{L^2} \sum_{m=0}^{M+1} r^{2(t+m)+1} \frac{1}{\sqrt{1+r^2}} \sum_{k=0}^\infty a_k \left(\frac{r^2}{1+r^2}\right)^k \frac{1}{(4^jr^2)^{t+m+\frac{1}{2}}}  \\
& \lesssim  \|A\|_{T^2_2} \|g\|_{L^2} \sum_{m=0}^{M+1} \frac{r^{2(t+m)+1}}{(4^jr^2)^{t+m+\frac{1}{2}}} (1+r^2(1-1))^{-\frac{1}{2}} \\
& \lesssim  2^{-j(2t+1)} \|A\|_{T^2_2} \|g\|_{L^2} \\
& \lesssim \frac{2^{-j\epsilon}}{V(2^jB)} \|g\|_{L^2}
  \end{split}\]
where we used the estimate \eqref{GUE} for the forth line.
\end{dem}

\begin{prop} \label{Main52}
Let $M\in \N$ and $\epsilon>0$. Then $E^1_{quad,\frac{1}{2}}(T_\Gamma) \subset H^1_{BZ2,M+\frac{1}{2}}(T_\Gamma) \cap H^2(T_\Gamma)$ and
$$\|G\|_{H^1_{quad,\frac{1}{2}}} \lesssim \|G\|_{H^1_{BZ2,M+\frac{1}{2},\epsilon}} \qquad \forall G \in H^1_{quad,\frac{1}{2}}(T_\Gamma) \cap E^2(T_\Gamma)$$
\end{prop}

\begin{dem}
Let $G\in E^1_{quad,\frac{1}{2}}(T_\Gamma)$. We set
$$F(.,l) = \sqrt{l+1} P^{l}d^*G.$$
By definition of $H^1_{quad,\frac{1}{2}}(T_\Gamma)$, one has that $F \in T^1_2(\Gamma)$. 
Moreover, Proposition \ref{danddstar} yields that $\Delta^{-\frac{1}{2}}d^* G \in L^2(G)$ and therefore $F\in T^2_2(\Gamma)$ with the $L^2$-boundedness of Littlewood-Paley functionals.

Thus, according to Lemma \ref{AtomicDecompTentSpaces}, there exist a scalar sequence $(\lambda_i)_{i\in \N} \in \ell^1(\N)$ and a sequence of $T^1_2$-atoms $(A_i)_{i\in \N}$ such that
$$F = \sum_{i=0}^\infty \lambda_i A_i \qquad \text{ in } T^1_2(\Gamma) \text{ and in } T^2_2(\Gamma)$$
and
$$\sum_{i\in \N} |\lambda_i| \lesssim \|F\|_{T^1_2} = \|G\|_{H^1_{quad,\frac{1}{2}}}.$$
Choose $\eta$ as in Lemma \ref{molecules}. Using Lemma \ref{L2convergenceCor}, since $\Delta^{-\frac{1}{2}} d^* G\in L^2(\Gamma)$,
\[ \begin{split}
\Delta^{-\frac{1}{2}} d^* G& =  \pi_{\eta,\frac{1}{2}} F(.,l)  \\
& = \sum_{i = 0}^{+\infty} \lambda_i \pi_{\eta,\frac{1}{2}} (A_i)
  \end{split}\]
where the sum converges in $L^2(\Gamma)$.
Recall that $d \Delta^{-1} d^* = Id_{H^2(T_\Gamma)}$. Moreover, $d\Delta^{-\frac{1}{2}}$ is bounded from $L^2(\Gamma)$ to $L^2(T_\Gamma)$ (see Proposition \ref{danddstar}). Then
\begin{equation} \label{L2DecomMol2}
 G = \sum_{i=0}^{+\infty} \lambda_i d\Delta^{-\frac{1}{2}} \pi_{\eta,\frac{1}{2}} (A_i)
\end{equation}
where the sum converges in $L^2(T_\Gamma)$. 
According to Lemma \ref{molecules2}, $d\Delta^{-\frac{1}{2}} \pi_{M,\frac12} (A_i)$ are $(BZ2,M+\frac{1}{2},\epsilon)$-molecules and 
then \eqref{L2DecomMol2} would provide a $(BZ2,M+\frac{1}{2},\epsilon)$-representation of $f$ if the convergence  held  in $L^1(\Gamma)$.
By uniqueness of the limit, it remains to prove that $\sum \lambda_i d\Delta^{-\frac{1}{2}}\pi_{\eta,\frac{1}{2}} (A_i)$ converges in $L^1$. 
Indeed,
\[\begin{split}
   \sum_{i\in \N} |\lambda_i| \left\|d\Delta^{-\frac12}\pi_{\eta,\frac{1}{2}} (A_i)\right\|_{L^1(T_\Gamma)} & \lesssim \sum_{i\in \N} |\lambda_i| \\
& < +\infty
 \end{split}\]
where the first line comes from Corollary \ref{BoundedMolecules4} and the second one because $(\lambda_i)_{i\in\N} \in \ell^1(\N)$.
\end{dem}

\subsection{Proof of Theorems \ref{Main3}, \ref{Main5} and \ref{Main4}}

\begin{dem} \em (Theorem \ref{Main3}) \em 

Let $\beta>0$, $M\in \N^*\cap (\frac{d_0}{4},+\infty)$ and $\epsilon>0$. Propositions \ref{Main31} and \ref{Main32} yield the continuous embeddings
$$H^1_{BZ1,M,\epsilon}(\Gamma) \cap L^2(\Gamma) \subset E^1_{quad,\beta}(\Gamma) \subset  H^1_{BZ2,M,\epsilon}(\Gamma) \cap L^2(\Gamma).$$
However, Theorem \ref{Main2} states that $H^1_{BZ1,M,\epsilon}(\Gamma) = H^1_{BZ2,M,\epsilon(\Gamma)}$. Thus, we deduce
\begin{equation} \label{equalityE1E1}H^1_{BZ1,M,\epsilon}(\Gamma) \cap L^2(\Gamma) = E^1_{quad,\beta}(\Gamma) =  H^1_{BZ2,M,\epsilon}(\Gamma) \cap L^2(\Gamma) \end{equation}
with equivalent norms. In particular, $E^1_{quad,\beta}(\Gamma)\subset L^1(\Gamma)$.\par
\noindent  Let us now prove that the completion of $E^1_{quad,\beta}(\Gamma)$ in $L^1$ exists. To that purpose, it is enough (see Proposition 2.2 in \cite{AMM}) to check that, for all Cauchy sequences $(f_n)_n$ in $E^1_{quad,\beta}(\Gamma)$ that converges to 0 in $L^1(\Gamma)$, $f_n\rightarrow 0$ for the $\|.\|_{H^1_{quad,\beta}}$ norm. Equivalent norms in \eqref{equalityE1E1} implies that $(f_n)_n$ is a Cauchy sequence in $H^1_{BZ\kappa,M,\epsilon}(\Gamma)$ that converges to 0 in $L^1(\Gamma)$. Since $H^1_{BZ\kappa,M,\epsilon}(\Gamma)$ is complete, it follows that $f_n\rightarrow g$ for some $g\in  H^1_{BZ\kappa,M,\epsilon}(\Gamma)$, but then also for the $L^1$-norm, which entails that $g=0$.  Thus, $f_n\rightarrow 0$ for the norm $H^1_{BZ\kappa,M,\epsilon}(\Gamma)$ and so for the norm $\|.\|_{H^1_{quad,\beta}}$ (the norms being equivalent on $E^1_{quad,\beta}(\Gamma)$). 

Therefore, the completion $H^1_{quad,\beta}(\Gamma)$ of $E^1_{quad,\beta}(\Gamma)$ exists and is defined by
$$H^1_{quad,\beta}(\Gamma) = \{f\in F, \, \text{there exists $(f_n)_n$ Cauchy sequence in $E^1_{quad, \beta}(\Gamma)$ such that $f_n \to f$ in $L^1(\Gamma)$}\}.$$
The fact that $H^1_{quad,\beta}(\Gamma) = H^1_{BZ\kappa,M\epsilon}(\Gamma)$ is then a straightforward consequence of $\eqref{equalityE1E1}$ and the fact that the space $H^1_{BZ\kappa,M\epsilon}(\Gamma)\cap L^2(\Gamma)$ is dense in $H^1_{BZ\kappa,M\epsilon}(\Gamma)$.
\end{dem}

\begin{dem} \em (Theorem \ref{Main5}) \em 

Let $M\in \N\cap (\frac{d_0}{4}-\frac12,+\infty)$ and $\epsilon>0$. Propositions \ref{Main51} and \ref{Main52} yield the continuous embeddings
$$H^1_{BZ2,M+\frac12,\epsilon}(T_\Gamma) \cap H^2(T_\Gamma) \subset E^1_{quad,\frac12}(T_\Gamma) \subset  H^1_{BZ2,M+\frac 12,\epsilon}(T_\Gamma) \cap H^2(T_\Gamma),$$
from which we deduce the equality of the two spaces, with equivalent norms.

Since $\mathbb H^1_{BZ2,M+\frac12,\epsilon}(T_\Gamma)$ is dense in $ H^1_{BZ2,M+\frac12,\epsilon}(T_\Gamma) \subset L^1(T_\Gamma)$ and is included in $ H^1_{BZ2,M+\frac12,\epsilon}(T_\Gamma) \cap H^2(\Gamma)$, 
it follows that $H^1_{BZ2,M+\frac12,\epsilon}(T_\Gamma)$ is the completion in $L^1(T_\Gamma)$ of $H^1_{BZ2,M+\frac12,\epsilon}(T_\Gamma) \cap H^2(\Gamma)$ and thus also of $E^1_{quad,\frac12}(T_\Gamma)$ with the same arguments than those used in the proof of Theorem \ref{Main3}.

Moreover, notice that if $F \in H^2(T_\Gamma)$,
$$F \in E^1_{quad,\beta}(T_\Gamma)  \Longleftrightarrow  \Delta^{-\frac{1}{2}}d^*F \in  E^1_{quad,\beta}(\Gamma).$$
Indeed,  the implication $\Delta^{-\frac{1}{2}}d^*F \in  E^1_{quad,\beta}(\Gamma)\Rightarrow F \in E^1_{quad,\beta}(T_\Gamma)$ is obvious, and the converse is due to Proposition \ref{danddstar}.  
As said in Theorem \ref{Main3}, the spaces $E^1_{quad,\beta}(\Gamma)$ are all equivalent once $\beta >0$; 
and so are the spaces $E^1_{quad,\beta}(T_\Gamma)$. Consequently, for all $\beta>0$, the completion of $E^1_{quad,\beta}(T_\Gamma)$ in $L^1(T_\Gamma)$ exists and is the same  as the one of  $E^1_{quad,\frac12}(T_\Gamma)$.
\end{dem}

\begin{dem} \em (Theorem \ref{Main4}) \em 

 Just use Proposition \ref{Main41} instead of Proposition \ref{Main31} (in the proof of Theorem \ref{Main3}), and Proposition \ref{Main42} instead of Proposition \ref{Main51} (in the proof of Theorem \ref{Main5}).
\end{dem}

\noindent Let us state and prove now  item $b)$ of Remark \ref{remresults}. We first  introduce $E^1_{BZ\kappa,M,\epsilon}(\Gamma)$ defined by
{\small $$E^1_{BZ\kappa,M,\epsilon}(\Gamma) := \left\{ f\in L^2(\Gamma), \ \sum_{j=0}^\infty \lambda_j a_j \text{ is a molecular $(BZ_\kappa,M,\epsilon)$-representation of $f$ and the sum converges in $L^2(\Gamma)$ } \right\}$$}
and outfitted with the norm
{\small $$\|f\|_{E^1_{BZ\kappa,M,\epsilon}} = \inf\left\{ \sum_{i\in \N}|\lambda_i|, \ \sum_{j=0}^\infty \lambda_j a_j \text{ is a molecular $(BZ_\kappa,M,\epsilon)$-representation of $f$ and the sum converges in $L^2(\Gamma)$ } \right\}.$$}
In the same way, we define $E^1_{BZ2,M+\frac12,\epsilon}(T_\Gamma)$ by
{\small $$E^1_{BZ2,M+\frac12,\epsilon}(T_\Gamma) := \left\{ f\in H^2(T_\Gamma), \ \sum_{j=0}^\infty \lambda_j a_j \text{ is a mol. $(BZ_2,M+\frac12,\epsilon)$-representation of $f$ and the sum converges in $L^2(T_\Gamma)$ } \right\}$$}
and we equipped it with the norm
{\small $$\|f\|_{E^1_{BZ\kappa,M+\frac12,\epsilon}} = \inf\left\{ \sum_{i\in \N}|\lambda_i|, \ \sum_{j=0}^\infty \lambda_j a_j \text{ is a mol. $(BZ_2,M+\frac12,\epsilon)$-representation of $f$ and the sum converges in $L^2(T_\Gamma)$ } \right\}.$$}

\begin{cor} \label{Main30}
Let $\Gamma$ be a weighted graph satisfying \eqref{DV} and \eqref{LB}. 
\begin{enumerate}[(i)]
 \item If $\kappa \in \{1,2\}$, $\epsilon \in (0,+\infty)$ and $M\in \N^*\cap(\frac{d_0}{4},+\infty)$, then 
$$E^1_{BZ\kappa,M,\epsilon}(\Gamma) = H^1_{BZ\kappa,M,\epsilon}(\Gamma) \cap L^2(\Gamma) = E^1_{quad,1}(\Gamma)$$
with equivalent norms. As a consequence, the completion of $E^1_{BZ\kappa,M,\epsilon}(\Gamma)$ in $L^1(\Gamma)$ exists and is equal to $H^1(\Gamma) = H^1_{BZ\kappa,M,\epsilon}(\Gamma)$.

 \item If $\epsilon \in (0,+\infty)$ and $M\in \N\cap(\frac{d_0}{4}-\frac12,+\infty)$, then 
$$E^1_{BZ2,M+\frac12,\epsilon}(T_\Gamma) = H^1_{BZ2,M+\frac12,\epsilon}(T_\Gamma) \cap L^2(\Gamma) = E^1_{quad,\frac12}(T_\Gamma)$$
with equivalent norms. As a consequence, the completion of $E^1_{BZ2,M,\epsilon}(T\Gamma)$ in $L^1(T_\Gamma)$ exists and is equal to $H^1(T_\Gamma) = H^1_{BZ2,M+\frac12,\epsilon}(T_\Gamma)$.

 \item If the Markov kernel $p(x,y)$ satisfies the pointwise gaussian bound \eqref{UE}, then $M$ can be choosen  arbitrarily  in $\N^*$ in (i) and in $\N$ in (ii).
\end{enumerate}
\end{cor}

\begin{dem}
The proof  consists  in noticing,  as the proofs show,  that the $(BZ_\kappa,M,\epsilon)$ (resp. $(BZ_2,M+\frac12,\epsilon)$) representation of $f\in E^1_{quad,1}(\Gamma)$ 
(resp. $F\in E^1_{quad,\frac12}(T_\Gamma)$) constructed in Proposition \ref{Main32} (resp. \ref{Main52})  also converges  in $L^2(\Gamma)$ (resp. $L^2(T_\Gamma)$).

Therefore, we proved in Propositions \ref{Main32} and \ref{Main52} that 
$$E^1_{quad,1}(\Gamma) \subset E^1_{BZ\kappa,M,\epsilon}(\Gamma) \subset H^1_{BZ\kappa,M,\epsilon}(\Gamma) \cap L^2(\Gamma)$$
and
$$E^1_{quad,\frac12}(T_\Gamma) \subset E^1_{BZ2,M+\frac12,\epsilon}(T_\Gamma) \subset H^1_{BZ2,M+\frac12,\epsilon}(T_\Gamma) \cap L^2(T_\Gamma).$$
We end then the proof as in Theorems \ref{Main3}, \ref{Main5} and \ref{Main4}.
\end{dem}

\appendix

\section{A  covering  lemma}

\begin{lem} \label{lemVitali}
 Let $B$ a ball of radius $r \in \N^*$ and $\alpha \geq 1$. There exists a collection of pairwise disjoint balls $(B_i)_{i\in \I_\alpha}$ of radius $r$ such that 
$$ \bigcup_{i\in I_\alpha} B_i \subset \alpha B \subset \bigcup_{i\in I_\alpha}3B_i.$$

\end{lem}

\begin{dem} It is a classical fact and we provide a proof for completeness.  Let $B$ be a ball of radius $r$ and of center $x_0$.
Let $(B_i)_{i\in I_\alpha}$ be a set of disjoint balls included in $\alpha B$ and of radius $r$. 
Assume that $(B_i)_{i\in I_\alpha}$ is maximal, that is, for every ball $B_0$ of radius $r$, either $B_0$ is not included in $\alpha B$, or there exist $i\in I_\alpha$ such that $B_0 \cap B_i \neq \emptyset$. 
Let us prove that 
\begin{equation}
\alpha B \subset \bigcup_{i\in I_\alpha}3B_i.
\end{equation}
Let $x \in \alpha B$ and let us prove that the ball $B(x,2r)$ intersects one of the $B_i$'s. Assume the opposite. 
There exists a path $x_0,x_1,\dots,x_{n-1},x$ joining $x_0$ to $x$ and of length $n = d(x,x_0)<\alpha r$. 
Then the balls $B(x_{\max\{0,n-r\}},r)$ is included in $B(x,2r)$ and in $\alpha B$, that is the set $(B_i)_{i\in I_\alpha}$ is not maximal. 
By contradiction, there exists $i\in I_\alpha$ such that $B(x,2r) \cap B_i \neq \emptyset$, that implies $x\in 3B_i$. 
%
\end{dem}

\begin{cor} \label{corVitali2}
There exist $M\in \N$ and $C>0$ such that for all balls $B$ of radius $r$ and all $j\geq 1$, 
there exists a covering $(B_i)_{i\in I_j}$ of $C_j(B)$ such that
\begin{enumerate}[(i)]
 \item each ball $B_i$ is of radius $r$,
 \item the covering is included in $\tilde C_j : = C_{j-1}(B) \cup C_j(B) \cup C_{j+1}(B)$ (with the convention $C_0(B) = \emptyset$), that is
   $$\bigcup_{i\in I_j} B_i \subset \tilde C_j$$
 \item each point is covered by at most $M$ balls $B_i$.
 \item the number of balls $\#I_j$ is bounded by $C2^{j(d_0+1)}$
\end{enumerate}
\end{cor}

\begin{dem}
Let $B$ be a ball of radius $r$ and $j\geq 1$. Notice that (iv) is a consequence of the three first points. Indeed,
\[\begin{split}
   \#I_j & = \frac{1}{V(2^jB)} \sum_{i\in I_j} V(2^jB) \\
& \leq \frac{1}{V(2^jB)} \sum_{i\in I_j} V(2^{j+3}B_i) \\
& \lesssim \frac{2^{j(d_0+1)}}{V(2^jB)} \sum_{i\in I_j} V(B_i) \\
& \qquad \leq M 2^{j(d_0+1)} \frac{1}{V(2^jB)} V(2^{j+2}B) \\
& \lesssim 2^{j(d_0+1)}.
  \end{split}\]
where the second line is a consequence of (i) and (ii), the third one holds thanks to Proposition \ref{propDV}, 
and the forth one is due to (ii) and (iii).

Let us now prove the first three conclusions of the corollary.

Assume that $r\in \{1,2\}$. Then the collection of balls $(B(x,r))_{x\in C_j(B)}$ satisfies (i), (ii) and (iii). 
Indeed, only (iii) for $r=2$ is not obvious, but is a consequence of the uniform local finiteness of $\Gamma$.

Assume now that $r\geq 3$. Let $s\in \left[\frac{r}{5}, \frac{r}{3}\right] \cap \N$. 
By Lemma \ref{lemVitali} (with $\alpha = 2^{j+1}\frac{s}{r}$), there exists a collection $(\tilde B_i)_{i\in I_\alpha}$ of balls of radius $s$ such that
$$\bigcup_{i\in I_\alpha} \tilde B_i \subset 2^{j+1} B \subset \bigcup_{i\in I_\alpha}3\tilde B_i.$$
We set 
$$I_j = \{i\in I_\alpha, \ 3\tilde B_i \cap C_j(B) \neq \emptyset\}$$
and then $B_i = \frac{r}{s} \tilde B_i$. Let us check that the collection of balls $(B_i)_{i\in I_j}$ satisfies the conclusions of the corollary.
(i) is a consequence of the construction. (ii) is true since 
$$\bigcup_{i\in I_j} B_i \subset \{x\in \Gamma, \ d(x,C_j(B)) < 2s\}.$$
For the point (iii), define for $x\in \Gamma$,
$$I_x = \{i\in I_j, B(x,s) \cap B_i \neq \emptyset \}.$$
Since all $\tilde B_i$ are disjoints, one has then 
$$\sum_{i\in I_x} V(\tilde B_i) \leq V(x,r+s) \leq V(x,6s).$$
However, notice that $B(x,6s) \subset V(12\tilde B_i)$ for all $i\in I_x$. Hence, with the doubling property,
$$V(x,6s) \gtrsim \sum_{i\in I_x} V(12\tilde B_i) \gtrsim \sum_{i\in I_x} V(x,6s)$$
and therefore, $\#I_x \lesssim 1$.
\end{dem}

\section{Exponential decay of some functions}

\begin{lem} \label{ExpDecay} For all $m \in [0,+\infty)$, there exists $C_m,c>0$ such that for all $t\geq 0$ and $k\in \N$, one has 
 $$\left( \frac{1+k}{1+t}\right)^m \left( \frac{t}{1+t}\right)^k \leq C_m \exp\left(-c \frac{k}{1+t} \right).$$
\end{lem}

\begin{dem}
 First check that the function 
$$\varphi(t) \in \R^*_+ \mapsto \left(1 - \frac{1}{1+t} \right)^{1+t}$$
satisfies $0<\varphi(t)<1$ for all $t>0$ and $\ds \lim_{t\to \infty} \varphi(t) = e^{-1}<1$. Then there exists $c>0$ such that $\varphi(t) \in (0,e^{-c})$ for all $t>0$.
From here, one has
\[\begin{split}
   \left( \frac{1+k}{1+t}\right)^m \left( \frac{t}{1+t}\right)^k & \leq \left(1 +  \frac{k}{1+t}\right)^m \left( \frac{t}{1+t}\right)^k \\
& \quad = \left(1 + \frac{k}{1+t}\right)^m \exp\left( \frac{k}{1+t} \ln\varphi(t)\right) \\
& \leq  \left( 1 + \frac{k}{1+t}\right)^m \exp\left( -c\frac{k}{1+t}\right)\\
& \leq C_m \exp\left(-\frac{c}{2} \frac{k}{1+t} \right).
  \end{split}\]
\end{dem}

\bibliographystyle{plain}

\bibliography{../biblio}

\begin{thebibliography}{10}

\bibitem{al}
C.~Arhancet and C.~Le~Merdy.
\newblock Dilation of {R}itt operators on ${L}^p$-spaces.
\newblock 2011.
\newblock Available as http://arxiv.org/ abs/1106.1513.

\bibitem{AMM}
P.~Auscher, A.~McIntosh, and A.~J. Morris.
\newblock Calder\`on reproducing formulas and applications to {H}ardy spaces.
\newblock Available as http://arxiv.org/abs/1304.0168.

\bibitem{AMR}
P.~Auscher, A.~McIntosh, and E.~Russ.
\newblock Hardy spaces of differential forms and {R}iesz transforms on
  {R}iemannian manifolds.
\newblock {\em J. Geom. Anal.}, 18(1):192--248, 2008.

\bibitem{BRuss}
N.~Badr and E.~Russ.
\newblock Interpolation of {S}obolev spaces, {L}ittlewood-{P}aley inequalities
  and {R}iesz transforms on graphs.
\newblock {\em Publ. Mat.}, 53:273--328, 2009.

\bibitem{BZ}
F.~Bernicot and J.~Zhao.
\newblock New abstract {H}ardy spaces.
\newblock {\em J. Funct. Anal.}, 255:1761--1796, 2008.

\bibitem{BD}
T.~A. Bui and X.~T. Duong.
\newblock Hardy spaces associted to the discrete laplacians on graphs and
  boundedness of singular integrals.
\newblock {\em Trans. Amer. Math. Soc.}, 255:1761--1796, 2014.

\bibitem{CMS}
R.~R. Coifman, Y.~Meyer, and E.~M. Stein.
\newblock Some new function spaces and their applications to harmonic analysis.
\newblock {\em J. Funct. Analysis}, 62:304--335, 1985.

\bibitem{CW}
R.~R. Coifman and G.~Weiss.
\newblock Extensions of {H}ardy spaces and their use in analysis.
\newblock {\em Bull. Amer. Math. Soc.}, 83(4):569--645, 1977.

\bibitem{CoulGrigor}
T.~Coulhon and A.~Grigor'yan.
\newblock Random walks on graphs with regular volume growth.
\newblock {\em Geom. Funct. Anal.}, 8(4):656--701, 1998.

\bibitem{CGZ}
T.~Coulhon, A.~Grigor'yan, and F.~Zucca.
\newblock The discrete integral maximum principle and its applications.
\newblock {\em Tohoku Math. J.}, 57:559--587, 2005.

\bibitem{CSC}
T.~Coulhon and L.~Saloff-Coste.
\newblock Puissances d'un op\'erateur r\'egularisant.
\newblock {\em Ann. Inst. H. Poincar\'e Probab. Statist.}, 26(3):419--436,
  1990.

\bibitem{Delmotte1}
T.~Delmotte.
\newblock Parabolic {H}arnack inequality and estimates of {M}arkov chains on
  graphs.
\newblock {\em Revista Matem\`atica Iberoamericana}, 15(1):181--232, 1999.

\bibitem{Dungey}
N.~Dungey.
\newblock A note on time regularity for discrete time heat kernels.
\newblock {\em Semigroups forum}, 72(3):404--410, 2006.

\bibitem{DY2}
X.~T. Duong and L.~X. Yan.
\newblock Duality of {H}ardy and {BMO} spaces associated with operators with
  heat kernel bounds.
\newblock {\em J. Amer. Math. Soc.}, 18(4):943--973, 2005.

\bibitem{DY1}
X.~T. Duong and L.~X. Yan.
\newblock New function spaces of {BMO} type, the {J}ohn-{N}iremberg inequality,
  interpolation, and applications.
\newblock {\em Comm. Pure Appl. Math.}, 58(10):1375--1420, 2005.

\bibitem{FS2}
C.~Fefferman and E.~M. Stein.
\newblock {$H^{p}$} spaces of several variables.
\newblock {\em Acta Math.}, 129(3-4):137--193, 1972.

\bibitem{Fen1}
J.~Feneuil.
\newblock {L}ittlewood-{P}aley functionals on graphs.
\newblock {\em Math. Nachr.}, to appear.
\newblock Available as http://fr.arxiv.org/abs/ 1404.1353.

\bibitem{grigograph}
Alexander Grigoryan.
\newblock Analysis on graphs.
\newblock {\em Lecture Notes, University Bielefeld}, 2009.

\bibitem{HLMMY}
S.~Hofmann, G.~Lu, D.~Mitrea, M.~Mitrea, and L.~Yan.
\newblock Hardy spaces associated to non-negative self-adjoint operators
  satisfying {D}avies-{G}affney estimates.
\newblock {\em Mem. Amer. Math. Soc.}, 214(1007):vi+78, 2011.

\bibitem{HoMar}
S.~Hofmann and J.~M. Martell.
\newblock {$L^p$} bounds for {R}iesz transforms and square roots associated to
  second order elliptic operators.
\newblock {\em Publ. Mat.}, 47(2):497--515, 2003.

\bibitem{HoMay}
S.~Hofmann and S.~Mayboroda.
\newblock Hardy and {BMO} spaces associated to divergence form elliptic
  operators.
\newblock {\em Math. Ann.}, 344:37--116, 2009.

\bibitem{Meyer}
Y.~Meyer.
\newblock {\em Ondelettes et op\'erateurs. {II}}.
\newblock Actualit\'es Math\'ematiques. [Current Mathematical Topics]. Hermann,
  Paris, 1990.
\newblock Op{\'e}rateurs de Calder{\'o}n-Zygmund. [Calder{\'o}n-Zygmund
  operators].

\bibitem{RudinFA}
W~Rudin.
\newblock {\em Functional {A}nalysis}.
\newblock International Series in Pure and Applied Mathematics. McGraw-Hill,
  Inc., New York, second edition, 1991.

\bibitem{Russ}
E.~Russ.
\newblock Riesz tranforms on graphs for $1\leq p\leq 2$.
\newblock {\em Math. Scand.}, 87(1):133--160, 2000.

\bibitem{Russ2}
E.~Russ.
\newblock The atomic decomposition for tent spaces on spaces of homogeneous
  type.
\newblock In {\em C{MA}/{AMSI} {R}esearch {S}ymposium ``{A}symptotic
  {G}eometric {A}nalysis, {H}armonic {A}nalysis, and {R}elated {T}opics''},
  volume~42 of {\em Proc. Centre Math. Appl. Austral. Nat. Univ.}, pages
  125--135. Austral. Nat. Univ., Canberra, 2007.

\bibitem{SW}
E.~M. Stein and G.~Weiss.
\newblock On the theory of harmonic functions of several variables. {I}. {T}he
  theory of {$H^{p}$}-spaces.
\newblock {\em Acta Math.}, 103:25--62, 1960.

\end{thebibliography}

\end{document}